\documentclass[a4paper,10pt]{amsart}
\usepackage{amssymb,amsfonts,amsmath,amsopn,amstext,amsthm,amscd,color,
graphicx,latexsym,mathrsfs,todonotes,url, verbatim,tikz-cd,enumitem}
\usepackage[all]{xy}
\usepackage{geometry}
\newtheorem{theorem}{Theorem}[section]

\theoremstyle{remark}

\newtheorem{example}[theorem]{\bf Example}

\newtheorem{remark}[theorem]{Remark}

\theoremstyle{definition}
 
\newtheorem{definition}[theorem]{Definition}

\newtheorem{assumption}[theorem]{Assumption}

\theoremstyle{plain}

\newtheorem{lemma}[theorem]{Lemma}

\newtheorem{corollary}[theorem]{Corollary}
\newtheorem{proposition}[theorem]{Proposition}

\numberwithin{equation}{section}

\usepackage[commandnameprefix=ifneeded]{changes}

\usepackage{hyperref}

\newcommand{\bbB}{{\mathbb B}}

\newcommand{\bbD}{{\mathbb D}}

\newcommand{\bbG}{{\mathbb G}}

\newcommand{\bbN}{{\mathbb N}}

\newcommand{\bbP}{{\mathbb P}}
\newcommand{\bbQ}{{\mathbb Q}}

\newcommand{\bbT}{{\mathbb T}}

\newcommand{\bbZ}{{\mathbb Z}}

\newcommand{\frd}{{\mathfrak d}}

\newcommand{\frg}{{\mathfrak g}}

\renewcommand{\frm}{{\mathfrak m}}

\newcommand{\frz}{{\mathfrak z}}

\newcommand{\calA}{{\mathcal A}}
\newcommand{\calB}{{\mathcal B}}
\newcommand{\calC}{{\mathcal C}}
\newcommand{\calD}{{\mathcal D}}

\newcommand{\calF}{{\mathcal F}}
\newcommand{\calG}{{\mathcal G}}
\newcommand{\calH}{{\mathcal H}}

\newcommand{\cI}{{\mathcal I}}

\newcommand{\calL}{{\mathcal L}}
\newcommand{\calM}{{\mathcal M}}

\newcommand{\calO}{{\mathcal O}}

\newcommand{\calR}{{\mathcal R}}
\newcommand{\calS}{{\mathcal S}}
\newcommand{\calT}{{\mathcal T}}

\newcommand{\ovg}{{\overline g}}

\newcommand{\ovi}{{\overline i}}

\newcommand{\ovD}{{\overline D}}

\newcommand{\ovG}{{\overline G}}

\newcommand{\ovV}{{\overline V}}

\newcommand{\ovX}{{\overline X}}

\newcommand{\ovZ}{{\overline Z}}

\makeatletter
\newcommand*\bdot{{\mathpalette\bdot@{.9}}}
\newcommand*\bdot@[2]{\mathbin{\vcenter{\hbox{\scalebox{#2}{$\m@th#1\bullet$}}}}}
\makeatother

\newcommand{\an}{{\rm an}}

\newcommand{\coker}{{\rm coker}}

\newcommand{\dR}{{\rm dR}}

\newcommand{\End}{{\rm End}}
\newcommand{\Ext}{{\rm Ext}}

\newcommand{\Fil}{{\rm Fil}}

\newcommand{\Gal}{{\rm Gal}}
\newcommand{\GL}{{\rm GL}}
\newcommand{\gr}{{\rm gr}}
\newcommand{\GS}{{\rm GS}}

 \newcommand{\fib}{{\rm fib}}
\newcommand{\Hom}{{\rm Hom}}
\newcommand{\uHom}{\underline{\Hom}}

\newcommand{\Ind}{{\rm Ind}}

\newcommand{\la}{{\it la}}

\newcommand{\Mod}{{\rm Mod}}

\newcommand{\pr}{{\rm pr}}

\newcommand{\PGL}{{\rm PGL}}

\newcommand{\Rep}{{\rm Rep}}

\newcommand{\sm}{{\rm sm}}

\newcommand{\ST}{{\rm ST}}

\newcommand{\trace}{\mathrm{tr}}

\newcommand{\WD}{\mathrm{WD}}
\newcommand{\cInd}{\mathrm{c\text{-}Ind}}
\newcommand{\Res}{\mathrm{Res}}
\newcommand{\dis}{\mathrm{dis}}
\newcommand{\Cone}{\mathrm{Cone}}
\newcommand{\Tate}{\mathrm{Tate}}
\newcommand{\BZ}{\ensuremath{\mathrm{BZ}}}
\newcommand{\CS}{\ensuremath{\mathrm{CS}}}

\newcommand{\AffRingQp}{{\mathrm{AffRing}^b_{\mathbb{Q}_p}}}
\newcommand{\AffQp}{{\mathrm{Aff}^b_{\mathbb{Q}_p}}}
\newcommand{\Sh}{{\mathrm{Sh}}}

\newcommand{\Map}{\mathrm{Map}}
\newcommand{\Anring}{\mathrm{AnRing}}
\newcommand{\AnialgR}{\mathrm{AniAlg}_{R_\square}}

\DeclareMathOperator{\Lie}{{Lie}}
\DeclareMathOperator{\sgn}{sgn}

\DeclareMathOperator{\Spa}{Spa}

\DeclareMathOperator{\colim}{colim}
\DeclareMathOperator{\AnSpec}{AnSpec}
\begin{document}

\title{Bernstein--Zelevinsky duality for equivariant vector bundles on Drinfeld spaces}

\author{Zhenghui Li}
\address{Institute for Advanced Study, 1 Einstein Dr, Princeton, NJ 08540, U.S.A.}
\email{zhenghuili@ias.edu}
\author{Benchao Su}
\address{IMJ-PRG, Sorbonne Universit\'e, 4 place Jussieu,
75005 Paris, France}
\email{benchao.su@imj-prg.fr}
\author{Zhixiang Wu}
\address{School of Mathematical Sciences, University of Science and Technology of China, 96 Jinzhai Road, 230026 Hefei, China}
\email{zhixiang.wu@ustc.edu.cn}

\begin{abstract}
We study duality and finiteness of locally analytic representations arising from compactly supported coherent cohomology of Drinfeld spaces. We prove that Grothendieck--Serre duality for equivariant vector bundles realizes Bernstein--Zelevinsky duality of locally analytic representations. We follow the strategy of Fargues for $\ell$-adic sheaves and obtain resolutions of locally analytic representations by compactly induced representations. The resolutions motivate us to introduce the notion of locally analytic representations \emph{of compact presentation}, which generalizes the notion of finitely presented representations. These representations are related to prim objects on the classifying stack in the sense of $6$-functor formalism and their extension groups to admissible locally analytic representations are finite-dimensional. As an application, we establish Serre duality for extension groups of locally analytic representations of $\GL_2(\bbQ_p)$ attached to non-trianguline de Rham Galois representations under the $p$-adic local Langlands correspondence.
\end{abstract}
\maketitle
\tableofcontents

\section{Introduction}
Let $p$ be a prime number and let $L$ be a finite extension of $\bbQ_p$. Let $n\geq 2$, and let $\calH$ be the set of all $L$-rational hyperplanes in the analytic projective space $\bbP^{n-1,\an}_L$ over $L$. The Drinfeld space \cite{drinfeld1974elliptic}
\begin{equation}\label{equationDrinfeldspace}
	\Omega_{\rm Dr}:=\bbP^{n-1,\an}_{L}\setminus \cup_{H\in \calH}H
\end{equation}
and its \'etale covers, the Drinfeld towers, have played a crucial role in the Langlands program. These spaces admit natural actions of the $p$-adic Lie group $\GL_n(L)$ or $\PGL_n(L)$. Their cohomology realizes the local Langlands correspondence and the Jacquet--Langlands correspondence in both the $\ell$-adic case when $\ell\neq p$, e.g. \cite{harrisTaylor2001,farguesGenestierLafforgue2008,dat2007theorie}, and in the $p$-adic case when $ L=\bbQ_p$ and $n=2$ \cite{dospinescuLeBras2017revetements, colmezDospinescuNiziol2020}. 

In this paper, we consider compactly supported coherent cohomology of equivariant vector bundles on $\Omega_{\rm Dr}$. The study of coherent cohomology for $\Omega_{\rm Dr}$ was initiated by Morita \cite{moritaMurase1982} and later by Schneider--Teitelbaum \cite{schneiderTeitelbaum2002boundary}, which motivated the development of their theory of locally analytic representations of $p$-adic Lie groups \cite{ST02,ST03}. Let $G=\PGL_n(L)$ and let $\calF$ be a $G$-equivariant vector bundle on $\Omega_{\rm Dr}$. Since the rigid analytic space $\Omega_{\rm Dr}$ is (quasi-)Stein \cite[§1, Prop.~4]{SS91}, the compactly supported coherent cohomology $R\Gamma_c(\Omega_{\rm Dr},\calF)$ concentrates in degree $n-1=\dim \Omega_{\rm Dr}$ by Kiehl's Theorem B \cite{kiehl1967theorem} and Serre duality for Stein spaces \cite{chiarellotto1990duality,vanderPut1992}. The cohomology group
\begin{equation}\label{equationcohomology}
	H_c^{n-1}(\Omega_{\rm Dr},\calF)\simeq H^0(\Omega_{\rm Dr}, \calF^{\vee}\otimes \omega)^{\vee},
\end{equation}
where $\omega$ denotes the canonical line bundle of $\Omega_{\rm Dr}$, is a locally analytic representation of the $p$-adic Lie group $G$. In the case $n=2$, such representations are related to $2$-dimensional de Rham $p$-adic Galois representations by the $p$-adic local Langlands correspondence, according to the Breuil--Strauch conjecture proved by Dospinescu--Le Bras \cite{dospinescuLeBras2017revetements}. We refer to \S~\ref{subsec:Drinfeldspaces} for more examples of equivariant vector bundles on Drinfeld spaces and to \cite{schraen2011gl3,breuilQian2026} for some work on the cohomology groups. 

Our first result is the determination of the Bernstein--Zelevinsky dual of the $G$-representation \eqref{equationcohomology}. Let $E$ be a finite extension of $L$. For $V$ a locally analytic representation of $G$ over $E$, its Bernstein--Zelevinsky dual is given by \cite[Def.~5.4]{strauch2025bernsteinzelevinskydualitylocallyanalytic}
\begin{equation}\label{equationdefinitiondual}
	\bbD_{\BZ}(V):=R\Hom_G(V,\calC_c^{\la}(G,E))
\end{equation}
where $\calC_c^{\la}(G,E)$ denotes the space of compactly supported locally analytic functions on $G$, which is a $G$-bimodule via left and right $G$-translations. The definition extends the classical Bernstein--Zelevinsky duality (also called cohomological duality) for smooth representations of $G$ in \cite[IV.5]{bernstein1992notes} and \cite{Farguesdualite} in which ``locally analytic'' is replaced by ``smooth''. Note that we work in the solid setting of \cite{rodrigues2022solid,JC23} for the derived category of locally analytic representations.
\begin{theorem}[{Corollary~\ref{corBZdualDrinfeld}}]\label{theoremintroBZdual}
    Let $G=\PGL_n(L)$ and let $\calF$ be a $G$-equivariant vector bundle on $\Omega_{\rm Dr}$. Write $\bbD_{\rm GS}(\calF)=\calF^{\vee}\otimes \omega[n-1]$ for the Grothendieck--Serre dual of the vector bundle. Then 
    \[\bbD_{\BZ}(R\Gamma_c(\Omega_{\rm Dr},\calF))\simeq R\Gamma_c(\Omega_{\rm Dr},\bbD_{\rm GS}(\calF))\] 
    as locally analytic representations of $G$.
\end{theorem}
In particular, $\bbD_{\BZ}(H^{n-1}_c(\Omega_{\rm Dr},\calF))\simeq H^{n-1}_c(\Omega_{\rm Dr},\calF^{\vee}\otimes \omega)[1-n]$ and $\bbD_{\BZ}(\bbD_{\BZ}(H^{n-1}_c(\Omega_{\rm Dr},\calF)))\simeq H^{n-1}_c(\Omega_{\rm Dr},\calF)$. The above result is an analogue of the classical result for smooth $G$-representations and $\ell$-adic local systems on $\Omega_{\rm Dr}$ by Fargues \cite{Farguesdualite} (see \cite{mieda2014zelevinsky, hamann2021zelevinsky} for generalizations). There, compactly supported cohomology of the Drinfeld tower intertwines the Bernstein--Zelevinsky duality and the Poincar\'e duality.

The proof of Theorem~\ref{theoremintroBZdual} is an adaptation of the argument of Fargues in \cite{Farguesdualite}. There exists an affinoid open subspace $U\subset \Omega_{\rm Dr}$ which is stabilized by an open compact subgroup of $G$ such that the family of open subspaces $\{g(U)\}_{g\in G}$ forms a locally finite open covering of $\Omega_{\rm Dr}$. Then $ R\Gamma_c(\Omega_{\rm Dr},\calF)$ can be computed using two different \v{C}ech complexes with respect to the affinoid covering (corresponding respectively to $!$-descent and $*$-descent). The terms of the \v{C}ech complexes are finite sums  of $G$-representations of the form $\cInd_{K}^GH^0(V,\calF)$ or $\cInd_{K}^GH^{n-1}_c(V,\calF)$ for some open compact subgroup $K\subset G$ and affinoid open $V\subset \Omega_{\rm Dr}$, where $\cInd_K^G$ denotes the compact induction. The key property of the Bernstein--Zelevinsky dual is \cite[Thm.~5.3]{strauch2025bernsteinzelevinskydualitylocallyanalytic}:
\begin{equation}\label{equationdualitySmithinduction}
\bbD_{\BZ}(\cInd_K^G W)\simeq \cInd_K^G W^{\vee}
\end{equation}
when $W$ is a Smith space so that its $E$-linear dual $W^{\vee}$ is a Banach space. Using \eqref{equationdualitySmithinduction} and Serre duality, one has $\bbD_{\BZ}(\cInd_{K}^GH^{n-1}_c(V,\calF))\simeq \cInd_{K}^GH^{n-1}_c(V,\calF)^{\vee}\simeq \cInd_{K}^GH^{0}(V,\calF^{\vee}\otimes \omega)$. The result then follows from a comparison of the \v{C}ech complex for $R\Gamma_c(\Omega_{\rm Dr},\bbD_{\GS}(\calF))$ and the $\bbD_{\BZ}$-dual of that for $R\Gamma_c(\Omega_{\rm Dr},\calF)$. 
\begin{remark}
	The arguments for Theorem~\ref{theoremintroBZdual} work for more general spaces (Theorem~\ref{maintheorem}). However, to apply \eqref{equationdualitySmithinduction}, we are forced to take $U$ affinoid (instead of Stein). Consequently, we have to work within the framework of the analytic geometry of Clausen--Scholze \cite{scholze2019condensed,scholze2022complex} and Rodr{\'\i}guez Camargo \cite{camargo2024analytic,Juan2024solidgeometry} to use Serre duality for affinoid spaces. Along the way, we prove results on compactly supported cohomology of smooth affinoid spaces (e.g. Proposition~\ref{propcompactsmith}) which may be of independent interest.
\end{remark}

The bonus of the proof of Theorem~\ref{theoremintroBZdual} is infinitely many resolutions of the $G$-representation $H_c^{n-1}(\Omega_{\rm Dr},\calF)$ by representations of the form $\cInd_K^GW$ where $W$ is a Smith space representation of $K$. These $W$ are of the form $H_c^{n-1}(V,\calF)$ where $V\subset \Omega_{\rm Dr}$ is affinoid. The resolutions motivate us to introduce the notion of locally analytic representations \emph{of compact presentation}, which generalizes the notion of finitely presented (smooth) representations. We refer to Definition~\ref{deffinitepresentation} for the precise definition. Roughly speaking, a locally analytic representation of compact presentation is a locally analytic representation that admits a finite resolution by representations of the form $\cInd_K^GW$ where $W$ is a Smith space representation of $K$, and one has certain freedom to change $W$ to obtain other resolutions, corresponding to shrinking or enlarging the affinoid $V\subset \Omega_{\rm Dr}$ in the Drinfeld space cases. By the work of Kohlhaase--Schraen \cite{kohlhaase2012homological}, extended in \cite{strauch2025bernsteinzelevinskydualitylocallyanalytic}, locally analytic principal series representations are of compact presentation.
\begin{proposition}[Proposition~\ref{propositionDrinfeldcompactpresentation}]\label{intropropositionDrinfeldcompactpresentation}
	Let $\calF$ be a $G$-equivariant vector bundle on $\Omega_{\rm Dr}$. Then $H_c^{n-1}(\Omega_{\rm Dr},\calF)$ is a locally analytic representation of $G$ of compact presentation in the sense of Definition~\ref{deffinitepresentation}.
\end{proposition}
The above result is an analogue of other finiteness results of cohomology groups of Drinfeld spaces, e.g. with mod $p$ coefficients \cite[Thm.~0.5]{colmezDospinescuNiziol2023factorisation}. The following theorem is one of the main reasons that we introduce the notion of locally analytic representations of compact presentation.
\begin{theorem}[Theorem~\ref{thm:finite-presentation-perfect}]\label{theoremintrofiniteness}
	Let $E$ be a finite extension of $\bbQ_p$ and let $G$ be a $p$-adic Lie group. Let $V$ be a locally analytic $E$-representation of $G$ of compact presentation in the sense of Definition \ref{deffinitepresentation},
	and let $\Pi$ be an admissible locally analytic representation of $G$ in the sense of Schneider--Teitelbaum \cite{ST03}. Then
	\[\dim_E\Ext^i_{G}(V,\Pi)<\infty\] 
	for every $i\geq 0$ and $\Ext^i_{G}(V,\Pi)= 0$ for all but finitely many $i\geq 0$.
\end{theorem}
Again, the above result generalizes a similar statement for finitely presented and admissible smooth representations. However, the proof for locally analytic representations is much more involved, relying on the theory of compact operators between Banach spaces. More precisely, we use a $p$-adic Schwartz lemma (Proposition~\ref{propfinite}) as in the proof of Cartan--Serre--Kiehl type finiteness theorem \cite{kiehl1967endlichkeitssatz}.
\begin{remark}
	To the best of the authors' knowledge, finiteness results for extension groups of locally analytic representations as Theorem~\ref{theoremintrofiniteness} were previously only available for particular representations before, e.g. \cite{ding2022Ext1} and \cite[Thm.~1.3.3]{breuilQian2026}. 

	The proof of Theorem~\ref{theoremintrofiniteness} may remind experts of Emerton's finiteness (admissibility) result on the Jacquet--Emerton module \cite{emerton2006jacquet} and its derived version \cite{lee2025derivedjacquetemerton} for admissible representations. However, even for $V$ a locally analytic principal series representation, the finiteness of higher extension groups in Theorem~\ref{theoremintrofiniteness} seems to be new, as the second adjointness (as in \cite[Prop.~4.2]{breuil2015vers}) is not available for the derived Jacquet--Emerton module yet.
\end{remark}

\begin{remark}
	The classical Bernstein--Zelevinsky duality for finitely presented smooth representations can be studied from a more categorical point of view, for example from the $6$-functor formalism in \cite{heyer20246functor}. The notions of the Bernstein--Zelevinsky duality and compact presentation are related to ``prim dual'' and ``prim objects'' \cite[\S~4.4]{heyer20246functor} in the derived category of solid locally analytic representations, which should play an important role in the categorical $p$-adic local Langlands correspondence, see \cite[\S~6.2.1]{emerton2025introduction}. Some results of this paper are motivated by this perspective, and have analogues in the categorical setting \cite[Lem.~4.4.9, Lem.~4.4.19]{heyer20246functor}. This work serves as an effort to understand prim locally analytic representations in a representation-theoretic way:
	\begin{theorem}[{Theorem~\ref{theoremprimobject}}]
		Suppose that $V$ is a locally analytic representation of a $p$-adic Lie group. 
		\begin{enumerate}
			\item The representation $V$ is prim if and only if $V$ is compact and its underlying solid $\bbQ_p$-vector space is nuclear. In this case, $\bbD_{\BZ}(V)$ coincides with its prim dual.
			\item If $V$ is of compact presentation (Definition~\ref{deffinitepresentation}), then $V$ is prim.
		\end{enumerate}
	\end{theorem}
	In particular, by Proposition~\ref{intropropositionDrinfeldcompactpresentation}, compactly supported coherent cohomology of Drinfeld spaces gives rise to prim representations. This result was previously obtained by Rodríguez Camargo \cite[Rem.~6.2.56]{camargo2026locallyanalyticnotes} in a different way.
\end{remark}
The last part of the paper focuses on the case $n=2$ and $L=\bbQ_p$. Let $E$ be a finite extension of $\bbQ_p$ and let $\rho:\Gal(\overline{\bbQ_p}/\bbQ_p)\rightarrow \GL_2(E)$ be a de Rham Galois representation of Hodge--Tate weights $0,1$. We assume that the Weil--Deligne representation $\mathrm{WD}(\rho)$ attached to $\rho$ is absolutely irreducible and $\det(\rho)$ is the $p$-adic cyclotomic character $\epsilon$. By \cite{colmez2016representations,dospinescuLeBras2017revetements}, the locally analytic representation $\pi(\rho)$ attached to $\rho$ sits in a short exact sequence
\[0\rightarrow \pi_{\sm}\rightarrow \pi(\rho)\rightarrow \pi_c\rightarrow 0\]
where $\pi_{\sm}$ is the smooth representation of $G=\PGL_2(\bbQ_p)$ attached to $\mathrm{WD}(\rho)$ by the classical local Langlands correspondence and $\pi_c\simeq H^1_c(\Omega_{\rm Dr},\calF)$ for some $G$-equivariant vector bundle $\calF$ on $\Omega_{\rm Dr}$. Thus we may apply Theorem~\ref{theoremintroBZdual} and Theorem~\ref{theoremintrofiniteness} to $\pi_c$ to calculate the Bernstein--Zelevinsky dual of $\pi(\rho)$ and obtain the following result on extension groups.
\begin{theorem}\label{thmintrodualextgroup}
	Let $\pi(\rho)$ be the locally analytic representation of $\PGL_2(\bbQ_p)$ attached to $\rho$. Then one has $\bbD_{\BZ}(\pi(\rho))\simeq \pi(\rho)[-1]$. Moreover, for $0\leq i\leq 3$, there is an isomorphism of finite-dimensional $E$-vector spaces
	\[\Ext^i_{G}(\pi(\rho),\pi(\rho))\simeq \Ext^{3-i}_{G}(\pi(\rho),\pi(\rho))^{\vee}.\]
\end{theorem}
The duality of $\Ext$ groups is established using the Bernstein--Zelevinsky duality together with another duality $\bbD_{\ST}(-)\simeq (-)^{\vee, R\la}$, the Schneider--Teitelbaum duality for admissible locally analytic representations in \cite{ST05} and \cite[\S~4.2]{JC23}. The last duality is a generalization of the smooth duality (contragredient) for smooth representations. Similar to the smooth case (see e.g. \cite{bondal1990representable}, \cite[\S~0.1]{gaitsgoryYomDin2018deligne}, \cite[Rem.~7.54]{zhu2025tame} and \cite{hansen2026Serrefunctor}), the shift of the composition of the two duality functors $\calS(-):=\bbD_{\ST}(\bbD_{\BZ}(-))[\dim G]$ behaves like a Serre functor (Corollary~\ref{cor:ST-BZ-comparison}): we have
\begin{equation}\label{equationintroSerreduality}
	R\Hom_G(V_1,V_2)^{\vee}\simeq R\Hom_G(V_2,\calS(V_1))
\end{equation}
for suitable locally analytic representations $V_1$ and $V_2$. The isomorphism \eqref{equationintroSerreduality} is a formal consequence of various properties of $\bbD_{\BZ}$ and $\bbD_{\ST}$. The representations in Theorem~\ref{thmintrodualextgroup} are simultaneously of compact presentation and admissible, so that \eqref{equationintroSerreduality} can apply. 
\begin{remark}
	In contrast to compact presentation, the admissibility of coherent cohomology of the Drinfeld spaces is not formal and is not known in general. See different proofs for the admissibility of $\pi_c$ in \cite{dospinescuLeBras2017revetements,pan2026completed,ardakovWadsley2023globalsections}. In \S~\ref{subsectionSTdual}, we also determine the Schneider--Teitelbaum dual of $\pi(\rho),\pi_c$, using Colmez's construction of these representations from $(\varphi,\Gamma)$-modules over the Robba ring \cite{colmez2016representations} together with Ding's study of their translations \cite{ding2023change}.
\end{remark}
\begin{remark}
	In a future work, we will use results in the paper to calculate all extension groups appearing in Theorem~\ref{thmintrodualextgroup}, and then determine the derived category generated by the relevant locally analytic representations, in the spirit of Pa{\v{s}}k{\=u}nas' work \cite{pavskunas2013image}. The duality in Theorem~\ref{thmintrodualextgroup} will correspond to the duality of Ext groups for modules over the Galois deformation rings.
\end{remark}

\noindent\textbf{Convention and Notation.} Rigid analytic spaces are viewed as adic spaces. All adic spaces are assumed to be sheafy. For an adic space $X$, we write $|X|$ for the underlying topological space.

If $A$ is an affinoid algebra, write $A^{\circ}$ for the subring of power bounded elements and $A^{\circ\circ}\subset A^{\circ}$ for the ideal of topologically nilpotent elements.

For a partially proper rigid analytic variety, we assume it is separated and countable at infinity. The last condition means that the space $X$ can be written as a countable increasing union $X=\cup_n X_n$ such that each $X_n$ is a quasi-compact open subspace of $X$ and the inclusion $X_n\subseteq X_{n+1}$ factors through the universal compactification of $X_n$ \cite[Def.~5.1.1]{huber1996etale}.

A strict open immersion of rigid analytic spaces $j:U\hookrightarrow X$ is an open immersion such that $j$ factors through the universal compactification of $U$.

We use the language of condensed mathematics \cite{scholze2019condensed}. Throughout the paper, we fix an uncountable solid cutoff
cardinal $\kappa$ \cite[Def.~2.9.11]{mann2022p} and ignore the related set-theoretical issue. We denote the site of  extremally disconnected sets by $\mathrm{ExtDisc}$. For an (animated) analytic ring $A$, we write $D(A)$ for the stable $\infty$-category of $A$-modules and $\underline{A}$ for the underlying condensed ring. If $A$ is static ($\underline{A}$ sits in degree $0$), then we write $\Mod(A)$ for the abelian category of $A$-modules. For a complete Huber pair $(A,A^+)$, we use $(A,A^+)_{\square}$ to denote the associated analytic ring in \cite[\S~3]{andreychev2021pseudocoherent}.  

Let $E$ be a finite extension of $\bbQ_p$, which will be our coefficient field, and let $\calO_E$ be its ring of integers. We write $E_{\square}$ for the analytic ring $(E,\calO_E)_{\square}$. We write $(R)\uHom_{E}(-,-)$ for the internal hom in $\Mod(E_{\square})$ (in $D(E_{\square})$) and $\otimes_{E_{\square}}^{(L)}$ for the (derived) tensor product over $E_{\square}$. If $\calD$ is a solid $E$-algebra, we write $\Mod_{\square}(\calD)$ (resp. $ D_{\square}(\calD)$) for the category of (left) $\calD$-module objects in $\Mod(E_{\square})$ (resp. $D(E_{\square})$) and $\otimes_{\calD}^{(L)\square}$ for the relative tensor products.

Let $\Delta$ denote the simplex category, with objects $[n]=\{0,\ldots,n\}$ and morphisms the order-preserving maps. If $Y\to X$ is a morphism of spaces, we write $(Y^{n/X})_{[n]\in \Delta^{\rm op}}$ or $(Y^{\bullet/X})_{\bullet\in \Delta^{\rm op}}$ for its \v{C}ech nerve, the simplicial object whose $n$th term $Y^{n/X}$ is the $(n+1)$-fold fibre product $Y\times_X \cdots \times_X Y$.

We freely use the results on solid functional analysis in \cite[Appendix A]{bosco2021drinfeldproetale}, \cite[\S~3]{rodrigues2022solid} and \cite[\S~2.2]{colmez2023arithmeticduality}. Throughout, Banach spaces, Fr\'echet spaces and spaces of compact type (over $\bbQ_p$) denote their solid version, which are equivalent to the corresponding locally convex spaces by the results in \textit{loc. cit.}. Especially, strict exact sequences of locally convex spaces of these types correspond to exact sequences of solid $E$-spaces. Recall Smith space over $E$ are linear duals (via $\uHom_E(-,E)$) of Banach spaces and they are compact project generators of $\Mod(E_{\square})$ \cite[Lem.~3.8]{rodrigues2022solid}.

If $H$ is a $p$-adic Lie group, and $E$ a finite extension of $\bbQ_p$, we write $\calC^{\la}(H,E)$ for the space of locally analytic functions on $H$. If $H$ is compact, write $\calD^{\la}(H,E)\simeq \calC^{\la}(H,E)^{\vee}$ for locally analytic distribution algebra. In general take $H_0\subset H$ an open compact subgroup and set $\calD^{\la}(H,E)=E[H]\otimes_{E[H_0]}\calD^{\la}(H_0,E)$. We view $\calD^{\la}(H,E)$ as a solid $E$-algebra. Let $\Rep^{\la}_{E_\square}(H)^{\heartsuit}$ be the abelian category of locally analytic representations of $G$ over solid $E$-spaces and $\Rep^{\la}_{E_\square}(H)$ be its derived category, which is a full subcategory of $D_{\square}(\calD^{\la}(G,E))$ \cite[Def.~3.3.1, Prop.~3.3.6]{JC23}. The embedding $\Rep^{\la}_{E_\square}(H)\hookrightarrow D_{\square}(\calD^{\la}(G,E))$ admits a right adjoint $(-)^{R\la}$, taking the derived locally analytic vectors \cite[Cor.~3.3.7]{JC23}. Admissible locally analytic $H$-representations in the sense of \cite{ST03} are considered as solid locally analytic representations by
\cite[Prop.~4.42]{rodrigues2022solid} and
\cite[Prop.~4.2.7]{JC23}.

If $V$ is a left $H$-module for a group, then we also view $V$ as a right $H$-module via the involution $\iota:h\mapsto h^{-1}$.
\\
\\
\noindent\textbf{Overview.} We first give a review of the basics on the analytic geometry and $6$-functor formalism in \S~\ref{sectionanalytic}, then we use the machinery to prove the main duality theorem (Theorem~\ref{maintheorem}) in \S~\ref{sectionequivariantvectorbundle}. In \S~\ref{sectionfintiepresentation}, we define and study locally analytic representations of compact presentation. In \S~\ref{sec:duality-supercuspidal}, we study supercuspidal locally analytic representations of $\GL_2(\bbQ_p)$. In Appendices~\ref{sec:STdual} and \ref{sec:translation-duality}, we give a review on the Schneider--Teitelbaum duality and translation functors.
\\
\\
\noindent\textbf{Declaration of AI use.}
We began to use AI in the later stages of the project, for literature search, proof exploration, {\TeX} assistance, and proofreading. The main mathematical contribution of AI is an argument leading to the formulation and the proof of Lemma~\ref{lem:strict-filtered-compact}, a key homological lemma used in the proof of Theorem~\ref{thm:finite-presentation-perfect}. The argument of AI is more direct than our original strategy; see Remark~\ref{remarktworesolution}. AI also contributed to the proof of Proposition~\ref{prop:DlaGRla}, found the reference for Proposition~\ref{propositionalternative}, and suggested further corrections and improvements.
\\
\\
\noindent\textbf{Acknowledgements.} The third author would like to thank Matthias Strauch for discussions and the collaboration on \cite{strauch2025bernsteinzelevinskydualitylocallyanalytic} that led to this project. He would also like to thank Juan Esteban Rodríguez Camargo for very helpful discussions. The second author thanks the Tianyuan Mathematics Research Center for the opportunity to present some of the results of this paper in July 2026. The first author would like to acknowledge Simons Foundation and Institute for Advanced Study for their support, and thank Peking University for hospitality during his visit.

Part of the paper was written when the first author was in his Ph.D program at Sorbonne Universit\'e. This program has received funding from the European Union's Horizon 2020 research and innovation programme under the Marie Skłodowska-Curie grant agreement No 945332 and partially funded by the project “Group schemes, root systems, and related representations” funded by the European Union - NextGenerationEU through Romania’s National Recovery and Resilience Plan (PNRR) call no. PNRR-III-C9-2023- I8, Project CF159/31.07.2023. The second author was partially supported by the National Natural Science Foundation of China under Grant Nos.~12231001 and~12321001. This work emerged while the third author was a postdoc at the University of Münster, funded by the Deutsche Forschungsgemeinschaft (DFG, German Research Foundation) – Project-ID 427320536 – SFB 1442, as well as under Germany's Excellence Strategy EXC 2044 390685587, Mathematics Münster: Dynamics–Geometry–Structure.

\section{Analytic geometry and Serre duality}\label{sectionanalytic}
In this section, we give a brief review of the rigid analytic geometry and the six functor formalism. As we will use the Serre duality for affinoid spaces, we will work in the framework of the analytic geometry developed by Clausen--Scholze in \cite{scholze2022complex}.

\subsection{Analytic adic space as categorified locale}

In the next two subsections we recall the Serre duality on smooth rigid analytic varieties proved in \cite{camargo2024analytic} and \cite{Juan2024solidgeometry}. 

For a presentably symmetric monoidal stable $\infty$-category $C$ (see \cite[\S~1]{camargo2024analytic} for the notion), let $\calS(C)^{\rm op}$ be the poset of isomorphic classes of idempotent algebras in $C$ such that $A\leq A'$ if there is an arrow $A'\to A$ \cite[Cons.~5.2]{scholze2022complex}. By \cite[Prop.~5.3]{scholze2022complex}, the poset $\calS(C):=(\calS(C)^{\rm op})^{\rm op}$ is a frame (i.e. a category that is similar to the category of open subsets of a topology space) and can be thought as open complements of $\calS(C)^{\rm op}$. For $Z\in \calS(C)^{\rm op}$ corresponding to an idempotent algebra $A$, we denote $C(Z):=\Mod_A(C)=\{X\in C\mid X\simeq X\otimes A\}\subset C$, the complement open of $Z$ as $Z^c\in \calS(C)$ and $C(Z^c):=C/C(Z)$, the Verdier quotient. We have the following localization functors:

\begin{definition}[{\cite[Def.~2.2.2]{camargo2024analytic}}]\label{openclosed}
	Let $U=Z^c\in \calS(C)$ with associated idempotent algebra $A$, we have the following functors:
	\begin{enumerate}
		\item Upper star functors: $i_Z^*: C\to C(Z): M\mapsto M\otimes A$, and $j_U^*:C\to C(U)$ the natural projection.
		\item Lower star functors: $i_{Z*}:C(Z)\to C$ the forgetful functor, and $j_{U*}:C(U)\to C$ such that $j_{U*}j^*_UM=\underline{\Hom}_C([1\to A],M)$ for all $M\in C$. Here $[1\rightarrow A]$ denotes the fiber of the map $1\rightarrow A$.
		\item Upper shriek functors: $i^!_Z:C\to C(Z),M\mapsto \underline{\Hom}_C(A,M)$ and $j_U^!:=j_U^*$.
		\item Lower shriek functors:  $i_{Z!}:=i_{Z*}$ and $j_{U!}:C(U)\to C$ such that $j_{U!}j_U^*M=[1\to A]\otimes M$.
	\end{enumerate}	
	A symmetric monoidal colimit preserving functor $f^*:C\to D$ of presentably closed symmetric monoidal stable $\infty$-categories is called an open (resp. closed) immersion if it is equivalent to a functor of the form $j_U^*:C\to C(U)$ (resp. $i_Z^*:C\to C(Z)$).
\end{definition}
The above functors satisfy the usual adjunction properties. For $M\in C$ and $Z\in \calS(C), U=Z^c$, we get the following localization sequence for free:
\begin{equation}\label{localizationsequence}
	\begin{aligned}
		j_{U!}j_{U}^*M\to &M\to i_{Z*}i_Z^*M\to,\\
		i_{Z*}i_Z^!M\to &M\to j_{U*}j_U^*M\to.
	\end{aligned}
\end{equation}

\begin{definition}[{\cite[Def.~7.1]{scholze2022complex}}]
	A categorified locale is a triple $(X,C,f)$ consisting of a locale $X$, a presentably symmetric monoidal stable $\infty$-category $C$ and a morphism of locale $f:\calS(C)\to X$. A morphism of categorified locales $F:(X,C,f)\to (Y,D,g)$ consists of a morphism of locale $F:X\to Y$ and a colimit preserving exact symmetric monoidal functor $G:D\to C$ commuting with $f$ and $g$. It can be enriched into an $\infty$-category (See \cite[Def.~4.2]{aoki2025sheavesspectrumadjunction}).
\end{definition}

Let $(A,A^+)$ be a sheafy Tate Huber pair and $X=\Spa(A,A^+)^{\rm op}$ the locale of open subspaces of $\Spa(A,A^+)$. By \cite[Prop.~2.3.2]{camargo2024analytic}, we can form a categorified locale $(X,D((A,A^+)_\square),f)$ as following: For each $U\in X$ a rational localization of $\Spa(A,A^+)$, the base change functor $j_U^*: D((A,A^+)_\square)\to D((\calO(U),\calO^+(U))_\square)$ is an open localization. So it corresponds to an idempotent algebra $A_Z\in D((A,A^+)_\square)$. The map $f$ is given by $U\mapsto Z^c\in \calS(D((A,A^+)_\square))$. This can be globalized to a fully faithful functor $X\mapsto X^{\it lo}=(|X|,D(X),f)$ from the 1-category of analytic adic spaces over $\bbQ_p$ to the $\infty$-category of categorified locales \cite[Cor.~2.3.4]{camargo2024analytic}: for an analytic adic space $X$ with $U_{\bullet} $ an hypercover of $X$ by open affinoid subspaces, then 
\begin{equation}\label{equationformulaDX}
	D(X)\simeq\varprojlim_{[n]\in \Delta} D(U_{n})
\end{equation}
which is independent of the choice of the hypercover.
\begin{lemma}\label{lemmafactor} 
Assume $X$ is an analytic adic space, $U$ is an open subset of $X$ and $Z$ is a closed subset of $X$.
	\begin{enumerate}
	    \item If $Z\subseteq U$ then there is a natural map $j_{U*}1_{U}\to A_Z$ in $D(X)$. 
		\item If $U\subseteq Z$ then there is a natural map $A_Z\to j_{U*}1_{U}$ in $D(X)$.
	\end{enumerate}
\end{lemma}
\begin{proof}
	Let $A_Z,A_{U^c}$ be the idempotent algebra in $\calS(D(X))^{\rm op}$ corresponding to the closed subspaces $Z$ and $U^c$.

	(1) We need to give a natural map 
	$$\uHom([1\to A_{U^c}],1)\to A_Z.$$
	Since $Z\subseteq U$, we have $A_Z\otimes A_{U^c}=0$ \cite[Prop.~5.3(1)]{scholze2022complex} and we get natural maps
	$$\uHom([1\to A_{U^c}],1)\to \uHom([1\to A_{U^c}],A_Z)\simeq \uHom_{A_Z}([1\to A_{U^c}]\otimes A_Z,A_Z)\simeq A_Z,$$
	where the first map is induced by the map $1\rightarrow A_Z$ corresponding to $Z\subset U$ and for the last isomorphism we used the assumption that $\otimes$ is colimit preserving and exact.

	(2) We need to give a natural map
	$$A_Z\to \uHom([1\to A_{U^c}],1).$$
	By adjunction, it suffices to give a natural map $[1\to A_{U^c}]\otimes A_Z\to 1$. Since $U\subseteq Z$, this implies $[1\to A_{U^c}]\otimes [1\to A_Z]=0$ by \cite[Prop.~5.3(3)]{scholze2022complex} which means $[1\to A_{U^c}]\otimes A_Z\simeq [1\to A_{U^c}]$. Recall that $[1\rightarrow A_{U^c}]=j_{U!}j_{U}^*1$. By adjunction there is a natural map $[1\rightarrow A_{U^c}]\to 1$. 
	\end{proof}

\begin{remark}\label{remarkmapsfactorizationclosedopen}
	Suppose that we have inclusions
	$U_1\subseteq Z\subseteq U_2\subset X$ with open immersions $j_{U_1}:U_1\to X$ and $j_{U_2}: U_2\to X$ and $j_{12}:U_1\to U_2$. A diagram chasing shows that the natural map 
	$j_{U_2*}1_{U_2}\to j_{U_1*}1_{U_1}$, coming from the natural map $1_{U_2}\to j_{12*}j^*_{12}1_{U_2}$ and $j_{U_1*}1_{U_1}=j_{U_2*}j_{12*}j^*_{12}1_{U_2}$, factors as $j_{U_2*}1_{U_2}\to A_Z\to  j_{U_1*}1_{U_1}$.
\end{remark}

The categorified locale associated to an analytic adic space is a derived Tate adic space in \cite[\S~2]{camargo2024analytic} which we recall below. We first recall the definition of bounded affinoid analytic rings. Let $\mathrm{AnRing}_{Z_{\square}}$ be the $\infty$-category of analytic rings over $\bbZ_{\square}$, see \cite[\S~1]{camargo2024analytic}.
\begin{definition}[{\cite[Def.~2.6.1, Def.~2.6.6, Def.~2.6.10]{camargo2024analytic} or \cite[Def.~7.1.1]{Juan2024solidgeometry}}]\label{definitionaffinoidanalyticring}
Let $A=(\underline{A},D(A))\in \mathrm{AnRing}_{Z_{\square}}$. Let $R_\square:=(\bbZ((\pi)),\bbZ[[\pi]])_\square$ be the analytic ring associated to the Huber pair and let $\AnialgR$ be the $\infty$-category of animated $R_{\square}$-algebras.
\begin{enumerate}
	\item The subring of solid elements is a discrete animated ring defined by
	$$A^+:=\Map_{\Anring_{\bbZ_\square}}(\bbZ[T]_\square,A).$$
	\item  The analytic ring $A$ is called an affinoid analytic ring if the induced map $(\underline{A},\pi_0(A^+))_\square\to A$ is an equivalence of analytic rings.
 \item Assume further that $A$ is over $R_\square$, the condensed subring of power bounded elements is defined as for $S\in \mathrm{ExtDisc}$,
 $$A^\circ(S):=\Map_{\AnialgR}(R_\square\langle\bbN[S]\rangle,A)$$
 and the condensed subring of bounded elements is defined as $A^b:=(A^\circ)[1/\pi]$.
 \item Assume $A$ is over $R_\square$. Then $A$ is called bounded if  the natural map $A^b\to \underline{A}$ is an equivalence.
\end{enumerate}
	The category of bounded affinoid analytic rings over $\bbQ_p$ is denoted by $\AffRingQp$. Let $\AffQp$ be the opposite category of $\AffRingQp$. For $A\in \AffRingQp$, the corresponding object in $\AffQp$ is denoted by $\AnSpec A$ and is called the analytic spectrun of $A$.
\end{definition}

\begin{remark}
 We will only consider analytic rings over $\bbQ_p$. This avoids the discussion between $E_{\infty}$-rings and commutative animated rings (See \cite[Rem.~2.1.4]{camargo2024analytic}). 
  \end{remark}
  
Let $A\in \AffRingQp$, by \cite[Thm.~2.7.2, Def.~2.7.3]{camargo2024analytic}, there is a categorified locale 
\[(|\Spa A|,D(A),\rho_A)\] 
associated to it. We note that if $A=(R,R^+)_\square$ is associated to a Tate-Huber pair, then $|\Spa A|\cong |\Spa(R,R^+)|$ (\cite[Prop.~2.7.8(3)]{camargo2024analytic}). The category $\AffRingQp$ is endowed with the analytic topology given by \cite[Def.~2.7.20]{camargo2024analytic}, which is subcanonical (\cite[Lem.~2.7.21]{camargo2024analytic}). 

Let $\Sh_{\an}(\AffQp)$ be the category of sheaves of anima on the site $\AffQp=\mathrm{AffRing}_{\bbQ_p}^{b,\rm op}$ with respect to the analytic topology.
 \begin{definition}[{\cite[Def.~2.7.22]{camargo2024analytic}}]
 	A derived Tate adic space (over $\bbQ_p$) is a sheaf $X\in \Sh_{\an}(\AffQp)$ which admits an open analytic cover by representable sheaves.
 \end{definition}
 
 \begin{example}\label{examtatespace}
Let $X$ be an analytic adic space over $\bbQ_p$, we can and will (implicitly) view it as a derived Tate adic space as the following: Let $\mathrm{Affadic}_{\bbQ_p}$ be the site of analytic affinoid adic spaces over $\bbQ_p$ endowed with analytic topology. By \cite[Prop.~2.7.8]{camargo2024analytic} and \cite[Prop.~4.12]{andreychev2021pseudocoherent}, the functor 
\begin{equation}\label{equationfunctoraffinoidtobounded}
	\mathrm{Affadic}_{\bbQ_p}\to \AffQp; \Spa(A,A^+)\mapsto \AnSpec ((A,A^+)_\square)	
\end{equation}
is a continuous functor between sites. Thus it extends to a colimit preserving functor $(-)_{\Tate}:\Sh_{\mathrm{\an}}(\mathrm{Affadic}_{\bbQ_p})\to \Sh_{\an}(\AffQp)$ (i.e. left Kan extension). We view $X$ as an object in $\Sh_{\mathrm{an}}(\mathrm{Affadic}_{\bbQ_p})$ and define the derived Tate adic space associated to it to be $X_\Tate$.
 
 By \cite[Def.~2.7.22]{camargo2024analytic}, one can associate a categorified locale $(X_{\Tate})^{lo}=(|X_{\Tate}|,D(X_{\Tate}))$ to $X_{\Tate}$ and by  \cite[Cor.~2.7.23]{camargo2024analytic} we have an isomorphism of categorified locales $(X_{\Tate})^{lo}\cong X^{lo}$, where $D(X)$ is given in \eqref{equationformulaDX}.
 \end{example}

 \begin{lemma}\label{lemmacompute}
	Let $X$ be an analytic adic space. Assume $Z\subseteq X$ is a closed subset such that $Z=\cap_i U_i$ where $\{U_i\}$ is a cofiltered family of open subsets of $X$ and for each $i$ there is $i'>i$ with $Z_i$ closed such that $U_{i'}\subseteq Z_i\subseteq U_i$. Then $A_Z=\colim_i j_{i*}1_{U_i}$. Here $A_Z$ is the idempotent algebra in $D(X)$ corresponding to $Z$.
\end{lemma}

\begin{proof}
	 By \cite[Prop.~5.3(4)]{scholze2022complex}, $A_Z=\colim_i A_{Z_i}$ and the condition says that the system $\{A_{Z_i}\}_i$ and $j_{i*}1_{U_i}$ are cofinal to each other (see Lemma \ref{lemmafactor}). The conclusion follows.
\end{proof}
\subsection{Six functor formalism and descent for quasi-coherent sheaves} 
To have Serre duality, we need a six functor formalism of quasi-coherent sheaves on analytic adic spaces. We first recall a general setup of six functor formalism in \cite[\S~A.5]{mann2022p} (see also a summary in \cite[\S~3.1]{camargo2024analytic}).

A geometric setup \cite[Def.~A.5.1]{mann2022p} is a pair $(\calC,E)$ consisting of an $\infty$-category $\calC$ with finite limits and $E$ a collection of homotopy classes of morphisms in $\calC$ containing all isomorphisms and stable under pullbacks and compositions. Let $\mathrm{Corr}(\calC,E)$ be the $\infty$-category of correspondences associated to $(\calC,E)$ \cite[Def.~A.5.2]{mann2022p}. Objects of $\mathrm{Corr}(\calC,E)$ are those of $\calC$, and a morphism from $X$ to $Y$ is a correspondence
\begin{equation}\label{equationcorrespondence}
	\begin{tikzcd}
		& V\ar[ld,"f"']\ar[rd,"g"]\\
		X && Y
	\end{tikzcd}		
\end{equation}
in $\calC$ with $g\in E$. 

A three-functor formalism \cite[Def.~A.5.6]{mann2022p} on $(\calC,E)$ is a lax symmetric monoidal functor 
\[D:\mathrm{Corr}(\calC,E)\to \mathrm{Cat}_\infty,\]
where $\mathrm{Cat}_\infty$ is the $\infty$-category of $\infty$-categories with the product symmetric monoidal structure. Taking $f$ or $g$ to be the identity in \eqref{equationcorrespondence} one obtains functors $ g_!:D(X)\rightarrow D(Y)$ and $f^*:D(X)\rightarrow D(Y)$ respectively. A six-functor formalism \cite[Def.~A.5.7]{mann2022p} is a three-functor formalism for which $f^*$, $f_!$ and $\otimes$ admit right adjoints, denoted by $f_*$, $f^!$ and $\underline{\Hom}$ respectively.

To construct three or six functor formalisms, \cite[Def.~A.5.9]{mann2022p} introduced the notion of a suitable decomposition of $E$, which consists of a pair $I,P\subset E$ satisfying analogue properties of the classes of \'etale morphisms and proper morphisms, and every $f\in E$ is of the form $f=p\circ i$ for some $i\in I$ and $p\in P$. Let $\calC=\AffQp=\mathrm{AffRing}_{\bbQ_p}^{b,\mathrm{op}}$. 
\begin{definition}[{\cite[Def.~3.2.2]{camargo2024analytic}}]\label{definitionsuitabledecomposition}
	We have the following classes of morphisms in $\calC$:
	\begin{enumerate}
		\item Let $I$ be the class of morphisms $f:\AnSpec B\to \AnSpec A$ such that $f^*:D(A)\to D(B)$ is an open localization in the sense of Definition~\ref{openclosed};
		\item Let $P$ be the class of morphisms $f:\AnSpec B\to \AnSpec A$ such that the analytic structure of $B$ is induced from $A$ (i.e.\ $D(B)=\mathrm{Mod}_{\underline{B}}(D(A))$, the $\underline{B}$-module objects in $D(A)$);
		\item Let $E$ be the class of morphisms of the form $p\circ j$ with $j\in I$ and $p\in P$.
	\end{enumerate}
\end{definition}
By \cite[Lem.~3.2.5]{camargo2024analytic}, $(I,P)$ gives a suitable decomposition of $E$. The functor $D:\calC^{\rm op}\to\mathcal{P}r^L:\AnSpec A\mapsto D(A)$ enhances to a six-functor formalism
\begin{align}\label{equationsixfunctorformalismaffinoid}
	D:\mathrm{Corr}(\calC,E)\to\mathcal{P}r^L.
\end{align}
applying \cite[Prop.~A.5.10]{mann2022p}. Here $\mathcal{P}r^L$ denotes the $\infty$-category of presentable stable $\infty$-categories  and colimit-preserving functors.
\begin{example}\label{examplecanonicalcompactifications}
	Let $f:(A,A^+)\rightarrow  (B,B^+)$ be a morphism of analytic Huber pairs. 
	\begin{enumerate}
		\item If $f$ is an open immersion, then $f$ lies in the class $I$ by \cite[Prop.~2.3.2]{camargo2024analytic}.
		\item The map $f$ factors through $f':(A,A^+)\rightarrow (B,B^+_A):=(B,\widetilde{(f(A^+)+B^{\circ\circ})})$, where $\widetilde{(-)}$ denotes the normalization. By \cite[Prop.~3.32]{andreychev2021pseudocoherent}, one has $(B,B^+_A)_{\square}\simeq (B,f(A^+))_{\square}$ whose analytic ring structure is induced from $(A,A^+)_{\square}$ via $A\rightarrow B$. Hence the map of affinoid analytic spaces induced by $f'$ lies in the class $P$.
	\end{enumerate}
\end{example}
\begin{remark}\label{remarkuniversalcompletion}
	The map in the class $P$ in the above example is related to Huber's universal compactification of adic spaces \cite[Def.~5.1.1]{huber1996etale}. Let $f:X\rightarrow Y$ be a morphism between analytic adic spaces. The universal compactification $\ovX^{/Y}$ of $f$ satisfies that for any maps $X\to X'\xrightarrow{f'} Y$ with $f'$ partially proper, the map $X\to X'$ factors through $\ovX^{/Y}$, which exists if $f$ is a separated and taut morphisms between rigid analytic spaces \cite[Thm.~5.1.5]{huber1996etale}. If $X=\Spa(B,B^+)\rightarrow Y=\Spa(A,A^+)$ is a map between affinoid rigid spaces, then $\ovX^{/Y}=\Spa(B,B^+_A)$ where the latter is as in Example~\ref{examplecanonicalcompactifications}.
\end{remark}	

There is a $D$-topology on $\calC$ defined as the following: 
\begin{definition}[{\cite[Def.~5.13]{scholze2025six}}] \label{defDtopo}
Assume $(\calC,E,D)$ is a six functor formalism. Let $f: X\to Y$ be a morphism in $\calC$. Let $\{X^{\bullet/ Y}\}_{\bullet\in \Delta}$ be the \v Cech nerve of the morphism $f$. For any morphism $Y'\to Y$ in $\calC$, let $f': X':=X\times_Y Y'\to Y'$ be the pullback of $f$.  
	\begin{enumerate}
		\item The map $f$ forms a canonical cover if for any $Z\in \calC$ and any $Y'\to Y$ in $\calC$, the functor $\Hom_{\calC}(-,Z)$ satisfies descent along $f'$, i.e. the following map is an equivalence
		$$\Hom_{\calC}(Y',Z)\stackrel{\simeq}\longrightarrow \lim_{\bullet\in\Delta} \Hom_{\calC}(X'^{\bullet/ Y'},Z).$$
		\item The map $f$ satisfies universal $*$-descent if for any $Y'\to Y$ in $\calC$, the functor $D^*$ satisfies descent along pullback morphism $f'$, i.e. the following map is an equivalence of categories where transition maps are given by $*$-pullback
		$$D(Y')\stackrel{\simeq}\longrightarrow \lim_{\bullet\in\Delta} D^*(X'^{\bullet/Y'}).$$
		\item Assume $f\in E$. The map $f$ satisfies universal $!$-descent if for any $Y'\to Y$ in $\calC$, the functor $D^!$ satisfies descent along pullback morphism $f'$, i.e. the following map is an equivalence of categories where transition maps are given by $!$-pullback
		$$D(Y')\stackrel{\simeq}\longrightarrow \lim_{\bullet\in\Delta} D^!(X'^{\bullet/Y'}).$$
		\item The map $f$ is a $D$-cover if it is a canonical cover, satisfying universal $*$-descent and universal $!$-descent. The $D$-topology on $\calC$ is generated by $D$-covers.
	\end{enumerate}
\end{definition}

\begin{remark}\label{remarkDtopo}
	Assume $f:X\to Y$ lies in $E$ and satisfies $!$-descent. Let $\calF\in D(Y)$, we have an equivalence in $D(Y)$
	$$\colim_{[n]\in \Delta^{\mathrm{op}}} f_{n!}f_n^!\calF\stackrel{\sim}\longrightarrow \calF$$ 
	 where $f_n: X^{n/Y}\to Y$. This is because for $\calG\in D(Y)$ we have
	$$\Hom_{D(Y)}(\colim_{[n] \in \Delta^{\mathrm{op}}} f_{n!}f_n^!\calF,\calG)\simeq \lim_{[n] \in \Delta} \Hom_{D(X^{n/Y})}(f^!_n\calF,f^!_n\calG)\simeq \Hom_{D(Y)}(\calF,\calG).$$
	Let $f:X\to Y$ be a morphism in $\calC$ and satisfies $*$-desecent. Assume $D(Y)$ is complete and $\calF\in D(Y)$, the same argument as above shows that we have an equivalence in $D(Y)$:
	$$\calF\stackrel{\sim}\longrightarrow \lim_{[n] \in \Delta} f_{n*}f_n^*\calF.$$
\end{remark}
From now on, the $*$- and $!$-descent will refer to the above definition and remark.
\begin{definition}[{\cite[Def.~3.2.10]{camargo2024analytic}}]\label{definitionTatestack}
	The $\infty$-category of Tate stacks over $\bbQ_p$ is the category $\Sh_D(\AffQp)$ of sheaves in anima of $\AffQp$ for the $D$-topology.
\end{definition}

By \cite[Theorem 3.4.11]{heyer20246functor}, we can extend the six functor formalism \eqref{equationsixfunctorformalismaffinoid} to $(\Sh_{\an}(\calC),E')=(\Sh_{\an}(\AffQp),E')$ and $(\Sh_{D}(\calC),E'')$ where $E'$ (resp. $E''$) is a collection of morphisms in $\Sh_{\an}(\calC)$ (resp. $\Sh_{D}(\calC)$) containing $E$, local on the target and source (cf. \cite[Cor.~3.2.6]{camargo2024analytic}). We call morphisms in $E'$ (resp. $E''$) $!$-able morphisms.

\begin{example}\label{exampleopencoverDcover}
Let $X$ be an analytic adic space over $\bbQ_p$, we can view it as a Tate stack as the following (See Example \ref{examtatespace} and \cite[\S~1.3]{mikami2026finiteness}): Let $\mathrm{Affadic}_{\bbQ_p}$ be the site of analytic affinoid adic spaces over $\bbQ_p$ endowed with analytic topology. By \cite[Prop.~4.12]{andreychev2021pseudocoherent} and \cite[Lem.~3.2.9]{camargo2024analytic}, the functor \eqref{equationfunctoraffinoidtobounded} is also a continuous functor where $ \AffQp$ is now endowed with $D$-topology. Thus it extends to a colimit preserving functor $(-)_\square: \Sh_{\mathrm{an}}(\mathrm{Affadic}_{\bbQ_p})\to \Sh_D(\AffQp)$. We view $X$ as an object in $\Sh_{\mathrm{an}}(\mathrm{Affadic}_{\bbQ_p})$ and define the Tate stack associated to it to be $X_\square\in \Sh_D(\AffQp)$. 
 
If $\{j_i:U_i\to X\}_{i\in I}$ is an analytic cover of $X$, then $\{j_{i\square}: U_{i\square}\to X_\square\}_{i\in I}$ is a $D$-cover: Put $U=\sqcup_i U_i$. We note that by \cite[Exam.~1.55]{mikami2026finiteness}, we have $(U_i\times_X U_j)_\square\simeq U_{i\square}\times_{X_\square} U_{j\square}$, and similarly for higher intersections. Thus, by construction and colimit preservation, $X_{\square}\simeq\colim_{[n]\in \Delta^{\rm op}}U_{\square}^{n/X_\square}$, so $\{j_{i\square}: U_{i\square}\to X_\square\}_{i\in I}$ is a canonical cover. By \cite[Rem.~1.15]{mikami2026finiteness} (cf.\ Lem.~1.56 there), it satisfies universal $D^*$-descent. In particular, as \eqref{equationformulaDX},
\begin{equation*}
	D(X)=\varprojlim_{[n]\in \Delta}D^{*}(U^{n/X}).
\end{equation*}
By \cite[Lem.~1.57]{mikami2026finiteness}, the map $j_{i\square}$ is $D$-\'etale \cite[Def.~4.6.1]{heyer20246functor}. Hence by \cite[Lem.~4.6.4(i)]{heyer20246functor} we have $j_{i\square}^!=j_{i\square}^*$ (see also \cite[\S~4.7]{meyerWagner2024qHodge}). So $\{j_{i\square}: U_{i\square}\to X_\square\}_{i\in I}$ also satisfies universal $D^!$-descent by \cite[Lem.~4.7.1]{heyer20246functor}.
\end{example}
 
\begin{remark}
By $!$-descent in the above example, cohomology of an analytic adic space does not change whenever we view analytic adic spaces as analytic sheaves or $D$-sheaves, or a solid stack \cite[Def.~3.2.7]{camargo2024analytic}. We choose to use analytic topology because it is not clear if a derived Tate adic space is a Tate stack and we need the locale in our proof of comparison of two compactly supported cohomologies in \S~\ref{subseccompare}, where the topology spaces matter. Note that, a priori, the closed subsets $Z\subset X$ in Lemma~\ref{lemmacompute} or Lemma~\ref{lemmafactor} need not be adic spaces nor analytic stacks.
\end{remark}
  
All discussions in \cite[\S~3.6]{camargo2024analytic} are on analytic rings, thus apply in our setting and we can deduce the Serre duality. Let $E$ be a finite extension of $\bbQ_p$. Let $f:X\rightarrow *=\Spa(E,\calO_E)$ be a smooth rigid analytic spaces over $E$. Then the above six functor formalism supplies a functor $f_!:D(X)\rightarrow D(E_{\square})$. Let $\omega_X$ be the dualizing complex of $X$. For $\calF\in D(X)$, we write 
\[R\Gamma_c(X,\calF):=f_!\calF\] 
and $R\mathscr{H}om(\calF,\omega_X)$ for the internal hom in $D(X)$, where $\omega_X=f^!E$.
\begin{theorem}[{\cite[Thm.~3.6.15]{camargo2024analytic}}]\label{theoremserreduality}
  	Assume $X$ is a smooth rigid analytic variety over $E$ and $\calF\in D(X)$ is a quasi-coherent sheaf on $X$, then we have a quasi-isomorphism
  	\[R\uHom_E\left(R\Gamma_c(X,\calF),E\right)\simeq R\Gamma(X,R\mathscr{H}om(\calF,\omega_X))\]
	in $D(E_\square)$.
  \end{theorem}  
  
\begin{remark}\label{remarkgroupaction}
	We remark on the group actions on (compactly supported) cohomologies and duality. Assume $G$ is a $p$-adic Lie group acting continuously on a smooth rigid analytic variety $X$ and $\calF$ is a $G$-equivariant vector bundle (see \S~\ref{subsectionquotientstacks} for more discussions). For each $g\in G$ we have an isomorphism $\sigma_g:g^*\calF\simeq \calF$. It induces the following maps:
	\begin{equation*}
	\begin{aligned}
		\tau_g: R\Gamma(X,\calF)\stackrel{\simeq}\to R\Gamma(X,g_*g^*\calF)&\simeq R\Gamma(X,g^*\calF)\stackrel{\sigma_g}\longrightarrow R\Gamma(X,\calF),\\
		R\Gamma_c(X,\calF)\stackrel{\simeq}\leftarrow R\Gamma_c(X,g_!g^!\calF)&\simeq R\Gamma_c(X,g^*\calF)\stackrel{\sigma_g^{-1}}\longleftarrow R\Gamma_c(X,\calF):\tau_{g,c}.
	\end{aligned}
	\end{equation*}
	Here we used that $g_*\simeq g_!,g^*\simeq g^!, g_*g^*\simeq \mathrm{id}$ and $g_!g^!\simeq \mathrm{id}$, and the middle isomorphisms are induced by the isomorphism $g\cdot:g^{-1}X\xrightarrow{\sim} X$. This tells that the Serre duality is compatible with the $G$-actions.
\end{remark}

\subsection{Comparison of two compactly supported cohomologies}\label{subseccompare}
Let $X$ be a partially proper rigid analytic variety over $\bbQ_p$. For a vector bundle $\calF$ on $X$, we can view it as an object in $D(X)$ and define the (compactly supported) cohomology $R\Gamma^{\CS}_{(c)}(X,\calF)$ using the six functor formalism in the last subsection. On the other hand we can view $\calF$ as a solid abelian sheaf on $X$ and define the (compactly supported) cohomology $R\Gamma_{(c)}^{\it cl}(X,\calF)$ for the analytic topology of $X$ \cite{abeLazda2025proper} (for condensed version see \cite[Appendix A]{li2026locallyanalyticvectorscompleted}). Here the superscripts ``$\CS$'' and ``$\mathit{cl}$'' signify ``Clausen--Scholze'' and ``classical'' respectively.

In this subsection, we are going to show that these two compactly supported cohomologies coincide (Corollary \ref{corocoincide}).

\begin{remark}
	A quicker way to compare two compactly supported cohomologies (in smooth case) is to use the Serre duality in both settings (for classical setting see \cite[\S~5]{achinger2025compactlysupportedpadicproetale}). But we prefer to give a direct proof.
\end{remark}

The classical compactly supported cohomology has the following property:

\begin{proposition}\label{propforDresolution} 
Suppose $X$ is a partially proper rigid analytic variety and $\{U_\alpha\}$ is a locally finite covering (i.e. for any $x\in X$, there exists an open neighbourhood $U_x$ such that $U_x\cap U_{\alpha}=\emptyset$ for all but finitely many $\alpha$) of $X$ by affinoid open subspaces. For a vector bundle $\calF$ on $X$, we have naturally
	$$R\Gamma_c^{\it cl}(X,\prod_{\alpha} j_{\alpha *}j^{-1}_{\alpha}\calF)\simeq \oplus_{\alpha} R\Gamma^{\it cl}(U_\alpha,j_{\alpha}^{-1}\calF),$$
where $j_{\alpha}:U_{\alpha}\hookrightarrow X$ is the open immersion and $j_{\alpha*}$ denotes the non-derived pushforward.
\end{proposition}
\begin{proof}
	Note that partially proper adic spaces are taut so that the closure of $U_\alpha$ in $X$ is still quasi-compact. Sicne the covering is locally finite, the direct product $\prod_{\alpha} j_{\alpha *}j^{-1}_{\alpha}\calF$ is (locally) a direct sum. By an argument on the compactness of the support and \cite[Def.~A.15]{li2026locallyanalyticvectorscompleted}, it is not hard to see that for any abelian sheaf $\calF_{\alpha}$ on $U_\alpha$, we have
	$$\Gamma_c(X,\prod j_{\alpha*}\calF_{\alpha})\simeq \oplus_{\alpha} \Gamma(U_\alpha,\calF_\alpha).$$
	Choose an injective resolution of $j^{-1}_{\alpha}\calF\to \cI_{\alpha}^\bullet$. Since $j^{-1}$ is exact, we know $j_{\alpha*}\cI_{\alpha}^\bullet$ is a complex of injective sheaves. Moreover, by Tate's acyclicity and our convention that $X$ is separated, we have $H^i(V\cap U_{\alpha},j_{\alpha}^{-1}\calF)=0$ for any affinoid open $V\subset X$ and $i\geq 1$. Hence $Rj_{\alpha*}\calF=j_{\alpha*}\calF$. Consequently, $j_{\alpha *}j^{-1}_{\alpha}\calF\to j_{\alpha*}\cI_{\alpha}^\bullet$ is an injective resolution. Therefore we have
	$$R\Gamma_c^{\it cl}(X,\prod j_{\alpha *}j^{-1}_{\alpha}\calF)\simeq \Gamma_c(X,\prod j_{\alpha*}\cI_{\alpha}^\bullet)\simeq \oplus_{\alpha} \Gamma(U_\alpha,\cI_{\alpha}^\bullet)\simeq \oplus_{\alpha} R\Gamma(U_\alpha,j_{\alpha}^{-1}\calF).$$
	The result follows.
\end{proof}

\begin{lemma}\label{comparisoncohomology}
	Let $X$ be a separated analytic adic space over $\bbQ_p$ and $\calF$ be a vector bundle over $X$. We have a natural quasi-isomorphism of solid $\bbQ_p$-vector spaces $$R\Gamma^{\CS}(X,\calF)\simeq R\Gamma^{\it cl}(X,\calF),$$
	where we view $\calF\in D(X)$ for the left-hand side and view $\calF$ as a solid abelian sheaf on $|X|$ for the right-hand side.	
\end{lemma}
\begin{proof}
	Let $\{U_i=\Spa(A_i,A_i^+)\}_{i\in I}$ be an open covering of $X$ by affinoid subspaces such that $\calF|_{U_i}$ is free. Let $Y=\sqcup_{i\in I} U_i$. By $*$-descent (Example~\ref{exampleopencoverDcover}), $R\Gamma^{\CS}(X,\calF)$ is computed by 
	\begin{equation}\label{complexcechnerve}
	\lim_{[n]\in \Delta} \calF(Y^{n/X})
	\end{equation}
	where the index is the \v{C}ech nerve of the open covering $Y\to X$. However, by \v{C}ech to derived comparison (\cite[\href{https://stacks.math.columbia.edu/tag/03F9}{Tag 03F9}]{stacks-project}) and Tate's acyclicity $R\Gamma^{\it cl}(U_i,\calF|_{U_i})=H^0(U_i,\calF|_{U_i})[0]$, the classical cohomology $R\Gamma^{\it cl}(X,\calF)$ is also computed by the complex \eqref{complexcechnerve}.
\end{proof}

To compare compactly supported cohomology, we first consider local cases.

\begin{proposition}\label{naiveinterior}
	Let $X=\Spa(A,A^{\circ})$ be an affinoid rigid analytic variety and let $j:X^0\subseteq X$ be an open immersion such that $X^0$ is a Stein variety (we call such an $X^0$ a naive interior of $X$). Assume $\calF$ is a vector bundle on $X$. Then we have 
	$$R\Gamma^{\CS}_c(X^0,\calF)\simeq R\Gamma_c^{\it cl}(X^0,\calF).$$
\end{proposition}

\begin{proof}
	Suppose $X^0=\cup_{i\geq 0} V_i$ where $V_i=\Spa(A_i,A^+_i)$ form a system of strict inclusions of open affinoids. Let $j_i:V_i\to X$ be the open immersion. Let $Z:=|X\backslash X^0|$ be the closed subset of $X$. Let $\ovV_i$ be the closure of $V_i$. Then $Z=\cap_i |X\backslash \overline{V_i}|$ satisfies the condition in Lemma \ref{lemmacompute}. Let $W_i=X\backslash \overline{V_i}$, $l_i: W_i\to X$ and $i:Z\to X$. 

	Now we view these objects as open or closed objects in $X^{lo}$. The map $V_i\to V_{i+1}$ factors through $V_i^{cp}=\AnSpec(A_i,\bbZ)$ (actually the universal compactification of $V_i\rightarrow *$ is homeomorphic to $\ovV_i$ since $X^0$ is partially proper, cf. \cite[Prop.~3.4.5(1)]{achingerLaraYoucis2023geometric}), and likewise so does $j_i$, and we have a decomposition (see Remark~\ref{remarkuniversalcompletion})
	$$V_i\stackrel{j'_i} \longrightarrow V_i^{cp}\stackrel{j''_i}\longrightarrow X\stackrel{f}\longrightarrow *.$$
	Since $j''_i$ lies in $P$, we have $j''_{i*}=j''_{i!}$, and hence 
	$$R\Gamma^{\CS}_c(V_i,j_i^*\calF)=(fj''_i)_*j'_{i!}j_i^*\calF=f_*(j''_ij'_i)_!j_i^*\calF=R\Gamma^{\CS}(X,j_{i!}j^*_i\calF).$$ 
	Let $c_i:Z_i\subseteq X$ be the complement of $V_i$ in $X$. By \eqref{localizationsequence} we have an exact triangle
	$$j_{i!}j^*_i\calF\to \calF\to c_{i*}c^*_i\calF\to .$$
	Taking global sections on $X$ we get an exact triangle
	$$R\Gamma^{\CS}_c(V_i,j^*_i\calF)\to R\Gamma^{\CS}(X,\calF)\to R\Gamma^{\CS}(Z_i,c^*_i\calF)\to .$$
	Here we use $R\Gamma^{\CS}(Z_i,c^*_i\calF)$ to denote $A_{Z_i}\otimes_{(A,A^{\circ})_\square}^L M$ as an object in $D(\bbQ_{p,\square})$, where $A_{Z_i}$ is the idempotent algebra corresponding to $Z_i$ and $M\in D((A,A^{\circ})_\square)$ is the complete $A$-module corresponding to $\calF$. Taking the colimit along $i$ and applying Lemma \ref{lemmacompute}, we get a triangle
	$$R\Gamma^{\CS}_c(X^0,j^*\calF)\to R\Gamma^{\CS}(X,\calF)\to \colim_i R\Gamma^{\CS}(W_i,l_i^*\calF)\to .$$
	
	On the other hand, for the classical compactly supported cohomology, we have $R\Gamma^{\it cl}_c(X^0,\calF)\simeq R\Gamma^{\it cl}(X,j_!\calF)$. This is because $\Gamma_c(X^0,\calF)\simeq \Gamma(X,j_!\calF)$; deriving the functor and noting that the (underived) extension-by-zero functor $j_!$ is exact, we get a triangle
	$$R\Gamma^{\it cl}_c(X^0,j^{-1}\calF)\to R\Gamma^{\it cl}(X,\calF)\to  R\Gamma^{\it cl}(Z,i^{-1}\calF)\to .$$
	There is a natural map (similarly to Remark~\ref{remarkmapsfactorizationclosedopen}) $\colim_i Rl_{i*}l^{-1}_i\calF\to i_*i^{-1}\calF$ which induces (for the first isomorphism below we use that $X$ is quasi-compact)
	$$\iota: \colim_i R\Gamma^{\it cl}(W_i,l_i^{-1}\calF)\simeq R\Gamma^{\it cl}(X,\colim_i Rl_{i*}l^{-1}_i\calF)\to R\Gamma^{\it cl}(Z,i^{-1}\calF).$$
	Note that $\{W_i\}_i$ is cofinal in the system of open neighbourhoods of $Z$, so by \cite[Lem.~2.12]{abeLazda2025proper} (for the condensed version using \cite[Prop.~A.13]{li2026locallyanalyticvectorscompleted}) the map $\iota$ is a quasi-isomorphism. This implies the desired result by Lemma~\ref{comparisoncohomology}.
\end{proof}

\begin{corollary}\label{corocoincide}
	Let $X$ be a partially proper rigid analytic variety over $\bbQ_p$ and let $\calF$ be a vector bundle on $X$. Then
	$$R\Gamma^{\CS}_c(X,\calF)\simeq R\Gamma^{\it cl}_c(X,\calF).$$
\end{corollary}
\begin{proof}
	By our convention, a partially proper rigid analytic variety $X$ is countable at infinity and can be covered by naive interiors $\{X_i\}_{i\in \bbN}$ as in Proposition~\ref{naiveinterior}. By $!$-descent (Example~\ref{exampleopencoverDcover}), we have
	$$R\Gamma_c^{\CS}(X,\calF)\simeq \colim_{[n]\in \Delta^{\rm op}} \oplus_{i_0,\cdots,i_n\in \bbN} R\Gamma_c^{\CS}(X_{i_0}\times_X \cdots \times_X X_{i_n},\calF).$$
	Classical compactly supported cohomology also admits analytic descent. This is because we have a resolution of $\calF$ by $\oplus_{i_0,\cdots,i_n} j_{i_0,\cdots,i_n!}j^{-1}_{i_0,\cdots,i_n}\calF$ where $j_{i_0,\cdots,i_n}:X_{i_0}\times_X\cdots\times_X X_{i_n}\rightarrow X$ (the exactness can be checked on stalks by \cite[Cor.~A.10]{li2026locallyanalyticvectorscompleted}, cf. \cite[\S~4.7]{abeLazda2025proper}). Note that $j_{i!}$ is exact and $R\Gamma^{\it cl}_c(X,-)$ commutes with direct limits (\cite[Cor. 4.11]{abeLazda2025proper}, the proof works in solid case, see \cite[proof of Lemma A.21]{li2026locallyanalyticvectorscompleted}). Hence
	$$R\Gamma_c^{\it cl}(X,\calF)\simeq \colim_{[n]\in\Delta^{\rm op}} \oplus_{i_0,\cdots,i_n\in \bbN}R\Gamma_c^{\it cl}(X_{i_0}\times_X \cdots \times_X X_{i_n},\calF).$$
	So the result follows from Proposition \ref{naiveinterior}.
\end{proof}

\subsection{Two functional analysis results}\label{subsectionfunctionalanalysis}

In this subsection, we show that compactly supported cohomology of an equidimensional smooth affinoid rigid analytic variety concentrates in one degree and is a Smith space (Proposition \ref{propcompactsmith}). This does not follow directly from Serre duality (Theorem \ref{theoremserreduality}). We also collect a result on functional analysis for later usage (Proposition \ref{propfinite}).

Let $E$ be a finite extension of $\bbQ_p$. Let $X=\Spa(A,A^\circ)$ be a smooth affinoid rigid analytic variety of pure dimension $d$ over $E$ and $\calF$ be a vector bundle on $X$. 
\begin{lemma}\label{lemmaconserve}
 	Assume $L$ is a complex of Smith spaces which concentrates in degree $\leq 0$. Assume that $R\uHom(L,E)$ has only possible nonzero cohomology in degree zero. Then there is a quasi-isomorphism $B[0]\simeq R\uHom(L,E)$ where $B$ is a Banach space over $E$ and $L\simeq S[0]$ where $S=B^{\vee}:=\underline{\Hom}(B,E)$ is a Smith space over $E$.
 \end{lemma}
 
\begin{proof}
	Since Smith spaces are projective in $\Mod(E_{\square})$ \cite[Lem.~3.8]{rodrigues2022solid}, $R\uHom(L,E)\simeq L^{-\bullet,\vee}$ can be represented by a complex of Banach spaces concentrating  in degree $\geq 0$. Let $B:=\ker(L^{0,\vee}\to L^{-1,\vee})$ which is a Banach space. By assumption, the complex
	$$ 0\to B\to L^{0,\vee}\to L^{-1,\vee}\to \cdots \to L^{-s,\vee}\to \cdots $$
	 is an exact sequence of Banach spaces over $E$. Since $\underline{\Hom}(-,E)$ induces an exact anti-equivalence between Smith spaces and Banach spaces (cf. \cite[Thm.~3.40(1)]{rodrigues2022solid}), we get that
	 $$\cdots \to L^{-s}\to L^{-s+1}\to \cdots L^0\to B^{\vee}\to 0$$ 
	 is an exact sequence of Smith spaces.
\end{proof}

Let $m\geq 1$ and $T_E:= E\langle T_1,\cdots,T_m\rangle$ be the Tate algebra and $R_{\bbZ}=\bbZ[T_1,\cdots,T_m]$. Consider spaces 
\begin{equation}\label{equationspacesZbarZ}
	Z=\AnSpec((E\langle T_1,\cdots,T_m\rangle,\calO_E\langle T_1,\cdots,T_m\rangle)_\square), \quad \overline{Z}:=\AnSpec((E\langle T_1,\cdots,T_m\rangle,\bbZ)_\square).	
\end{equation}
Define a $(\bbZ[T_1,\cdots,T_m],\bbZ)_\square$-module
\[R_{\infty}:=(\bbZ((T_1^{-1}))/\bbZ[T_1])\otimes_{\bbZ_\square}\cdots\otimes_{\bbZ_\square} (\bbZ((T_m^{-1}))/\bbZ[T_m])\] 
and set (cf. \cite[Lem.~2.16]{rodrigues2022solid})
\[T_E^{\vee}:=(\calO_E[[T_1^{-1},\cdots,T_m^{-1}]])[1/p]\cdot T_1^{-1}\cdots T_m^{-1}\simeq \uHom_{E}(T_E,E).\] 
Consider the maps 
$g:Z\stackrel{j}\longrightarrow \overline{Z}\stackrel{g'}\longrightarrow \AnSpec E.$

\begin{lemma}\label{lemma:compact-support-tensor}
	We have an equivalence of functors $g_!(-)\simeq g'_*j_!(-)\simeq -\otimes^L_{(R_{\bbZ},\bbZ)_\square} R_{\infty}[-m]$ from $D(Z)$ to $D(E_\square)$.
\end{lemma}

\begin{proof}
	We have a pullback diagram
	\[
	\begin{tikzcd}
	\AnSpec((T_E,R_{\bbZ})_\square)
	\arrow[r]
	\arrow[d]
	& \AnSpec((T_E,\bbZ)_\square)
	\arrow[d]\\
	\AnSpec((R_{\bbZ},R_{\bbZ})_\square)
	\arrow[r]
	& \AnSpec((R_{\bbZ},\bbZ)_\square)
	\end{tikzcd}
	\]
	where each horizonal arrow is an open immersion in the sense of categorified locale by \cite[\S~8]{scholze2019condensed}. By \cite[Prop.~3.32]{andreychev2021pseudocoherent}, $\AnSpec((T_E,R_{\bbZ})_\square)\cong Z$. Assume the bottom horizonal open immersion corresponds to the idempotent algebra $B$, then it is easy to check the top horizontal open immersion corresponds to $B':=B\otimes^L_{(R_{\bbZ},\bbZ)_\square} (T_E,\bbZ)_\square$ (cf. \cite[Lem.~6.4(4)]{scholze2022complex}). By Definition~\ref{openclosed} and that $g'_*=g'_!$ is the forgetful functor, for $M\in D(Z)$, we have $g_!M=g'_*j_!(M)=[M\to M\otimes^L_{(T_E,\bbZ)_\square} B']=[M\to M\otimes_{(R_{\bbZ},\bbZ)_\square}^L B]\simeq h_!(M)$ where $h:\AnSpec((R_{\bbZ},R_{\bbZ})_\square)\to \AnSpec \bbZ_{\square}$. So it suffices to show $h_!(-)\simeq -\otimes^L_{(R_{\bbZ},\bbZ)_\square} R_{\infty}[-m]$. This is in \cite[Thm.~8.13]{scholze2019condensed}.
\end{proof}

\begin{lemma}\label{lemmatensor}
	We have a quasi-isomorphism in $D((R_{\bbZ},\bbZ)_\square)$:
	$$T_E\otimes^L_{(R_{\bbZ},\bbZ)_\square} R_{\infty}\simeq (\calO_E[[T_1^{-1},\cdots,T_m^{-1}]])[1/p]\cdot T_1^{-1}\cdots T_m^{-1} =T_E^{\vee}.$$
\end{lemma}

\begin{proof}
	Note that we have\begin{equation*}
		T_E\simeq E\langle T\rangle^{\otimes m}, \quad (R_{\bbZ},\bbZ)_\square\sim (\bbZ[T],\bbZ)^{\otimes m}_\square,\quad T_E^{\vee}\simeq(\calO_E[[T^{-1}]][\frac{1}{p}]\cdot T^{-1})^{\otimes m}.
	\end{equation*}
	So it suffices to show the case $m=1$. We have a resolution
	$$0\longrightarrow \bbZ[[X]][T]\stackrel{XT-1}\longrightarrow \bbZ[[X]][T]\longrightarrow \bbZ((T^{-1}))\longrightarrow 0.$$
	One computes that (using that $E\langle T\rangle=\calO_E\langle T\rangle[1/p]$,  $\otimes^L$ commutes with colimits and \cite[Prop.~6.3]{scholze2019condensed}):
	$$E\langle T\rangle\otimes_{(\bbZ[T],\bbZ)_\square}^L \bbZ[[X]][T]\simeq E\langle T\rangle\otimes_{\bbZ_\square}^L \bbZ[[X]]\simeq \calO_E[[X]]\langle T\rangle [1/p]$$
	Thus $E\langle T \rangle\otimes_{(\bbZ[T],\bbZ)_\square}^L (\bbZ((T^{-1}))/\bbZ[T])\simeq \Cone(E\langle T\rangle\rightarrow \calO_E[[T^{-1}]]\langle  T\rangle[1/p])\simeq \calO_E[[T^{-1}]]\cdot T^{-1}[1/p]$.
\end{proof}

\begin{proposition}\label{propcompactsmith}
	Let $X$ be a smooth affinoid rigid variety of pure dimension $d$. We have \[R\Gamma_c(X,\calF)\simeq H^d_c(X,\calF)[-d]\] 
	and $H^d_c(X,\calF)$ is a Smith space over $E$.
\end{proposition}

\begin{proof}
  First, it suffices to show $R\Gamma_c(X,\calF)$ can be represented by a complex of Smith spaces concentrated in degree $\leq d$: By Serre duality (Theorem~\ref{theoremserreduality}), the dual of $R\Gamma_c(X,\calF)[d]$ is quasi-isomorphic to a Banach space concentrated in degree $0$. Then by Lemma \ref{lemmaconserve}, $R\Gamma_c(X,\calF)[d]$ is quasi-isomorphic to a Smith space concentrate in degree $0$.
  
   If the statement holds for the structure sheaf $\calO$, since $\calF$ can be resolved by sheaves of the form $\calO^{\oplus m}$, we find $R\Gamma_c(X,\calF)$ can be represented by a complex of Smith spaces concentrated in degree $\leq d$. So it suffices to consider the case $\calF=\calO$.

	Now we show that it suffices to assume $A=E\langle T_1,\cdots,T_m\rangle/(f_1,\cdots,f_r)$ where $f_1,\cdots,f_r$ form a regular sequence in $E\langle T_1,\cdots,T_m\rangle$ such that $d=m-r$. We can find finitely many open affinoids $U_1,\cdots,U_n$ covering $X$ such that each $U_i$ are complete intersection. The existence of such cover can be seen using, for example, the Jacobian criterion \cite[Prop.~1.7.1, Cor.~1.6.10]{huber1996etale}. It is clear that a rational subdomain of a complete intersection is a complete intersection thus we can further assume $U_i$ are rational subdomains of $X$. In particular, any intersections of $U_i$ are also complete intersections since intersections of rational subdomains are rational subdomains \cite[\S~3.3, Prop.~14]{bosch2014formal}. Consider the open covering $\sqcup_i U_i\to X$. By $!$-descent and our assumption, we have
	$$[H^d_c(U_1\cap \cdots \cap U_n,\calO)\to \cdots \to \oplus_{1\leq i<j\leq n} H^d_c(U_i\cap U_j,\calO)\to \oplus_{1\leq i\leq n} H^d_c(U_i,\calO)] \stackrel{\simeq}\longrightarrow R\Gamma_c(X,\calO)$$
	where the term $\oplus_{1\leq i\leq n} H^d_c(U_i,\calO)$ lies in degree $d$. Thus $R\Gamma_c(X,\calO)$ can be represented by a complex of Smith spaces concentrated in degree $\leq d$ if the statement is proved for complete intersections.
	
	Now assume that $A=E\langle T_1,\cdots,T_m\rangle/(f_1,\cdots,f_r)$ is a complete intersection and $\calF$ is free of rank one. Put $X=\AnSpec((A,A^\circ)_\square)$ and set $Z,\ovZ$ as in \eqref{equationspacesZbarZ}. We have a commutative diagram
	\[
	\begin{tikzcd}
	X
	\arrow[r,"i"]
	\arrow[rd,"f"']
	& Z
	\arrow[r,"j"]
	\arrow[d,"g"]
	& \overline{Z}
	\arrow[ld,"g'"]\\
	& \AnSpec E &
	\end{tikzcd}
	\]
	where $i$ and $g'$ lie in the class $P$ of Definition~\ref{definitionsuitabledecomposition}. Hence $R\Gamma_c(X,\calO)= f_!\calO\simeq g'_*j_!i_*\calO$. By Lemma~\ref{lemma:compact-support-tensor} and Lemma~\ref{lemmatensor}, we have $$g'_*j_!i_*\calO\simeq \mathrm{Kos}(f_1,\cdots,f_r; T_E)\otimes^L_{(T_E,\bbZ)_\square} T_E^{\vee}[-m],$$
	where $\mathrm{Kos}$ denotes the Koszul complex so that $A\simeq \mathrm{Kos}(f_1,\cdots,f_r; T_E)$ \cite[\href{https://stacks.math.columbia.edu/tag/062F}{Tag 062F}]{stacks-project}. The above sequence is formally the dual (via $\uHom_{E}(-,E)$) of $\mathrm{Kos}(f_1,\cdots,f_r; T_E)$ by the self duality of the Koszul complexes \cite[§~A.2, Ch.~IV]{serre1975algebre}. Since $\uHom_{E}(-,E)$ is exact on strict exact sequence of Banach spaces \cite[Prop.~3.5]{rodrigues2022solid}, $ \mathrm{Kos}(f_1,\cdots,f_r; T_E)\otimes^L_{(T_E,\bbZ)_\square} T_E^{\vee}$ is quasi-isomorphic to $\uHom_E(A,E)[r]$. We finished the proof.
\end{proof}
Recall the \emph{trace-class functor} (\cite[Def.~A.46]{bosco2021drinfeldproetale} and see \cite[Def.~A.13]{bosco2023rationalpadicHodge} for the derived version)
\begin{equation}\label{equationdefinitiontraceclassfunctor}
	(-)^{\trace}: \Mod(E_{\square})\rightarrow \Mod(E_{\square}), W\mapsto (S\in \mathrm{ExtDisc}\mapsto (E_{\square}[S]^{\vee}\otimes^{\square}_EW)(*)).
\end{equation}
If $W\in \Mod(E_{\square})$, then a map $P\rightarrow W$ from a Smith space $P$ is \emph{trace-class} if and only if the map factors through $W^{\rm tr}\rightarrow W$. The last map is an isomorphism if and only if $W$ is nuclear \cite[Def.~3.14]{rodrigues2022solid}. See \cite[Lem.~A.47]{bosco2021drinfeldproetale} and also \cite[Lem.~A.14]{bosco2023rationalpadicHodge}.

By the flatness in \cite[Prop. A.27]{bosco2021drinfeldproetale}, the functor $(-)^{\trace}$ is exact on Smith spaces, and commutes with direct sum. Moreover, if $W$ is Smith, then by the proof of \cite[Lem. A.58]{bosco2021drinfeldproetale}, $W^{\trace}=W^B$ is the associated Banach space \cite[Def.~3.34(3)]{rodrigues2022solid}: choosing a lattice $W_0$ of $W$ then 
\begin{equation}\label{equationdefinitionassociatedBanach}
	W^B:=\lim_n (\underline{W_0(*)_{\dis}}/\varpi_E^n)[1/\varpi].
\end{equation}
\begin{corollary}\label{corollarytraceclass}
		Suppose there is a strict open immersion of smooth affinoid rigid analytic varities $j:X_1\to X_2$ pure of dimension $d$ such that there is a Stein variety $X^0_1$ with open immersions $X_1\subset X_1^0\subset X_2$. Let $\calF$ be a vector bundle on $X_2$. Then the induced map $j_c:R\Gamma_c^{\CS}(X_1,\calF)\to R\Gamma_c^{\CS}(X_2,\calF)$ is a trace-class map in $D(E_\square)$.
\end{corollary}

\begin{proof}
	 The compactly supported cohomologies of $X_1$, $X_1^0$, and $X_2$ are concentrated in degree $d$ (by Proposition~\ref{propcompactsmith}), and $H^d_c(X_1,\calF)$ is a Smith space. After shifting by $d$, the map $j_c$ factors as $$j_c:H^d_c(X_1,\calF)\to H^d_c(X_1^0,\calF)\to H^d_c(X_2,\calF).$$ 
	 Since $X_1^0$ is Stein, $H^d_c(X_1^0,\calF)$ is nuclear by \cite[Cor.~A.50, Thm.~A.43]{bosco2021drinfeldproetale} being a direct limit of Banach spaces (dual of the Fr\'echet space $H^0(X_1^0,\calF^{\vee}\otimes \Omega_{X_{1}^0})$), the first map in the factorization is trace-class, so is $j_c$ by \cite[Rem.~A.45]{bosco2021drinfeldproetale}.
\end{proof}

Now we come to the second result, this is a condensed version of the $p$-adic Schwarz lemma \cite[Satz~1.5]{kiehl1967endlichkeitssatz}.

\begin{definition}[{\cite[Def.~3.33]{rodrigues2022solid}}]\label{definitioncompactmorphism}
 	For two Banach $E$-vector spaces $V,W$, a morphism $f:V\to W$ is called compact if the dual map between Smith spaces $f^\vee: W^\vee\to V^\vee$ is trace-class where $V^\vee:=\underline{\Hom}(V,E)\in \Mod(E_\square)$.
 \end{definition}
 
 \begin{remark}
	Note that classical compact morphism $f:W\to V$ between Banach spaces is a compact morphism in condensed setting. One can check that in this case, $f^\vee: V^\vee\to W^\vee$ factors through $(W^\vee)^B$.  
\end{remark}
 
 \begin{lemma}\label{lemmatrace}
 	Assume $W^\bullet,V^\bullet$ are two bounded complexes of solid Smith $E$-vector spaces. Suppose we have compatible morphisms $f^i: W^i\to V^i$ which are trace-class. Then regarding $f^\bullet$ as a morphism in $D(E_\square)$, it is trace-class (\cite[Def.~8.1]{scholze2022complex}). In particular, if $f^\bullet$ induces a quasi-isomorphism, then $W^\bullet$ (so is $V^\bullet$) is a perfect complex as an object in $D(E_\square)$.
 \end{lemma}
 
 \begin{proof}
 	By our assumption, $f^i$ factors as $W^i\stackrel{f^i_1}\longrightarrow V^{i,B}\stackrel{f^i_2}\longrightarrow V^i$. It is clear that $V^{\bullet,B}$ form a complex of Banach spaces and $f^\bullet_2$ is a morphism of complexes. By uniqueness of the factorization, $f^\bullet_1$ is also a morphism of complex. Now we view $W^\bullet,V^{\bullet,B},V^\bullet$ as objects in $D(E_\square)$. Since $V^{\bullet,B}$ is a nuclear complex, we know $f^\bullet_1$ is trace-class using \cite[Lem.~A.14]{bosco2023rationalpadicHodge} (here we use boundedness so that $W^{\bullet}$ is compact), thus $f^\bullet=f_2^\bullet\circ f_1^\bullet$ is trace-class by \cite[Lem.~8.2]{scholze2022complex}.
 	
 	For the second claim, we know $W^\bullet$ is a compact object by \cite[Lem.~5.40]{andreychev2021pseudocoherent} in $D(E_\square)$ and $\mathrm{id}_{W^\bullet}$ is equivalent to a trace-class map. By \cite[Prop.~9.4]{scholze2022complex}, $W^\bullet$ is a dualizable object in $D(E_\square)$ and by \cite[Cor.~5.51.1]{andreychev2021pseudocoherent} it is isomorphic to a perfect complex.
 \end{proof}

\begin{proposition}\label{propfinite}
	Assume $f^\bullet: W^\bullet\to V^\bullet$ is a morphism of bounded complexes of Banach $E$-vector spaces such that each $f^i$ is a compact morphism and $f^\bullet$ induces a quasi-isomorphism. Then $W^\bullet$ (so is $V^\bullet$) is perfect as an object in $D(E_\square)$.
\end{proposition}
 
\begin{proof}
	Let $V_s^i=V^{-i,\vee}$ and $W_s^i=W^{-i,\vee}$ be Smith spaces. There is a morphism of bounded complexes of Smith spaces $f^\bullet_s: V_s^\bullet\to W^\bullet_s$ whose dual (taking $\uHom_{E}(-,E)$ termwisely) is $f^\bullet$. By our assumption, $f^i_s$ are trace-class. We claim that $f^\bullet_s$ induces a quasi-isomorphism. Then by Lemma \ref{lemmatrace}, this implies $V_s^\bullet$ is a perfect complex thus $W^\bullet=R\uHom(W_s^\bullet,E)$ is also a perfect complex.
	
	To show $f^\bullet_s$ induces a quasi-isomorphism, let $C:=\Cone(f_s^\bullet)$. Since $W^\bullet\simeq R\uHom(W_s^\bullet,E)$, $V^\bullet\simeq R\uHom(V_s^\bullet,E)$ and $\Cone(f)\simeq 0$, we get that $R\uHom(C,E)\simeq 0$. Applying Lemma \ref{lemmaconserve}, we get $C\simeq 0$. This finished the proof.
\end{proof}

\section{Equivariant vector bundles}\label{sectionequivariantvectorbundle}
In this section, we relate Bernstein--Zelevinsky duality of $p$-adic Lie groups with the Grothendieck--Serre duality of equivariant vector bundles (Theorem~\ref{maintheorem}) and discuss equivariant vector bundles on the Drinfeld spaces in \S~\ref{subsec:Drinfeldspaces}. 

Let $G$ be a $p$-adic Lie group and let $K\subset G$ be an open compact subgroup. Let $E$ be a finite extension of $\bbQ_p$.
\subsection{Bernstein--Zelevinsky duality of locally analytic representations}
We first recall the definition of the Bernstein--Zelevinsky duality of locally analytic representations. For $V\in \Rep^{\la}_{E_\square}(G)$, its Bernstein--Zelevinsky dual (\cite[Def.~5.4]{strauch2025bernsteinzelevinskydualitylocallyanalytic}) is given by:
\begin{equation}\label{definitionBZdual}
	\bbD_{\BZ}(V):=R\uHom_{\calD^{\la}(G,E)}(V,\calC_c^{\la}(G,E))\in D_{\square}(\calD^{\la}(G,E)).
\end{equation}

\begin{lemma}\label{lemmadualla}
	Assume $W$ is a bounded complex of locally analytic representations of $K$ such that underlying modules of each $W^i$ is a Smith space. Then $W^\vee\in \Rep^{\la}_{E_\square}(K)$ where $W^\vee=R\underline{\Hom}_E(W,E)$ is endowed with dual action of $K$. 
\end{lemma}

\begin{proof}
	Note that $R\underline{\Hom}_E(W_i,E)=\underline{\Hom}_E(W_i,E)$ using that Smith spaces are projective, the fact $\uHom_E(V,-)(S)=\uHom_E(V\otimes_{E_{\square}}E^{\square}[S],-)$) and \cite[Prop.~3.12]{rodrigues2022solid}. Hence $(W^\vee)_i=(W_{-i})^{\vee}$ is a Banach spaces. By \cite[Prop.~4.41, Thm.~4.36(2)]{rodrigues2022solid}, since $W_{-i}$ is a $\calD^{\la}(K,E)$-module, $(W^\vee)_i$ is a locally analytic representation. Hence $W^{\vee}\in \Rep^{\la}_{E_\square}(K)$ since the latter is a full subcategory of $D_{\square}(\calD^{\la}(K,E))$ stable under small colimits \cite[Prop.~3.3.3]{JC23}.
\end{proof}
We write 
\[\cInd_{K}^{G}(-)=\calD^{\la}(G,E)\otimes^{L\square}_{\calD^{\la}(K,E)}-:D_{\square}(\calD^{\la}(K,E))\rightarrow D_{\square}(\calD^{\la}(G,E))\] 
for the compact induction functor.
\begin{proposition}[{\cite[Thm.~5.3]{strauch2025bernsteinzelevinskydualitylocallyanalytic}}]\label{BZdual}
	Assume $W=W^{\bullet}$ is a bounded complex of locally analytic representations of $K$ such that underlying modules of each $W^i$ is a Smith space. Then we have a quasi-isomorphism in $D(\Mod_{E_\square}(D^{\la}(G,E)))$
	$$\bbD_{\BZ}(\cInd_K^G W)\simeq \cInd_K^G W^\vee.$$
\end{proposition}
\begin{proof} 
	Let $\chi=\Delta_G^*|_K$ be the right character of Appendix~\ref{sec:STdual}, and let $\iota$ be the involution $g\mapsto g^{-1}$ of $G$. We calculate that
	$$\begin{aligned}
	\bbD_{\BZ}(\cInd_K^G W)
	&\simeq R\underline{\Hom}_{\calD^{\la}(G,E)}(\cInd_K^G W,\calC_c^{\la}(G,E))\\
	&\simeq R\underline{\Hom}_{\calD^{\la}(K,E)}(W,\calC_c^{\la}(G,E))\\
	&\simeq R\underline{\Hom}_{\calD^{\la}(K,E)}(E,R\underline{\Hom}_E(W,\calC_c^{\la}(G,E))\\
	&\simeq R\underline{\Hom}_{\calD^{\la}(K,E)}(E,W^\vee\otimes^L_{E_\square} \calC_c^{\la}(G,E))\\
	&\simeq (\iota(W^{\vee})\otimes^L_{E_\square} \chi)\otimes_{\calD^{\la}(K,E)}^L \calC_c^{\la}(G,E)[-\dim G]\\
	&\simeq  ((\iota(W^{\vee})\otimes^L_{E_\square} \chi)\otimes^L_{\calD^{\la}(K,E)} \calC^{\la}(K,E)[-\dim G])\otimes^L_{\calD^{\la}(K,E)} \calD^{\la}(G,E)\\
	&\simeq  R\uHom_{\calD^{\la}(K,E)}(E, W^{\vee}\otimes^L_{E_\square} \calC^{\la}(K,E))\otimes^L_{\calD^{\la}(K,E)} \calD^{\la}(G,E)\\
	&\simeq W^\vee \otimes^L_{\calD^{\la}(K,E)}\calD^{\la}(G,E) \simeq \cInd_K^G W^\vee
\end{aligned}$$
where the second quasi-isomorphism is by adjunction; the third one follows from \cite[Prop.~A.2.1(4)]{JC23}; the fourth one is because $W$ is compact and $\calC_c^{\la}(G,E)$ is nuclear; the fifth and seventh isomorphisms use \cite[Prop.~3.2.12]{JC23}; the eighth is by Lemma \ref{lemmadualla} which tells that $W^\vee\simeq (W^\vee)^{R\la}=R\uHom_{\calD^{\la}(K,E)}(E, W^{\vee}\otimes^L_{E_\square} \calC^{\la}(K,E))$ (\cite[Def.~3.2.2]{JC23}).
\end{proof}
\subsection{Locally analytic vector bundles}\label{subsectionquotientstacks}
In this subsection, we let $G$ be a $p$-adic Lie group. Let $E$ be a finite extension over $\bbQ_p$ and let $X$ be an analytic adic space over $E$ with a continuous action of $G$. We discuss the quotient analytic stack $[X/G^{\la}]$ and vector bundles on it, viewing $X$ as an analytic stack over $E_{\square}$. From now on, by an analytic stack we mean a Tate stack (Definition~\ref{definitionTatestack}).

Let $H\subset G$ be a open compact subgroup. We also write $H$ for the analytic adic space $\Spa(\calC(H,E),\calC(H,\calO_E))$ and 
\[H^{\la}:=\AnSpec((\calC^{\la}(H,E),\calO_E)_{\square}).\] 
Since $G=\sqcup_{g\in G/H}g.H$, we have a group stack $G^{\la}=\sqcup_{g\in G/H}g.H^{\la}$. 

If $G$ acts on an analytic adic space $X$ over $E$, we say that the action of $G$ on $X$ is continuous if the (a priori set-theoretic) action map 
\begin{equation}\label{equationactionGX}
	\mathrm{act}:G\times X\rightarrow X
\end{equation}
extends to an action of the adic space $G$ on $X$, and the action is locally analytic if \eqref{equationactionGX} is continuous and extends to an action 
\begin{equation}\label{equation}
	\mathrm{act}^{\la}:G^{\la}\times X\to X
\end{equation}
where $X$ is viewed as an analytic stack over $E_{\square}$. 
\begin{definition}\label{definitionquotientstack}
	Let $X$ be an analytic adic space with a locally analytic action of $G$. We define an analytic stack
	\[[X/G^{\la}]:=\colim_{[n]\in \Delta^{\rm op}} (G^{\la})^n\times X,\]
	the geometric realization of the simplicial object $((G^{\la})^n\times X)_{[n]\in \Delta^{\rm op}}$ described in \cite[Def.~4.3.1]{JC23}.
\end{definition}
By the construction of the six functor formalism, we have (see discussions after \cite[Exam.~3.4.5]{heyer20246functor})
\begin{equation}\label{equationDXG}
	D([X/G^{\la}])=\lim_{[n]\in \Delta}D((G^{\la})^n\times X).
\end{equation}
If $X=*$, then by \cite[Thm.~4.3.3]{JC23}: 
\begin{equation}\label{equationstackrepresentation}
	D([*/G^{\la}])\simeq \Rep^{\la}_{E_\square}(G).
\end{equation} 
If $X=\AnSpec A$ where $A$ is a solid analytic affinoid ring (Definition~\ref{definitionaffinoidanalyticring}) and $G$ is compact, then by \cite[Prop.~2.47]{mikami2026finiteness} (taking $h=\infty$ in \textit{loc. cit.})
\[D([X/G^{\la}])\simeq \mathrm{Rep}_A^{\la}(G)\]
where the latter is the category of semilinear $A$-modules in \cite[Def.~2.46]{mikami2026finiteness}.

We will focus on classical examples of equivariant vector bundles. 
\begin{lemma}\label{lemmaecontinuousactionrigid}
	Let $X$ be a rigid analytic variety over $E$ equipped with an action of $G$. Suppose that for any affinoid open $U\subset X$, there exists an open compact subgroup $K_U\subset G$ which stabilizes $U$ and the induced action $K_U\times \calO(U)\rightarrow \calO(U)$ is continuous (this is the definition of continuous action in \cite[Def.~2.1, Lem.~2.2]{scholze2018lubintate}).  Then the $G$-action on $X$ extends to an action of $G$ or $G^{\la}$ on $X$ as adic spaces or analytic stacks respectively.
\end{lemma}
\begin{proof}
	To extend the action $G^{(\la)}\times X\rightarrow X$, it suffices to work locally and define actions $ G^{(\la)}\times U\rightarrow X$ for $U\subset X$ affinoid. Since $U$ is stabilized by a compact open subgroup of $G$, we may easily reduce to the case when $G$ is a compact group and $X=\Spa(A,A^{\circ})$ is an affinoid space. Since the action of $G$ on $X$ is continuous, the orbit map $A\rightarrow \Hom(G,A), v\mapsto (g\mapsto g.v)$ has image in $\calC(G,A)$ the continuous functions from $G$ to $A$, which gives a map $A\rightarrow \calC(G,A)\simeq \calC(G,E)\otimes_{E_{\square}}A$ (use \cite[(2.1.3)]{emerton2017locally} and \cite[Lem.~3.13]{rodrigues2022solid}). This makes $A$ a comodule over the Hopf algebra $\calC(G,E)$ (cf. \cite[Lem.~6.1.3]{JC23}). Taking the adic spectra we get the map $G\times X\rightarrow X$. To get the locally analytic statement, by \cite[Lem.~4.3.4]{JC23}, it suffices to show that $G$ acts locally analytically on $A$. This follows from Lemma~\ref{lemmacontla} below.
\end{proof}

\begin{lemma}\label{lemmacontla}
	Let $A$ be an affinoid algebra over $E$, and let $G$ be a $p$-adic Lie group acting continuously on $A$ by $E$-algebra automorphisms. Let $M$ be a finite $A$-module equipped with a continuous $A$-semilinear action of $G$ (with respect to the canonical Banach topology on $M$). Then the actions of $G$ on $A$ and on $M$ are locally analytic.
\end{lemma}
\begin{proof}
	Choose a surjection $E\langle T_1,\cdots,T_n\rangle\twoheadrightarrow A$ and let $A_0$ be the image of $\calO_E\langle T_1,\cdots,T_n\rangle$ which is an open lattice of $A$. We may suppose that $G$ is compact and after shrinking $G$, we may assume that $A_0$ is stabilized by $G$ (by continuity). The group $G$ acts continuously, hence smoothly, on the discrete space $A_0/\varpi_E$  which is finitely generated over $\calO_E/\varpi_E$. Thus there exists an open subgroup $H\subset G$ which fixes a set of generators of $A_0/\varpi_E$. Consequently, as  $G$ acts by algebra endomorphisms, $H$ acts trivially on $A_0/\varpi_E$. Applying \cite[Prop.~3.4.1]{JC23}, the $G$-action on $A$ is locally analytic. The proof for $M$ is similar, cf. \cite[Exam.~2.3.8]{camargo2026geometricSen}. 
\end{proof}
Next, suppose that $X$ is an analytic adic space with a continuous $G$-action $\mathrm{act}:G\times X\rightarrow X$. Let $\calF$ be a $G$-equivariant vector bundle on $X$. This means a locally free $\calO_X$-module with an isomorphism $\mathrm{act}^*\calF\simeq \pi^*\calF$ satisfying the usual cocycle condition, where $\pi:G\times X\rightarrow X$ is the natural projection. Similarly, we say the $G$-equivariant vector bundle $\calF$ is locally analytic if the $G$ action induces an isomorphism $(\mathrm{act}^{\la})^{*}\calF\simeq (\pi^{\la})^{*}\calF$ where $\pi^{\la}:G^{\la}\times X\to X$ is the natural projection.
\begin{lemma}\label{lemmaequivariantbundle}
	Let $X$ be a rigid analytic variety over $E$ equipped with a continuous action of $G$ satisfying the assumption in Lemma~\ref{lemmaecontinuousactionrigid}. Let $\calF$ be a $G$-equivariant vector bundle on $X$. Then $\calF$ is locally analytic and corresponds to an object in $D([X/G^{\la}])$.
\end{lemma}
\begin{proof}
	The vector bundle $\calF$ is locally analytic by Lemma~\ref{lemmacontla}. Assume that an open affinoid $U=\Spa(A,A^+)\subset X$ is stabilized by an open compact subgroup $K$ of $G$. The space $\calC^{\la}(K,E)$ is a filtered colimit of Banach spaces, hence is flat over $E_{\square}$ \cite[Cor.~A.65]{bosco2021drinfeldproetale}. Thus $(K^{\la})^{n}\times U=\AnSpec ((A,A^+)_{\square}\otimes_{E_{\square}} \calC^{\la}(K^n,E))$. It's then formal to check that $\calF(U)$ defines an object in $D(U/K^{\la})$ and $\calF$ defines an object in $D(X/G^{\la})$.  
\end{proof}
\begin{lemma}\label{lemmaDcover}
	Let $X$ be an analytic adic space equipped with a locally analytic action of $G$. The projection $p_X: X\to [X/G^{\la}]$ is a $D$-cover (Definition~\ref{defDtopo}). In particular, $p_X$ satisfies universal $*$-descent and universal $!$-descent.
\end{lemma}
\begin{proof}
	By \cite[Prop.~6.2.3]{camargo2024analytic}, the map $*\to [*/G^{\la}]$ is a $D$-cover. The result follows since universal $*$-descent and $!$-descent are stable under pullback by definition and $p_X$ is the pullback of $*\rightarrow [*/G^{\la}]$ along $[X/G^{\la}]\rightarrow [*/G^{\la}]$.
\end{proof}
\begin{lemma}\label{lemmapullbackcohomology}
	Let $f\colon X\rightarrow Y$ be a $G^{\la}$-equivariant morphism of analytic adic spaces equipped with locally analytic actions of $G$. The square
	\[\begin{tikzcd}
		X \ar[r, "f"] \ar[d, "p_X"'] & Y \ar[d, "p_Y"] \\
		{[X/G^{\la}]} \ar[r, "\bar f"'] & {[Y/G^{\la}]}
	\end{tikzcd}\]
	is cartesian.  If $f$ is $!$-able (i.e. in the class $E''$ after Definition~\ref{definitionTatestack}), then $\bar f$ is $!$-able, and $p_Y^*\bar f_!\simeq f_!p_X^*$ as functors $D([X/G^{\la}])\rightarrow D(Y)$. Moreover, the natural morphism $p_Y^*\bar f_*\to f_*p_X^*$ is an equivalence of functors $D([X/G^{\la}])\rightarrow D(Y)$ if $f$ is quasi-compact and separated.
\end{lemma}
\begin{proof}
	Both $p_X$ and $p_Y$ are pullbacks of $*\rightarrow [*/G^{\la}]$, hence the square is cartesian. By Lemma~\ref{lemmaDcover}, $p_Y$ satisfies universal $!$-descent. Since to check if $h_!$ is defined for a morphism $h$, one can check $!$-locally on the target, we see $\bar f$ is $!$-able if $f$ is \cite[Thm.~3.4.11]{heyer20246functor}. The equivalence $p_Y^*\bar f_!\simeq f_!p_X^*$ is the proper base change \cite[Prop.~3.1.8]{heyer20246functor} since $p_{Y}$ is a $D$-cover.

	We remain to prove the last statement. The natural morphism $p_Y^*\bar f_*\to f_*p_X^* $ arises as usual from the adjoint of $f^*p_{Y}^*\bar f_*\simeq p_X^{*}\bar f^*\bar f_*\rightarrow p_{X}^*$. For any map $\alpha:[n]\rightarrow [m]$ in $\Delta$, write $\psi_{\alpha}$ for the corresponding maps of the simplicial objects, and consider the following Cartesian diagram:
	\begin{equation}\label{equationdiagramsimplicalXY}
		\begin{tikzcd}
			(G^{\la})^m\times X \ar[d, "\psi_{X,\alpha}"] \ar[r, "f_m"'] & (G^{\la})^m\times Y \ar[d, "\psi_{Y,\alpha}"] \\
			(G^{\la})^n\times X\ar[r, "f_n"'] & (G^{\la})^n\times Y.
		\end{tikzcd}
	\end{equation}
	We claim that we have equivalences 
	\begin{equation}\label{equationequivalencebasechangefmnstarphi}
		f_{m*}\psi_{X,\alpha}^*\simeq \psi_{Y,\alpha}^*f_{n*}
	\end{equation} 
	This equivalence holds when $G$ is compact and both $X$ and $Y$ are affinoid by the proof of \cite[Cor.~2.42]{mikami2026finiteness}. We prove our case by reducing to the affinoid cases. We may assume $\alpha$ is a degeneracy map or a face map. 
	
	We prove the equivalence \eqref{equationequivalencebasechangefmnstarphi} in the case when $\psi_{X,\alpha}$ is the action map $G^{\la}\times X\rightarrow X$; the other cases are similar or easier. Take an affinoid cover $(U_i)_{i}$ with open compact stabilizers $K_i$ of $U_i$. Write $V_i=f_0^{-1}(U_i)$. Consider the map $\sqcup_{g\in G/K_i} K_ig^{-1}\times g.V_i \rightarrow \sqcup_{g\in G/K_i} K_ig^{-1}\times g.U_i$. We note that open immersions are $D$-\'etale by the discussion in Example~\ref{exampleopencoverDcover}, hence $*$-pushforwards commute with $*$-pullback along open immersions by \cite[Lem.~4.5.13]{heyer20246functor}. Since $ \sqcup_{g\in G/K_i} K_ig^{-1}\times g.U_i$ is an open covering of the preimage of $U_i$ in $G^{\la}\times Y$, and $f_{1*}$ commutes with direct product, we can check \eqref{equationequivalencebasechangefmnstarphi} after pulling back \eqref{equationdiagramsimplicalXY} to each $K_ig^{-1}\times g.U_i$. For this we can assume that $G$ is compact and $Y$ is affinoid. Let $i:V\hookrightarrow X$ be an open immersion of an affinoid subspace and we may assume that $V$ is stabilized by an open subgroup $K\subset G$. We can cover $G\times X$ by $Kg^{-1}\times g.V,g\in G/K$. By the affinoid case in \cite[Cor.~2.42]{mikami2026finiteness}, we have $(f_{1}|_{Kg^{-1}\times g.V})_*(\psi_{X,\alpha}|_{Kg^{-1}\times g.V})^*\simeq (\psi_{Y,\alpha}|_{Kg^{-1}\times Y})^*(f_{0}|_{V})_*$. Since $f$ is quasi-compact and separated, we may take a finite affinoid open cover $(U_i)_{i=1,\cdots,l}$ of $X$ and let $K$ be an open subgroup of $G$ stabilizing all $U_i$, as well as their affinoid intersections $U_J=\cap_{i\in J}U_{i}$ for any $J\subset \{1,\cdots,l\}$. Write $i_{J}:U_J\rightarrow X$ for the open immersion. Since an open covering is a $D$-cover, by descent, $\calF$ is isomorphic to the total complex of the alternating \v{C}ech complex  $[\prod_{1\leq j\leq l}i_{\{j\}*}\calF|_{U_j}\rightarrow \prod_{|J|=i}i_{J*}\calF|_{U_J}\rightarrow \cdots \rightarrow j_{\{1,\cdots,l\}*}\calF|_{U_1\cap\cdots \cap U_{l}}]$ which is a finite limit of the objects $i_{J*}\calF|_{U_J}$. Since both $\Psi_{\alpha}^*$ and $f_*$ commute with finite limits (as $D(X)$ is stable), the equivalence \eqref{equationequivalencebasechangefmnstarphi} follows.

	Hence the diagram obtained after applying $D^*(-)$ to \eqref{equationdiagramsimplicalXY} (with $*$-pullbacks) is right adjointable in the sense of \cite[Def.~4.7.4.13]{lurie2017higher}, with right adjoints the $*$-pushforwards. Since $\bar f_*$ is by definition of the right adjoint of $\bar f^*=\lim_{[n]}f_n^*$ \eqref{equationDXG} and right adjoints commute with limits, it is formal to verify that $\bar f_*=\lim_{[n]}f_{n*}$ (cf. \cite[Cor.~4.7.4.18]{lurie2017higher}). The result follows.
\end{proof}
\begin{remark}
	The last statement of Lemma~\ref{lemmapullbackcohomology} is wrong without the quasi-compact assumption. For example if we take $Y=*$ and $X=\Omega_{\mathrm{Dr},E}$ a Stein space and take $\calF$ a $G$-equivariant vector bundle on $\Omega_{\mathrm{Dr},E}$, then $H^0(\Omega_{\mathrm{Dr},E},\calF)$ is in general not locally analytic, but the linear dual of a locally analytic representation of compact type. Hence $f_*\calF$ does not lie in the image of $D([*/G^{\la}])\rightarrow D(*)$ via $*$-pullback.
\end{remark}
\subsection{Two resolutions of Fargues}\label{subsectionresolutionFargues}\label{subsectionconclude}
Let $X$ be a smooth partially proper rigid analytic variety over $E$ (countable at infinity) of pure dimension $d$. Suppose $G$ acts continuously on $X$ as in Lemma~\ref{lemmaecontinuousactionrigid}. Let $\calF$ be an equivariant vector bundle on $X$ defining an object in $D([X/G^{\la}])$ by Lemma~\ref{lemmaequivariantbundle}. We follow the strategy of Fargues in \cite{Farguesdualite} to calculate the Bernstein--Zelevinsky dual of $R\Gamma_c(X,\calF)\in \Rep^{\la}_{E_\square}(G)$ using \v{C}ech complexes.

We take $U\subset X$ be an affinoid open subset and we assume the following assumption:
\begin{assumption}\label{assumptionaffionid}
	We assume that the affinoid open $U\subseteq X$ is stablized by a compact open $K\subseteq G$ such that $GU=X$ and the set $\{g\in G/K: U\cap g.U\neq \emptyset\}$ is a finite set. 
\end{assumption}

For a compact open subgroup $K\subseteq G$, we get a map $q:[*/K^{\la}]\to [*/G^{\la}]$. It is not difficult to see that $q^*$ corresponds, under the equivalence \eqref{equationstackrepresentation}, to the restriction functor. So $q_*$ corresponds to the locally analytic induction functor $\Ind_K^G(-)$ and $q_!$ corresponds to locally analytic compact induction functor $\cInd_K^G$ (see discussions after \cite[Corollary 3.3.7]{JC23}, cf. \cite[\S~1.4]{heyer20246functor}). We have a commutative diagram of Tate stacks
\begin{equation}\label{equationdiamgramUKXG}
\begin{tikzcd}
{[U/K^{\la}]}
\arrow[r,"j"]
\arrow[d,"p_U"']
& {[X/K^{\la}]}
\arrow[r,"i"]
& {[X/G^{\la}]}
\arrow[d,"p_X"]\\
{[*/K^{\la}]}
\arrow[r,equals]
& {[*/K^{\la}]}
\arrow[r,"q"']
& {[*/G^{\la}]}
\end{tikzcd}
\end{equation}

All morphisms in the above diagram are $!$-able, since $!$-ability is $!$-local on the target and we have Lemma~\ref{lemmaDcover}. 

\begin{lemma}\label{lemmaDcoverUX}
	The map $f: [U/K^{\la}]\to [X/G^{\la}]$ satisfies universal $D^*$-descent and universal $D^!$-descent
\end{lemma}
\begin{proof}
	We apply \cite[Lem.~4.7.1]{heyer20246functor}. We need to check $f$ is $D$-suave \cite[Def.~4.5.1]{heyer20246functor} and $f^*$ is conservative. By Lemma~\ref{lemmaDcover}, the map $X\to [X/G^{\la}]$ is a $D$-cover. In particular, it satisfies universal $*$-descent and universal $!$-descent and is conservative.
	
	By our assumption, the map $G\times^KU=\sqcup_{g\in G/K} g.U\simeq [U/K^{\la}]\times_{[X/G^{\la}]} X\to X$ is an open cover. Thus it is $D$-suave (even $D$-\'etale) and conservative by discussions in Example~\ref{exampleopencoverDcover}. Now the claim follows because the property for a map to be $D$-suave can be checked after pulling back to a universal $D^*$-cover on the target (\cite[Lemma 4.5.7]{heyer20246functor}).
\end{proof}

Let $\{[U/K^{\la}]^{[n]}\}_{[n]\in \Delta^{\rm op}}$ be the \v{C}ech nerve of $f$. For $\alpha=(1,\overline{g_1},\cdots,\overline{g_n})\in G\backslash(G/K)^{n+1}$, we put $U_{\alpha}:=U\cap g_1.U\cdots \cap g_n.U$, $K_\alpha:=K\cap g_{1}Kg_1^{-1}\cdots \cap g_nKg_n^{-1}$. Note that by our assumption of the covering $\{g.U\}_{g\in G/K}$, for each $n\in \bbN$, the set $\{\alpha\in G\backslash(G/K)^{n+1}: U_\alpha\neq \emptyset\}$ is a finite set. 
\begin{lemma}\label{lemmaCechnerveUK}
	We have an equivalence of Tate stacks
	\begin{equation*}
		\begin{aligned}
			{[U/K^{\la}]^{[n]}}&=[U/K^{\la}]\times_{[X/G^{\la}]}\times\cdots\times_{[X/G^{\la}]}[U/K^{\la}]\\
			&\simeq \sqcup_{\alpha\in G\backslash(G/K)^{n+1}} [U_{\alpha}/K_{\alpha}^{\la}],
		\end{aligned}
	\end{equation*}
	where \(U_{\alpha}\), \(K_{\alpha}\) are as above.
\end{lemma}
\begin{proof}
	Pulling back along $X\to [X/G^{\la}]$  it amounts to show 
	\[(G\times^KU)^{n/X}\simeq \sqcup_{\alpha\in G\backslash(G/K)^{n+1}} G\times^{K_{\alpha}}U_{\alpha}\] which is formal.
\end{proof} 
Following the language of Fargues, we define functors: 
\begin{equation*}
	\begin{aligned}
		\Res^{X/G}_{U/K}:=j^*i^*:D([X/G^{\la}])\rightarrow D([U/K^{\la}]),\\
		\cInd^{X/G}_{U/K}:=i_!j_!:D([U/K^{\la}])\rightarrow D([X/G^{\la}]).
	\end{aligned}
\end{equation*}
By $!$-descent, Lemma~\ref{lemmaDcoverUX} and Lemma~\ref{lemmaCechnerveUK}, we get a resolution
\begin{equation}\label{equationresolutionsheafC}
	\colim_{[n]\in \Delta^{\rm op}} C^n(\calF)\xrightarrow{\sim} \calF
\end{equation} 
of $\calF$ in $D([X/G^{\la}])$, where 
$$C^n(\calF):=\bigoplus_{\alpha \in G\backslash (G/K)^{n+1}} \cInd^{X/G}_{U_{\alpha}/K_{\alpha}}\Res^{X/G}_{U_{\alpha}/K_{\alpha}} \calF.$$

We can also view $\calF$ as a solid abelian sheaf, and consider the classical compactly supported cohomology as in \S~\ref{subseccompare}. For $\calF_{\alpha}$ a $K_\alpha$-equivariant vector bundle on $U_\alpha$, let
$$\Ind_{U_{\alpha}/K_{\alpha}}^{X/G}\calF_{\alpha}:= \prod_{g\in G/K_{\alpha}} g^*j_{U_{\alpha}*}\calF_{\alpha}\simeq \bigoplus  g^*j_{U_{\alpha}*}\calF_{\alpha}$$
be a $G$-equivariant solid abelian sheaf on $X$ (note that $j_{U_{\alpha}*}\calF_{\alpha}=Rj_{U_{\alpha}*}\calF_{\alpha}$ as in the proof of Proposition~\ref{propforDresolution}). There is a resolution 
\begin{equation}\label{equationresolutionsheafD}
	\calF\to\lim_{[n]\in \Delta}D^n(\calF)	
\end{equation}
of $\calF$ given by 
\[D^n(\calF)=\bigoplus_{\alpha \in G\backslash (G/K)^{n+1}} \Ind^{X/G}_{U_{\alpha}/K_{\alpha}} \calF|_{U_\alpha}.\] In other words, $\calF$ is quasi-isomorphic to the cohomological complex (see \cite[Thm.~8.3.8]{weibel1994introduction})
\begin{equation}\label{equationresolutionsheafDcomplex}
	D^{\bullet}(\calF):=[D^0(\calF)\rightarrow D^1(\calF)\rightarrow \cdots\to D^n(\calF)\to\cdots]
\end{equation} 
One can check the above complex is a resolution by restricting to $U$ using $G$-translations, together with the very standard statements \cite[\href{https://stacks.math.columbia.edu/tag/02FU}{Tag 02FU}]{stacks-project}. 
\begin{remark}
 	The differentials of the complex $D^{\bullet}(\calF)$ is given using the following: It is easy to check that 
	$$\Hom_{K_\alpha}(\calF|_{U_\alpha},\calF_{\alpha})\cong \Hom_G(\calF,\Ind_{U_{\alpha}/K_{\alpha}}^{X/G}\calF_{\alpha})$$ 
	where the $\Hom$ sets are taken in the category of equivariant solid abelian sheaves. Let $\alpha=(1,\overline{g_1},\cdots,\overline{g_{n}})\in G\backslash(G/K)^{n+1}$ and $\beta=(\alpha,\overline{g_{n+1}})\in G\backslash(G/K)^{n+2}$, then we have maps
 	 $$\begin{aligned}
 	 \Hom_G(\Ind_{U_\alpha/K_\alpha}^{X/G}\calF|_{U_\alpha},\Ind_{U_{\beta}/K_{\beta}}^{X/G} \calF|_{U_{\beta}})& =\Hom_{G}(\Ind_{U_\alpha/K_\alpha}^{X/G}\calF|_{U_\alpha},\Ind_{U_\alpha/K_\alpha}^{X/G}\Ind_{U_{\beta}/K_{\beta}}^{U_\alpha/K_\alpha} \calF|_{U_{\beta}})\\
 	 &\leftarrow \Hom_{K_\alpha}(\calF|_{U_\alpha},\Ind_{U_{\beta}/K_{\beta}}^{U_\alpha/K_\alpha} \calF|_{U_{\beta}})\\
 	 &=\Hom_{K_\beta}(\calF|_{U_\beta},\calF|_{U_\beta}).
 	 \end{aligned}$$
 \end{remark} 
 
 \begin{lemma}\label{lemmacompactinduction}
	We have natural isomorphisms
	\begin{enumerate}
		\item $R\Gamma_c(X,\cInd^{X/G}_{U/K}\calF)\simeq \cInd_K^GR\Gamma_c(U,\calF)\simeq \cInd_K^GH_c^d(U,\calF)[-d]$ and 
		\item $R\Gamma_c^{\it cl}(X,\Ind_{U/K}^{X/G} \calF|_{U})\simeq \cInd_{K}^{G} R\Gamma(U,\calF|_U)\simeq \cInd_{K}^{G} H^0(U,\calF)$
	\end{enumerate}
	in $D_{\square}(\calD^{\la}(G,E))$.
\end{lemma}
\begin{proof}
	(1) By Lemma~\ref{lemmapullbackcohomology} (proper base change) and \eqref{equationdiamgramUKXG}, $R\Gamma_{c}(X,\calF)\in \Rep^{\la}_{E_\square}(G)$ corresponds to $p_{X!}\calF$, we have  $R\Gamma_c(X,\cInd^{X/G}_{U/K}\calF)\simeq R\Gamma_c(X,i_!j_!\calF)\simeq p_{X!}i_!j_!\calF\simeq q_!p_{U!}\calF\simeq \cInd_K^GR\Gamma_c(U,\calF).$ The concentration of $R\Gamma_c(U,\calF)$ in degree $d$ is Proposition~\ref{propcompactsmith}.

	(2) By definition $\Ind_{U/K}^{X/G} \calF|_{U}=\prod_{g\in G/K} g^*j_{U*}j_U^*\calF$. By Proposition \ref{propforDresolution}, the underlying (complex of) module of  $R\Gamma_c^{\it cl}(X,\Ind_{U/K}^{X/G} \calF|_{U})$ is $\oplus_{g\in G/K} R\Gamma(U,\calF)$. By Tate's acyclicity, $R\Gamma(U,\calF)$ concentrates in degree $0$ and by Lemma~\ref{lemmaecontinuousactionrigid}, the $K$-action on it is locally analytic. Since both sides concentrate in a single degree, the isomorphism is $\calD^{\la}(G,E)$-equivalent.
\end{proof}

\begin{proposition}\label{propositionresolutions}
	We have natural quasi-isomorphisms of complexes in $\Rep^{\la}_{E_\square}(G)$:
	\begin{enumerate}
		\item($!$-resolution) \begin{equation}\label{resolution}
			\colim_{[n]\in\Delta^{\rm op}} (\oplus_{\alpha\in G\backslash (G/K)^{n+1}}\cInd_{K_\alpha}^G  H_c^d(U_\alpha,\calF)[-d])\simeq R\Gamma_c(X,\calF);
		\end{equation}	
		\item ($*$-resolution)
		\begin{equation}\label{dualresolution}
			R\Gamma_c(X,\calF)\simeq \lim_{[n]\in\Delta}(\oplus_{\alpha\in G\backslash (G/K)^{n+1}} \cInd_{K_\alpha}^G H^0(U_\alpha,\calF)).
			\end{equation}
	\end{enumerate} 
\end{proposition}

\begin{proof}
	(1) The result follows by taking $R\Gamma_c(X,-)$ on the resolution \eqref{equationresolutionsheafC} which commutes colimits and applying Lemma~\ref{lemmacompactinduction}(1).

	(2) Applying $R\Gamma^{\it cl}_c(X,-)$ to the resolution \eqref{equationresolutionsheafD}, since each term of \eqref{equationresolutionsheafD} is acyclic for the functor by Lemma \ref{lemmacompactinduction}(2), we see 
	\[R\Gamma_c^{\it cl}(X,\calF)\simeq [\cInd_{K}^GH^0(U,\calF)\rightarrow \oplus_{\alpha=(1,\ovg_1)\in G\backslash (G/K)^2}\cInd_{K_{\alpha}}^G H^0(U_{\alpha},\calF)\to \cdots].\]
	We conclude by applying Corollary \ref{corocoincide}.
\end{proof}
\begin{theorem}\label{maintheorem}
	Let $\calF$ be a $G$-equivariant vector bundle on $X$ which is partially proper and countable at infinity. Suppose that there exists an affinoid open subspace $U\subset X$ satisfying Assumption~\ref{assumptionaffionid}. Then we have a quasi-isomorphism of $\calD^{\la}(G,E)$-modules 
	$$R\Gamma_c(X,\calF^\vee\otimes\omega_X)\simeq \bbD_{\BZ}(R\Gamma_c(X,\calF)),$$
	where $\calF^{\vee}$ is the dual vector bundle and $\omega_X$ is the dualizing complex on $X$ which is a line bundle concentrating in degree $-d$.
\end{theorem}
\begin{proof}

Taking the Bernstein--Zelevinsky dual \eqref{definitionBZdual} on \eqref{resolution}, we get following quasi-isomorphisms of $\calD^{\la}(G,E)$-representations:

$$\begin{aligned}
\bbD_{\BZ}(R\Gamma_c(X,\calF))& \simeq \lim_{[n]\in\Delta}(\oplus_{\alpha\in G\backslash (G/K)^{n+1}} \bbD_{\BZ}(\cInd_{K_\alpha}^G R\Gamma_c(U_{\alpha},\calF))) \\
&\simeq \lim_{[n]\in\Delta}(\oplus_{\alpha\in G\backslash (G/K)^{n+1}} \cInd^{G}_{K_{\alpha}} R\Gamma_c(U_\alpha,\calF)^\vee)\\
&\simeq \lim_{[n]\in\Delta}(\oplus_{\alpha\in G\backslash (G/K)^{n+1}} \cInd^{G}_{K_{\alpha}} R\Gamma(U_\alpha,\calF^\vee\otimes \omega_X))\\
&\simeq R\Gamma_c(X,\calF^\vee\otimes \omega_X)
\end{aligned}$$
where the first quasi-isomorphism is because the direct sum in each degree is essentially finite, the second one comes from Proposition \ref{BZdual}, using that $H^d_c(U_\alpha,\calF)$ is a locally analytic $K$-representation whose underlying module is a Smith space (Proposition~\ref{propcompactsmith}), the third one is Serre duality (Theorem~\ref{theoremserreduality}) and the last one is \eqref{dualresolution}.
\end{proof}
\begin{remark}\label{remarkalternateproof}
	One can obtain Theorem~\ref{maintheorem} without using classical compactly supported cohomology. In other words, \eqref{dualresolution} can be established using the $6$-functor formalism as follows. By $*$-descent and Assumption~\ref{assumptionaffionid}, a quasi-isomorphism $\calF\simeq \lim_{[n]\in\Delta}D^{n}(\calF)$ as in \eqref{equationresolutionsheafDcomplex} exists in $D(X)$, where $D^{n}(\calF)$ is given by $\oplus_{\alpha \in G\backslash (G/K)^{n+1}} \Ind^{X/G}_{U_{\alpha}/K_{\alpha}} \calF|_{U_\alpha}\simeq \oplus_{\beta \in (G/K)^{n+1}} j_{\beta*}\calF|_{U_\beta}$ (for $\beta=(\ovg_1,\ovg_1,\cdots,\ovg_n)$, $U_{\beta}=\cap_{i=0}^n g_{i}.U$ and $j_{\beta}:U_{\beta}\hookrightarrow X$). Replacing the \v{C}ech complex by the alternating one $D^{\mathrm{alt},\bullet}(\calF)$ \cite[\href{https://stacks.math.columbia.edu/tag/01FG}{Tag 01FG}]{stacks-project} (cf.\ the proof of Proposition~\ref{propositionDrinfeldcompactpresentation}), we may assume that $D^{\mathrm{alt},n}(\calF)=0$ for $n$ large enough. Since $R\Gamma_c(X,-)$ preserves finite limits of the stable $\infty$-categories, we see $R\Gamma_c(X,\calF)\simeq\varprojlim_{[n]\in\Delta}R\Gamma_c(X,D^{\mathrm{alt},n}(\calF))$. On the other hand, for each affinoid open $V=\Spa(A,A^{\circ})\hookrightarrow X$ consider the factorization $V\stackrel{j_V}{\rightarrow}\overline{V}=\Spa(A,\bbZ_p)\stackrel{\ovi_V}{\rightarrow}X\xrightarrow{p_X}*$ (Remark~\ref{remarkuniversalcompletion}). It is easy to see that $\ovi_{V,*}j_{V*}j_V^*\ovi_V^*\calF\simeq \ovi_{V,*}\ovi_V^*\calF$ using that $\calF$ is a vector bundle. We can factor $\ovi$ as $\overline{V}\stackrel{\Gamma_{\ovi}}\longrightarrow \overline{V}\times X\stackrel{\mathrm{pr}}\longrightarrow X$. Since $\Gamma_{\ovi*}=\Gamma_{\ovi!}$, $\mathrm{pr}_*=\mathrm{pr}_!$, we have $\ovi_{V!}=\ovi_{V*}$. Hence $R\Gamma_c(X, (\ovi_Vj_V)_*\calF|_V)=p_{X!}\ovi_{V,*}j_{V*}j_V^*\ovi_V^*\calF\simeq p_{X!}\ovi_{V,!}\ovi_V^*\calF$. This can be used to replace Lemma~\ref{lemmacompactinduction}(2). 
\end{remark}

\subsection{Application to the Drinfeld spaces}\label{subsec:Drinfeldspaces}
We apply Theorem~\ref{maintheorem} to the Drinfeld spaces and discuss examples of equivariant vector bundles. 

Let $L$ be a finite extension of $\bbQ_p$ with a uniformizer $\varpi_L$ and $E\supset L$. Let $n\geq 2$, we let $X=\Omega_{\mathrm{Dr},E}=\bbP_E^{n-1,\mathrm{an}}\setminus\cup_{H\in \calH}H$ be the Drinfeld space over $E$ of dimension $n-1$ \eqref{equationDrinfeldspace}. The group $G:=\PGL_n(L)$ acts on $X\subset \bbP_E^{n-1,\mathrm{an}}$ naturally and continuously.
\begin{lemma}\label{lemmasequenceaffinoid}
There exists a sequence of affinoid open subspaces $U_{r_1}\subset U_{r_2}\subset \cdots U_{r_k}\subset \cdots ,k\in\bbN$ of $\Omega_{\mathrm{Dr},E}$ such that $\Omega_{\mathrm{Dr},E}=G.U_{r_k}$ for all $k\geq 1$ and satisfying the following:
\begin{enumerate}
	\item There exists an open compact subgroup $K\subset G$ stabilizing $U_{r_k}$ and there exists a finite subset $\calS\subset G/K$ such that the set $\{g\in G/K\mid U_{r_k}\cap g U_{r_k}\neq \emptyset\}$ is contained in $ \calS$ for all $k$.
	\item For any $k_1<k_2$, $U_{r_{k_1}}$ is relatively compact in $U_{r_{k_2}}$ in the sense of \cite[\S~9.6.2]{BGR84}. 
\end{enumerate}	
\end{lemma}
\begin{proof}
	We can assume $E=L$. There is a $G$-equivariant map 
	\begin{equation}\label{equationresiduemapDrinfledspace}
		\lambda:|\Omega_{\mathrm{Dr},E}| \rightarrow  |\calB\calT|
	\end{equation}
	by \cite[\S~6, A)]{drinfeld1974elliptic}. Here $|\calB\calT|$ can be identified with the set of norms on $L^{n}$ up to similarity, and $\lambda$ sends $(z_0,\cdots,z_{n-1})\in \Omega_{\mathrm{Dr},E}\subset \bbP^{n-1}_{L}$ to the norm $(a_0,\cdots,a_{n-1})\in L^{n}\mapsto |\sum_i z_ia_i|$ (with $|p|=p^{-1}$). There is a metric $\rho$ on $|\calB\calT|$ such that $\rho(x,y)\leq 1$ for any $x,y\in\sigma$ for any simplex $\sigma$ of $\calB\calT$. Take a vertex $v_0$ ($0$-simplex) on $\calB\calT$ and for $r\in \bbQ$, set $V_r=\{a\in |\calB\calT|\mid \rho(a,v_0)\leq r \}\subset |\calB\calT|$ and 
	\[U_r:=\{z\in \Omega_{\mathrm{Dr},E}\mid \rho(v_0,\lambda(z))\leq r\}=\lambda^{-1}(V_r).\]
	By \cite[Prop.~6.1]{drinfeld1974elliptic}, $U_r$ is an affinoid open subspace of $\Omega_{\mathrm{Dr},E}$ and if $r<r'$, then $U_r$ is relatively compact in $ U_{r'}$. If $r>1$, then for any simplex $\sigma\in\calB\calT$, there exists $g\in G$ such that $g.\sigma\in V_r\subset V_2$ since $G$ acts transitively on the vertices of $\calB\calT$. Thus $G.V_r=|\calB\calT|$ and consequently $G.U_r=\Omega_{\mathrm{Dr},E}$. Take an open compact subgroup $K$ fixing $v_0$ (and we may assume $K=\PGL_n(\calO_L)$). Since the $G$-action is isometric, $K$ stabilizes $V_r$ for all $r$. If $\rho(g.v_0,v_0)>2r$, then $g(V_r)\cap V_r=\emptyset$. Take $\calS=\{g\in G/K\mid \rho(g.v_0,v_0)\leq 4\}$. Then $\calS$ is a finite set, and $g(V_r)\cap V_r=\emptyset$ for all $g\notin \calS$ and $r<2$. Consequently $g(U_r)\cap U_r=\emptyset$ for all $g\notin \calS$ and $r<2$. Now taking a sequence $r_1<r_2<\cdots <r_k<\cdots $ in  $\bbQ\cap (1,2)$, we get the required sequence of affinoid subspaces $U_{r_k}$.
\end{proof}

\begin{corollary}\label{corBZdualDrinfeld}
	Let $\calF$ be a $G$-equivariant vector bundle on $\Omega_{\mathrm{Dr},E}$. Then there is a quasi-isomorphism of $\calD^{\la}(G,E)$-modules
	$$R\Gamma_c(\Omega_{\mathrm{Dr},E},\bbD_{\GS}(\calF))\simeq \bbD_{\BZ}(R\Gamma_c(\Omega_{\mathrm{Dr},E},\calF)).$$
\end{corollary}
\begin{proof}
	By Lemma~\ref{lemmasequenceaffinoid}, Assumption~\ref{assumptionaffionid} is satisfied for $\Omega_{\mathrm{Dr},E}$. we can apply Theorem~\ref{maintheorem}.
\end{proof}
We now give examples of equivariant vector bundles on $\Omega_{\mathrm{Dr},E}$. Let $\varpi_L$ be a uniformizer of $L$. Let $\breve{L}$ be the completion of the maximal unramified extension of $L$. Let $D_L$ be the division algebra over $L$ of invariant $1/n\in \bbQ/\bbZ$ with the maximal order $\calO_{D_L}$. Let $D_L^{\times},\calO_{D_L}^{\times}$ be the group of units. Consider the Drinfeld tower
\begin{equation}\label{equationDrinfeldtower}
	\cdots\to\breve{\calM}_m\rightarrow \cdots \rightarrow \breve{\calM}_1\rightarrow \breve{\calM}_0\rightarrow \Omega_{\mathrm{Dr},\breve{L}}
\end{equation}
where each $\breve{\calM}_m$ is a rigid analytic space over $\breve{L}$ and is a finite \'etale $(\calO_{D_L}/\varpi_L^m)^{\times}$-torsor over $\breve{\calM}_{0}$. Theses spaces are essentially considered in \cite[\S~3]{drinfeld1976coverings}, except that in the moduli description of $\breve{\calM}_m$, we do not require the quasi-isogeny of special formal $\calO_{D_L}$-modules is of height $0$. By \cite[Thm.~3.72]{rapoportZink1996}, there is an isomorphism $ \breve{\calM}_0\simeq \Omega_{\mathrm{Dr},\breve{L}}\times \bbZ$ compatible with the Weil descent data described in \textit{loc. cit.}.

The tower $\{\breve{\calM}_m\}_m$ admits a continuous action of $(D_L^{\times}\times \GL_n(L))/L^{\times}$ where $L^{\times}\subset D_L^{\times}\times \GL_n(L):g\mapsto (g,g^{-1})$ is in the center. The period map $\{\breve{\calM}_m\}_m \to \Omega_{\mathrm{Dr},\breve{L}}$ is $D_L^{\times}\times \GL_n(L)$-equivariant for the natural action of $\GL_n(L)$ and trivial action of $D_L^{\times}$-action on $\Omega_{\mathrm{Dr},\breve{L}}$ exhibiting the Drinfeld tower as a pro-\'etale $D_L^{\times}$-torsor over $\Omega_{\mathrm{Dr},\breve{L}}$.

The Weil descent datum on $\breve{\calM}_m$ is not effective, but that on $\varpi_L^{\mathbb Z}\backslash\breve{\calM}_m$ is \cite[Thm.~3.49]{rapoportZink1996}. Denote its descent to $L$ by $\Sigma_m$ and put $\Sigma_{m,E}:=\Sigma_m\otimes_{L}E$. See \cite[\S1.1 and \S3.1]{dospinescuLeBras2017revetements} for the $\GL_2(\bbQ_p)$-case. We write $\ovG:=\GL_n(L)/\varpi_L^{\bbZ}$ and $ \ovD_L^{\times}:=D_L^{\times}/\varpi_L^{\bbZ}$. Then the tower
\begin{equation}\label{equationDrinfeldtowerrational}
	\cdots\to\Sigma_{m,E}\rightarrow \cdots \rightarrow \Sigma_{1,E}\rightarrow \Sigma_{0,E}\rightarrow \Omega_{\mathrm{Dr},E}
\end{equation}
exhibits $ \{\Sigma_{m,E}\}_{m}$ as a pro-\'etale $\ovD_L^{\times}$-torsor over $\Omega_{\mathrm{Dr},E}$.
\begin{lemma}
	There exists an affinoid open subspace $U\subset \Sigma_{m,E}$ such that Assumption~\ref{assumptionaffionid} is satisfied for the $\ovG$-action on $\Sigma_{m,E}$. The same statement holds for the $\ovG$-action on $\Omega_{\mathrm{Dr},E}$.
\end{lemma}
\begin{proof}
	For $\Omega_{\mathrm{Dr},E}$, this is a direct consequence of Lemma~\ref{lemmasequenceaffinoid}. For $\Sigma_{m,E}$, the proof is similar to Lemma~\ref{lemmasequenceaffinoid}, using the building and the discussions in \cite[\S~3.71]{rapoportZink1996}, or take the preimage of the affinoids in $\Omega_{\mathrm{Dr},E}$. 
\end{proof}
Therefore, Theorem~\ref{maintheorem} applies to $\ovG$-equivariant vector bundles on $\Sigma_{m,E}$, e.g. $\calO_{\Sigma_{m,E}}$ and the sheaf of differentials $\Omega_{\Sigma_{m,E}}^1$. The following example will be used for our application to the $\GL_2(\bbQ_p)$-case in \S~\ref{subsectionBZdual}. Let $\pi_m:\Sigma_{m,E}\rightarrow \Omega_{\mathrm{Dr},E}$ be the period maps. Since $\pi_m$ is $\ovG$-equivariant, $\pi_{m,*}\calO_{\Sigma_{m,E}}$ is naturally a $\ovG$-equivariant sheaf on $\Omega_{\mathrm{Dr},E}$.
\begin{definition}\label{definitionFM}
	Let $\tau$ be an irreducible smooth representation of $\ovD_L^{\times}$ over $E$ and $\tau^{\vee}=\Hom_E(\tau,E)$ be the dual representation. We define a $\ovG$-equivariant $\calO_{\Omega_{\mathrm{Dr},E}}$-module 
	\[\calF_{\tau}:=\Hom_{D_L^{\times}}(\tau^{\vee},\varinjlim_m\pi_{m,*}\calO_{\Sigma_{m,E}})\]
	which is a vector bundles by Lemma~\ref{lemmaFMvectorbundle} below
\end{definition}
\begin{lemma}\label{lemmaFMvectorbundle}
	The equivariant sheaf $\calF_{\tau}$ in Definition~\ref{definitionFM} is a vector bundle on $\Omega_{\mathrm{Dr},E}$. Moreover, there is an isomorphism $\calF_{\tau}^{\vee}\simeq \calF_{\tau^{\vee}}$ of equivariant vector bundles.
\end{lemma}
\begin{proof}
	The group $\ovD_L^{\times}$ is compact, hence the action of $\ovD_L^{\times}$ on $\varinjlim_m\pi_{m,*}\calO_{\Sigma_{m,E}}$ is semisimple and $\ovD_L^{\times}$ acts on $\tau$ through a finite quotient $\ovD_L^{\times}/(1+\varpi_L^{k}\calO_{D_L})$ for some $k\geq 1$.  Note that the map $\pi_m$ is finite \'etale, since the maps $\breve{\calM}_m/\varpi_L^{\bbZ}\rightarrow \Omega_{\mathrm{Dr},\breve{E}}$ is and $\pi_m$ is its descent (using \cite[Lem.~1.1.5]{conrad1999irreducible}). Hence each $\pi_{m,*}\calO_{\Sigma_{m,E}}$ is a vector bundle on $\Omega_{\mathrm{Dr},E}$. Then $\calF_{\tau}$ is a direct summand of the $\ovG$-equivariant vector bundle $\pi_{k,*}\calO_{\Sigma_{k,E}}\simeq (\pi_{m,*}\calO_{\Sigma_{m,E}})^{(1+\varpi_L^k\calO_{D_L})},m\geq k$, thus is a vector bundle.

	The trace pairing \cite[\href{https://stacks.math.columbia.edu/tag/0BVH}{Tag 0BVH}]{stacks-project} for the finite \'etale map $\pi_k$ induces a $(\ovG\times \ovD_L^{\times})$-equivariant perfect pairing. Hence $\pi_{k,*}\calO_{\Sigma_{k,E}}\simeq (\pi_{k,*}\calO_{\Sigma_{k,E}})^{\vee}$ equivariantly and $\calF_{\tau}^{\vee}\simeq \calF_{\tau^{\vee}}$.
\end{proof}
By Lemma~\ref{lemmaFMvectorbundle} and Theorem~\ref{maintheorem}, we get \
	\begin{equation}\label{equationBZdualFM}
		\bbD_{\rm BZ}(H^{n-1}_c(\Omega_{\mathrm{Dr},E},\calF_{\tau}))\simeq H^{n-1}_c(\Omega_{\mathrm{Dr},E},\calF_{\tau^{\vee}}\otimes \Omega^{n-1}_{\Omega_{\mathrm{Dr},E}})[1-n]	
	\end{equation}
as $\ovG$-representations.
\begin{remark}\label{remarkothercentralcharacter}
	The $\ovG$ or $\GL_n(L)$ representations in \eqref{equationBZdualFM} should play a role in the $p$-adic Langlands correspondence and the $p$-adic Jacquet--Langlands correspondence, see \S~\ref{sec:duality-supercuspidal} for the example when $n=2,L=\bbQ_p$.  For application, we need remove the following restrictions on the above discussion. 
	\begin{enumerate}
		\item The Bernstein--Zelevinsky dual in \eqref{equationBZdualFM} is for $\ovG$-representations rather than $\GL_n(L)$-representations. This can be resolved as follows. Let $G=\GL_n(L)$. We have isomorphisms $\calD^{\la}(G,E)\otimes_{E[\varpi_L^{\bbZ}]}^{L\square}E\simeq \calD^{\la}(\ovG,E)$ and $R\Hom_{E[\varpi_L^{\bbZ}]}(E,\calC_c^{\la}(G,E))\simeq \calC_c^{\la}(\ovG,E)[-1]$. Hence
		\begin{equation}
			\begin{aligned}
				R\uHom_{\calD^{\la}(G,E)}(V,\calC_c^{\la}(G,E))&\simeq R\uHom_{\calD^{\la}(\ovG,E)}(V,R\uHom_{\calD^{\la}(G,E)}(\calD^{\la}(\ovG,E),\calC_c^{\la}(G,E)))\\
				&\simeq R\uHom_{\calD^{\la}(\ovG,E)}(V,R\uHom_{E[\varpi_L^{\bbZ}]}(E,\calC_c^{\la}(G,E)))\\
				&\simeq R\uHom_{\calD^{\la}(\ovG,E)}(V,\calC_c^{\la}(\ovG,E))[-1].
			\end{aligned}
		\end{equation}
		Similar result holds Replacing $\ovG$ by $\PGL_n(L)$ using that 
		\begin{align*}
			\calD^{\la}(G,E)\otimes_{\calD^{\la}(L^{\times},E)}^{L\square}E&\simeq \calD^{\la}(\PGL_n(L),E)\\
			R\Hom_{\calD^{\la}(L^{\times},E)}(E,\calC_c^{\la}(G,E))&\simeq \calC_c^{\la}(\PGL_n(L),E)[-1].	
		\end{align*}
		These isomorphisms can be obtained using the resolutions of the trivial representation in \cite[Thm.~5.8]{rodrigues2022solid}.
		\item One may want to consider $D_L^{\times}$-representations that do not factor through $\ovD_L^{\times}$, so that $\varpi_L^{\bbZ}$ acts non-trivially on $H^{n-1}_c(\Omega_{\mathrm{Dr},E},\calF_{\tau})$. In practice, this can be resolved by twisting a character $G=\GL_n(L)\xrightarrow{\det}L^{\times}\xrightarrow{\chi}E^{\times}$ (after enlarging $E$) and noticing that \[\bbD_{\rm BZ}(V\otimes \chi\circ \det)\simeq \bbD_{\rm BZ}(V)\otimes \chi^{-1}\circ \det.\]
		One can also consider the duality theory of categories of $G$-equivariant vector bundles and $G$-representations with fixed central characters, as in Appendix~\ref{sec:translation-duality}. 
	\end{enumerate}
\end{remark}
\section{Finiteness for locally analytic representations}\label{sectionfintiepresentation}
Let $G$ be a $p$-adic Lie group. We define a notion of ``finite presentation'' for locally analytic representations of $G$ (Definition~\ref{deffinitepresentation}) and prove finiteness result for certain $\Ext$ groups (Theorem~\ref{thm:finite-presentation-perfect}). Let $E$ be a finite extension of $\bbQ_p$.
\subsection{Representations of compact presentation}
We define locally analytic representations of compact presentation and show that compactly supported coherent cohomologies of the Drinfeld spaces give such representations. We also show that locally analytic representations of compact presentation are prim (Theorem~\ref{theoremprimobject}).
\begin{definition}\label{deffinitepresentation}
	Let $G$ be a $p$-adic Lie group. Let $V\in \Rep^{\la}_{E_\square}(G)$. A \emph{compact presentation} of $V$ consists of a sequence $V_{r_k}^{\bullet},k\in \bbN$ of complexes in $\Rep^{\la}_{E_\square}(G)^{\heartsuit}$ with transition maps $\gamma_{k}:V_{r_k}^{\bullet}\to V_{r_{k+1}}^{\bullet}$ between complexes:
	\[
	\begin{tikzcd}[column sep=tiny,row sep=small] && \vdots \arrow[d] & \vdots \arrow[d] & \\
	V_{r_k}^{\bullet}:\arrow[d,"\text{\rotatebox[origin=c]{90}{$\sim$}}","\gamma_k"']&\cdots \arrow[r] & \bigoplus_{j=1}^{l_i}\cInd_{K_{i,j}}^G W_{r_k,i,j} \arrow[r]\arrow[d,"\gamma_k^i"] & \bigoplus_{j=1}^{l_{i+1}}\cInd_{K_{i+1,j}}^G W_{r_k,i+1,j} \arrow[r]\arrow[d,"\gamma_k^{i+1}"] & \cdots \\
	V_{r_{k+1}}^{\bullet}:&\cdots \arrow[r] & \bigoplus_{j=1}^{l_i}\cInd_{K_{i,j}}^G W_{r_{k+1},i,j} \arrow[r]\arrow[d] & \bigoplus_{j=1}^{l_{i+1}}\cInd_{K_{i+1,j}}^G W_{r_{k+1},i+1,j} \arrow[r]\arrow[d] & \cdots \\
	&& \vdots & \vdots &
\end{tikzcd}
\]
such that there are quasi-isomorphisms $V\simeq V_{r_k}^{\bullet}$ for all $k$ satisfying the following properties:
	\begin{enumerate}
		\item The complexes $V_{r_k}^{\bullet}$ are uniformly bounded: there exist integers $ a\leq b$ such that $V_{r_k}^i=0$ for all $i\notin [a,b],k\in\bbN$.
		\item For each $i\in [a,b]$, there exists a finite set of compact open subgroups $K_{i,1},\cdots, K_{i, l_i}\subset G$ such that for each $k\in \bbN$, there exist locally analytic representations $W_{r_k,i,j},j=1,\cdots,l_i$ of $K_{i,j}$ with an isomorphism
		\begin{align}\label{equationtermscompactinduction}
			V_{r_k}^i\simeq \bigoplus_{j=1}^{l_i} \cInd_{K_{i,j}}^G W_{r_k,i,j}.
		\end{align}
		\item For all $k\in\bbN,i\in [a,b],j\in [1,l_i]$, the underlying solid vector spaces $W_{r_k,i,j}$ are Smith spaces.
		\item For all $k,i,j$, there are \emph{trace-class} morphisms $\gamma_{k,i,j}:W_{r_k,i,j}\rightarrow W_{r_{k+1},i,j}$ of $E_{\square}$-spaces, which are $K_{i,j}$-equivariant and induce $\gamma_{k}^i:V_{r_k}^i\rightarrow V_{r_{k+1}}^i$ via compact induction and \eqref{equationtermscompactinduction}.
		\item the morphisms $\gamma_{k}=\gamma_k^{\bullet}$ are quasi-isomorphism for all $k$.
	\end{enumerate}
	An object $V\in \Rep^{\la}_{E_\square}(G)$ is said to be \emph{of compact presentation} or \emph{compactly presented} if it admits a compact presentation as above.
\end{definition}
\begin{remark}
	The word ``compact'' signifies that the Smith spaces $W$ are compact objects in $\Mod(E_{\square})$, the maps $\gamma_{r_k,i,j}$ are dual to compact operators of Banach spaces. Moreover, the underlying $E$-space of a compactly presented representation is similar to a compact type space (being nuclear by Theorem~\ref{theoremprimobject}).
\end{remark}
Let $L$ be a finite extension of $\bbQ_p$ with ring of integers $\calO_L$ and let $\bbG$ be a reductive group over $L$ with $G=\bbG(L)$. We give examples of compactly presented $G$-representations.
\begin{example}\label{rem:smooth-finite-presentation}
 Every finitely generated smooth $E$-representation of $G$ is of finite
presentation, i.e. is quasi-isomorphic to a bounded complex with terms finite direct sums of representations of the form $\cInd_{K}^GW$ where $W$ is finite-dimensional and $K\subset G$ is open compact. For admissible representations, this follows from the
Schneider--Stuhler resolution
\cite[Prop.~II.2.1, Thm~II.3.1]
{schneiderStuhler1997representation} and for general see \cite[Cor.~1.4]{bartelsLuck2023algebraic}. Since finite-dimensional spaces are Smith and their identity maps are trace-class, finitely presented representations are of compact presentation.
\end{example}
\begin{example}
	The work of Kohlhaase--Schraen \cite{kohlhaase2012homological} (slightly extended in \cite{strauch2025bernsteinzelevinskydualitylocallyanalytic}) shows that locally analytic principal series representations of $G$ are of compact presentation if $\bbG$ is split. Fix a splitting model of $\bbG$ with a Borel subgroup and maximal torus $\bbB\supset \bbT$ over $\calO_L$. Let $T=\bbT(L), B=\bbB(L)$. Let $\chi:T\to E^{\times}$ be a continuous (locally analytic) character and $\Ind_B^G\chi$ a locally analytic principal series. Let $I\subset \bbG(\calO_L)$ be the Iwahori subgroup, the preimage of $ \bbB(\calO_L/\varpi_L)\subset \bbG(\calO_L/\varpi_L)$. Fix good analytic open normal subgroups $\cdots\subset I_{k+1}\subset I_k\subset\cdots\subset I$ as in \cite[Def.~4.1.3, Prop.~4.1.6]{emerton2006jacquet} which form a neighbourhood basis of the identity and admit rigid analytic Iwahori decompositions. For $k$ large enough such that $\chi$ is rigid analytic on $T\cap I_k$ \cite[\S~3.4]{emerton2017locally}, consider the Banach $I$-representation
	\[
		\calA_k:=(\Ind_B^G\chi(BI))^{I_k\text{-an}}\subset\Ind_B^G\chi
	\]
	of $I_k$-rigid analytic functions supported on $BI$. Then \cite[Thm.~2.5]{kohlhaase2012homological} shows that there is a resolution of the form
	\begin{align}\label{equationPSKoszul}
		V_{r_k}^\bullet:=\wedge^{\bullet}E^{d}\otimes_{E}\cInd_I^G \calA_k\xrightarrow{\sim}\Ind_B^G\chi
	\end{align} where $d$ is the rank of $T$ (see \cite[Prop.~5.2]{strauch2025bernsteinzelevinskydualitylocallyanalytic} for the solid version). For any $k$, the inclusion $\calA_k\rightarrow \calA_{k+1}$ is $E$-compact and induces a quasi-isomorphism $V_{r_k}^\bullet\xrightarrow{\sim} V_{r_{k+1}}^\bullet$. Hence the collection $\{V_{r_k}^\bullet\}_{k}$ satisfies the requirements in Definition \ref{deffinitepresentation}, except for that $\calA_k$ are Banach rather than Smith. After changing suitably the choice of $\calA_k$, e.g. \cite[Lem.~2.17]{strauch2025bernsteinzelevinskydualitylocallyanalytic}, we can take $ \calA_k$ Smith and the map $\calA_k\rightarrow \calA_{k+1}$ trace-class. Then we get a compact presentation of $\Ind_B^G\chi$. 
	
	It is also shown in \cite{strauch2025bernsteinzelevinskydualitylocallyanalytic} that representations constructed using Orlik--Strauch's functor $\calF_{B}^G$ are of compact presentation.
\end{example}
\begin{example}
	Agrawal--Strauch generalized the Schneider--Stuhler complexes \cite{schneiderStuhler1997representation} for locally analytic representations in \cite{agrawal2024resolutions}. For an admissible locally analytic representation $V$ of $G$, they constructed a bounded complex $\mathcal{S}_{\bullet}(V)$ \cite[Def.~4.13]{agrawal2024resolutions} whose terms are finte sums of representations of the form $\cInd_{P_{\sigma}^{\dagger}}^G\uHom_{U_{\sigma}}(\calD_r(U_{\sigma},E),V)$ \cite[(4.16)]{agrawal2024resolutions} where $U_{\sigma}\subset P_{\sigma}^{\dagger}$ are open compact (modulo center) subgroups of $G$ indexed by facets of the Bruhat--Tits building and $\calD_r(U_{\sigma},E),r\in (1/p,1)$ are the usual Banach distribution algebras of Schneider--Teitelbaum \cite[p.162]{ST03}. The $E$-spaces $\uHom_{U_{\sigma}}(\calD_r(U_{\sigma},E),V)$ are Smith \cite[Prop.~3.25, Prop.~3.20]{agrawal2024resolutions}. It's likely that the maps $\uHom_{U_{\sigma}}(\calD_r(U_{\sigma},E),V)\rightarrow \uHom_{U_{\sigma}}(\calD_{r'}(U_{\sigma},E),V)$ are trace-class. If the natural map $\calS_{\bullet}(V)\rightarrow V$ is a quasi-isomorphism, these resolutions varying $r$ would give a compact presentation of $V$. The last assumption already holds in some situations, see \cite[\S~4.18]{agrawal2024resolutions}.
\end{example}
Mikami gave the following example of a compactly presented representation that is not of finite presentation for the simplest non-compact group $\bbZ$  in \cite[Exam.~A.2]{mikami2026relativelydiscrete}.
\begin{example}
	Let $G=\bbZ$, then by the Cartier duality we have $\Rep^{\la}_{\bbQ_{p,\square}}(G)\simeq D_{\square}(\bbQ_p[T^{\pm}])$, cf. \cite[Thm.~4.4.4]{JC23}. Let $W=(\prod_{n\geq 0}\bbZ_p)[\frac{1}{p}]$ which is Smith and $\cInd_{\{0\}}^{\bbZ}W=\bbQ_p[T^{\pm}]\otimes_{\bbQ_p}^{\square}W$. Let $\gamma\in \End_{\bbQ_p}(W)$ be the map $ (a_0,a_1,a_2,\cdots)\mapsto (0,a_0,pa_1,\cdots)$ and 
	\[\alpha: \cInd_{\{0\}}^{\bbZ}W\to \cInd_{\{0\}}^{\bbZ}W,\quad	1\otimes a\mapsto 1\otimes a-T^{-1}\otimes \gamma(a)\] 
	for $a\in W$. Let $V:=\coker(\alpha)$. Then $V\simeq \colim (W\xrightarrow{\gamma}W\xrightarrow{\gamma}\cdots)$ where $T$ acts on $W$ by $\gamma$. The map $\gamma$ is trace-class. Hence the following commutative diagram
	\begin{equation*}
		\begin{tikzcd}0 \arrow[r] &\cInd_{\{0\}}^{\bbZ}W\arrow[r,"\alpha"] \arrow[d, "\cInd_{\{0\}}^{\bbZ}\gamma"] & \cInd_{\{0\}}^{\bbZ}W \arrow[r] \arrow[d, "\cInd_{\{0\}}^{\bbZ}\gamma"] & V \arrow[r] \arrow[d, "\simeq"',"T\times "] & 0 \\
			0 \arrow[r] & \cInd_{\{0\}}^{\bbZ}W\arrow[r,"\alpha"] &\cInd_{\{0\}}^{\bbZ}W \arrow[r] & V\arrow[r] & 0.
		\end{tikzcd}
	\end{equation*}
	shows that $V$ is of compact presentation. We note that since $\bbQ_p[T^{\pm}]$ is nuclear over $\bbQ_p$, $V\in D_{\square}(\bbQ_p[T^{\pm}])$ is dualizable if and only if it is compact and nuclear if and only if it is compact with its underlying $\bbQ_p$-vector space nuclear, cf. Theorem~\ref{theoremprimobject}.
\end{example}
As another example, we show that compactly supported coherent cohomologies of the Drinfeld spaces are of compact presentation. 
\begin{proposition}\label{propositionDrinfeldcompactpresentation}
	Let $G=\PGL_{n}(L)$ and let $\calF$ be a $G$-equivariant vector bundle on the Drinfeld space $\Omega_{\mathrm{Dr},E}$ of dimension $n-1$ as in \S~\ref{subsec:Drinfeldspaces}. Then the compactly supported cohomology $H^{n-1}_c(\Omega_{\mathrm{Dr},E},\calF)$ is a locally analytic representation of $G$ of compact presentation in the sense of Definition~\ref{deffinitepresentation}.
\end{proposition}
\begin{proof}
	Let $U_{r_k}\subset \Omega_{\mathrm{Dr},E}$ be the affinoid open subspaces in Lemma~\ref{lemmasequenceaffinoid} and $K$ the open compact subgroup of $G$ that stabilizes all $U_{r_k}$. We write $K_{\alpha}=g_{0}Kg_0^{-1}\cap\cdots \cap g_{n}Kg_n^{-1}$ as in \S\ref{subsectionresolutionFargues}. Take $\alpha\in G\backslash(G/K)^{n+1}$. By Proposition~\ref{propcompactsmith}, $R\Gamma_c(U_{\alpha,r_k},\calF)\simeq H_c^{n-1}(U_{\alpha,r_k},\calF)[1-n]$ and $H_c^{n-1}(U_{\alpha,r_k},\calF)$ is a Smith space representation of $E$. The resolution \eqref{resolution} shows that $R\Gamma_c(\Omega_{\mathrm{Dr},E},\calF)$ can be computed using the covering $\{g.U_{r_k}\}_{g\in G/K}$ and is quasi-isomorphic to the following (\v{C}ech) complex 
	\begin{equation}\label{equationresolutionrkinfinite}
		\begin{aligned}
			\cdots\rightarrow &\oplus_{\alpha=(1,g_1,\cdots,g_i)\in G\backslash (G/K)^{i+1}}\cInd_{K_{\alpha}}^G H_c^{n-1}(U_{r_k,\alpha},\calF)\rightarrow\cdots \\
			&\rightarrow \oplus_{(1,g)\in G\backslash (G/K)^2}\cInd_{K\cap g Kg^{-1}}^G H_c^{n-1}(g.U_{r_k}\cap U_{r_k},\calF)\rightarrow \cInd_{K}^G H_c^{n-1}(U_{r_k},\calF).
			\end{aligned}	
	\end{equation}

	For each $i\geq 0$, define a set
	\[\calT_i:=\{S\subset G/K\mid |S|=i+1\}.\]
	Formally, by the proof of \cite[\S3, Prop.~2]{SerreFAC}, the complex \eqref{equationresolutionrkinfinite} is quasi-isomorphic to the following alternating \v{C}ech complex (a direct summand of \eqref{equationresolutionrkinfinite})
	\begin{equation}
		\begin{aligned}\label{equationresolutionrknormalized}
			V_{r_k}^{\bullet}:\cdots\rightarrow &\oplus_{\alpha\in G\backslash \calT_i}\cInd_{K_{\alpha}^{\rm alt}}^G\left( H_c^{n-1}(U_{r_k,\alpha},\calF)\otimes_E\mathrm{sgn}_{\alpha}\right)\rightarrow\cdots \\
			&\rightarrow \oplus_{\alpha\in G\backslash\calT_1}\cInd_{K_{\alpha}^{\rm alt}}^G \left(H_c^{n-1}(U_{r_k,\alpha},\calF)\otimes_E\mathrm{sgn}_{\alpha}\right)\rightarrow \cInd_{K}^G H_c^{n-1}(U_{r_k},\calF).
			\end{aligned}	
	\end{equation}
	Here in \eqref{equationresolutionrknormalized}, for $\alpha=\{g_0,\cdots,g_{i}\}\in \calT_i$, a subset of $G/K$ of size $i+1$, we let \[K_{\alpha}^{\rm alt}:=\{g\in G\mid g.\alpha=\alpha\subset G/K\},\] 
	which is an open compact subgroup of $G$ containing $K_{\alpha}=\cap_{a=0}^{i}g_aKg_a^{-1}$ and stabilizing $\cap_{a=0}^{i}g_a.U_{r_k}$, and $\sgn_{\alpha}(g)=\sgn(\sigma)$ if $g.(g_0,\cdots,g_i)=(g_{\sigma(0)},\cdots,g_{\sigma(i)})$ for $g\in K_{\alpha}^{\rm alt}$ where $\sgn(\sigma)$ is the sign of the permutation. 
	
	There exists a finite subset $\calS\subset G/K$ such that the set $g.U_{r_k}\cap U_{r_k}=\emptyset$ if $ g\notin \calS$ for all $k$. For each $i$, every element in $G\setminus\{\alpha\subset G/K\mid U_{r_k,\alpha}\neq \emptyset,\alpha\in\calT_i\}$ has at least one representative in 
	\[\calT_i':=\{S\subset \calS\mid |S|=i+1\}.\]
	Choose representatives $\alpha_{i,1},\cdots,\alpha_{i,l_i}\in \calT_i'$ for elements in $G\backslash \calT_i$ if such representatives exist. For $i\in [-|\calS|-1,0],j\in[1,l_i]$, set $K_{i,j}=K_{\alpha_{-i,j}}^{\rm alt}$, and $W_{r_k,i,j}:=H_c^{n-1}(U_{r_k,\alpha_{-i,j}},\calF)\otimes\sgn_{\alpha_{-i,j}}$. Let $\gamma_{k,i,j}:W_{r_k,i,j}\rightarrow W_{r_{k+1},i,j}$ be the natural map induced by the inclusion $U_{r_k,\alpha_{i,j}}\hookrightarrow U_{r_{k+1},\alpha_{i,j}}$. Since $\Omega_{\mathrm{Dr},E}$ is separated, $U_{r_k,\alpha}$ are affinoid \cite[\S~9.6.1, Prop.~6]{BGR84}. Since $U_{r_k}$ is relatively compact in $U_{r_{k+1}}$, it follows easily from the definition in \cite[\S~9.6.2]{BGR84} that $U_{r_k,\alpha_{i,j}}$ is relatively compact in $U_{r_{k+1},\alpha_{i,j}}$. Hence $\gamma_{k,i,j}$ is trace-class by Corollary~\ref{corollarytraceclass}. We conclude that the resolutions $V_{r_k}^{\bullet}$ give a compact presentation of $H_c^{n-1}(\Omega_{\mathrm{Dr},E},\calF)$.
\end{proof}

\begin{lemma}\label{lem:Banachification-analytic}
Let $K$ be a compact $p$-adic Lie group and let $W$ be a
soid locally analytic representation whose underlying solid $E$-vector space is Smith. Then $W$ is $(h^+)$-analytic for some $h>0$ (see the proof for the meaning of $(h^+)$-analytic). The associated Banach space $W^B\rightarrow W$ (see \eqref{equationdefinitionassociatedBanach}) carries a unique locally analytic $K$-action for which $W^B\rightarrow W$ is
$\calD^{\la}(K,E)$-linear. 
\end{lemma}

\begin{proof}
We fix certain notation which will also be used later. We fix a uniform pro-$p$ subgroup of $K$ so that we have the distribution algebras $\calD^{(h)}(K,E), \calD^{(h^+)}(K,E)$  and $\calD_{(h)}(K,E)$ in \cite[\S~4.1.4]{rodrigues2022solid} for $h\in \bbQ_{>0}$. Here $\calD_{(h)}(K,E)$ is Banach, $\calD^{(h)}(K,E)$ is Smith and $\calD^{(h^+)}(K,E)=\varinjlim_{h'>h}\calD^{(h)}(K,E)=\varinjlim_{h'>h}\calD_{(h)}(K,E)$ is LS. For $h'>h$, the map $\calD^{(h')}(K,E)\rightarrow \calD^{(h)}(K,E)$ is trace-class and there are natural maps $\calD_{(h')}(K,E)\rightarrow \calD^{(h)}(K,E)\rightarrow \calD_{(h)}(K,E)$ (\cite[Cor.~4.18]{rodrigues2022solid}). A locally analytic representation $W$ of $K$ is $(h^+)$-analytic if the action of $\calD^{\la}(K,E)$ on $W$ extends to $\calD^{(h^+)}(K,E)$ and we write $W^{(h^+)-\mathrm{an}}$ for the space of $(h^+)$-analytic vectors of $W$, see \cite[Def.~4.29, Thm.~4.36]{rodrigues2022solid}. 

By \cite[Prop.~4.41]{rodrigues2022solid}, if $W$ is Smith, then $W$ is $(h^+)$-analytic for some $h$, hence is a $\calD_{(h)}(K,E)$-module for some $h$. Applying the trace-class functor \eqref{equationdefinitiontraceclassfunctor} to the action map 
\[\calD_{(h)}(K,E)\otimes_{E}^{\square}W\rightarrow W\] 
we obtain a map
\[
	\calD_{(h)}(K,E)^\trace\otimes_{E}^{\square}W^\trace\rightarrow (\calD(K^{(h)},E)\otimes_{E}^{\square}W)^\trace
 \rightarrow W^\trace
\]
Since $\calD_{(h)}(K,E)$ is Banach, $\calD_{(h)}(K,E)=\calD_{(h)}(K,E)^\trace$ and we obtain a map $\calD_{(h)}(K,E)\otimes_E^{\square}W^B\rightarrow W^B$. It is formal to check that this makes $W^B$ a $\calD_{(h)}(K,E)$-module, hence it is a locally analytic $K$-representation by \cite[Thm.~4.36]{rodrigues2022solid}. 
\end{proof}
Let $G$ be a $p$-adic Lie group. Let $\Rep^{\la}_{E_{\square}}(G)^{\omega}\subset \Rep^{\la}_{E_{\square}}(G)$ be the full subcategory consisting of compact objects \cite[\href{https://stacks.math.columbia.edu/tag/07LS}{Tag 07LS}]{stacks-project} (see also \cite[Prop.~1.4.4.1(3)]{lurie2017higher}).
\begin{lemma}\label{lem:compact-induction-Smith}
Let $K\subset G$ be a compact open subgroup.
\begin{enumerate}
	\item Let $W\in\Rep^{\la}_{E_{\square}}(K)^{\heartsuit}$ with underlying solid $E$-vector space Smith. Then $\cInd_K^G W$ is a compact object of $\Rep^{\la}_{E_{\square}}(G)$.
	\item There exists a set of compact generators of $\Rep^{\la}_{E_{\square}}(G)$ consisting of representations of the form $\cInd_K^G W$ where $W$ is Smith. Any object in $\Rep^{\la}_{E_{\square}}(G)^{\omega}$ is a retract of a finite complex of finite sums of these compact generators.
	\item If $V\in \Rep^{\la}_{E_{\square}}(G)^{\omega}$, then $\bbD_{\BZ}(V)$ is locally analytic and its underlying solid $E$-vector space is nuclear.
\end{enumerate}
\end{lemma}
\begin{proof}
(1) By \cite[Prop.~2.37(2)]{mikami2026finiteness}, the statement holds if $G=K$. The general case follows from the Frobenius reciprocity.

(2) The Grothendieck abelian category $\mathrm{Rep}^{\la}(G)^{\heartsuit}$ has a set of compact generators $\cInd_K^GW$ where $W$ are Smith spaces of the form $\calD^{(h)}(K,E)\otimes_{E_{\square}}E_{\square}[S], S\in \mathrm{ExtDisc}$, cf. the proof of \cite[Prop.~3.3.6]{JC23} where it is shown that $\mathrm{Rep}^{\la}(G)^{\heartsuit}$ is generated by $\cInd_K^GW^+$ for $W^+$ of the form 
	\[\calD^{(h^+)}(K,E)\otimes_{E_{\square}}E_{\square}[S]=\colim_{h'>h}\calD^{(h')}(K,E)\otimes_{E_{\square}}E_{\square}[S],\quad S\in \mathrm{ExtDisc}.\] 
	Indeed, by \cite[Thm.~4.36(1)]{rodrigues2022solid}, one has 
	\begin{align*}
		R\Hom_{\calD^{\la}(K,E)}(\calD^{(h)}(K,E)\otimes_{E_{\square}}E_{\square}[S],C)\simeq C^{R(h)-\an}(S)
	\end{align*}
	for any $C\in \Rep^{\la}_{E_{\square}}(G)$ and the functor $(-)^{R(h)-\an}$ preserves colimits. The last assertion follows from \cite[{\href{https://stacks.math.columbia.edu/tag/0948}{Tag 0948}}]{stacks-project}. 

	(3) This follows from (2) and Proposition~\ref{BZdual}: retracts, extensions and colimits of nuclear spaces are nuclear \cite[Prop.~5.33, Lem.~5.40]{andreychev2021pseudocoherent} and $\bbD_{\BZ}(\cInd_K^GW)=\cInd_K^GW^{\vee}$ is nuclear if $W$ is Smith.
\end{proof}
\begin{theorem}\label{theoremprimobject}
	Let $G$ be a $p$-adic Lie group. Let $V\in\Rep^{\la}_{E_{\square}}(G)$. We view $V$ as an object in $D([*/G^{\la}])$ via \eqref{equationstackrepresentation}. Let $f:[*/G^{\la}]\to *$ be the natural map.
	\begin{enumerate}
		\item The representation $V$ is $f$-prim \cite[Def.~4.4.1]{heyer20246functor} if and only if $V$ is in $\Rep^{\la}_{E_{\square}}(G)^{\omega}$ and the underlying solid $E$-space of $V$ is nuclear. Moreover, if $V$ is $f$-prim, the prim dual $\mathrm{PD}_f(V)$ corresponds to the Bernstein--Zelevinsky dual $\bbD_{\BZ}(V)$ and the natural morphism $V\to \bbD_{\BZ}(\bbD_{\BZ}(V))$ is an isomorphism.
		\item If $V$ is of compact presentation, then $V$ is $f$-prim.
	\end{enumerate}
\end{theorem}
\begin{proof}
	(1) We will use the criterion of prim objects in \cite[Lem.~4.4.6]{heyer20246functor}. Let $\pi_i,i=1,2$ be the projection $[*/(G\times G)^{\la}]\simeq [*/G^{\la}]\times [*/G^{\la}]\to [*/G^{\la}]$ to the $i$-th factor and let $\Delta:[*/G^{\la}]\to [*/(G\times G)^{\la}]$ be the diagonal map. We first show that for any $\Pi\in \Rep^{\la}_{E_{\square}}(G)$ with underlying solid $E$-space nuclear and $P\in \Rep^{\la}_{E_{\square}}(G)^{\omega}$ is compact, the natural morphism
	\begin{equation}\label{equationtheoremprimobject1}
		f_!(\pi_{2*}\uHom(\pi_1^*P,\Delta_!\mathbf{1})\otimes \Pi)\rightarrow f_*\uHom(P,\Pi)
	\end{equation}
	in \textit{loc. cit.} (see also \cite[Prop.~6.9]{scholze2025six}) is an isomorphism. 
	
	Let $K\subset G$ be an open compact subgroup. By (2) of Lemma~\ref{lem:compact-induction-Smith}, it suffices to prove that \eqref{equationtheoremprimobject1} is an isomorphism for $P$ of the form $\cInd_K^GW$ with $W$ Smith. Write $q:[*/K^{\la}]\to [*/G^{\la}]$ so that  $\cInd_K^GW$ corresponds to $q_!W\in D([*/G^{\la}])$ and $q^*$ corresponds to the restriction of $G$-representations to $K$. Since $\Delta_!:D([*/G^{\la}])\to D([*/(G\times G)^{\la}])$ is the compact induction functor as in \cite[\S~5.4]{heyer20246functor}, $\Delta_!\mathbf{1}$ corresponds to the representation $\calC_c^{\la}(G,E)$ and the object $\pi_{2*}\uHom(\pi_1^*P,\Delta_!\mathbf{1})$ corresponds to the $G$-representation $\bbD_{\BZ}(P)^{R\la}$, which is canonically isomorphic to $q_!W^{\vee}$ by Proposition~\ref{BZdual}. By the projection formula, we have $q_!W^{\vee}\otimes \Pi\simeq q_!(W^{\vee}\otimes q^*\Pi)$. By \cite[Prop.~2.39]{mikami2026finiteness}, we have $(f\circ q)_!=(f\circ q)_*$. Hence the lefthand side of \eqref{equationtheoremprimobject1} for $P=q_!W$ is given by $(f\circ q)_*( W^{\vee}\otimes \Pi|_K)$ calculated by $R\uHom_{\calD^{\la}(K,E)}(E,W^{\vee}\otimes_{E_{\square}}^L \Pi|_K)$. On the other hand, the righthand side of \eqref{equationtheoremprimobject1} is given by $R\uHom_{\calD^{\la}(G,E)}(E, R\uHom_{E}(P,\Pi)^{R\la})=R\uHom_{\calD^{\la}(G,E)}(P,\Pi)\simeq R\uHom_{\calD^{\la}(K,E)}(W,\Pi|_K)$. Then \eqref{equationtheoremprimobject1} becomes the natural isomorphism
	\[R\uHom_{\calD^{\la}(K,E)}(E,W^{\vee}\otimes_{E_{\square}}^L \Pi|_K)\simeq R\uHom_{\calD^{\la}(K,E)}(E,R\uHom_E(W, \Pi|_K))\]
	which follows from that $\Pi|_K$ is nuclear and $W$ is Smith. Hence \eqref{equationtheoremprimobject1} is an isomorphism when $P$ is compact and $\Pi$ is nuclear.

	Now \cite[Lem.~4.4.6]{heyer20246functor} says that an object $P\in D([*/G^{\la}])$ is $f$-prim if and only if the morphism \eqref{equationtheoremprimobject1} is an isomorphism after applying $\Hom(\mathbf{1},-)$. Hence $P$ is $f$-prim if $P$ is compact and nuclear by taking $\Pi=P$. The same lemma implies that $\mathrm{PD}_f(P)=\bbD_{\BZ}(P)^{R\la}=\bbD_{\BZ}(P)$ in this case.

	Conversely, suppose that $V$ is $f$-prim. Then \textit{loc. cit.} tells that $\mathrm{PD}_f(V)=\bbD_{\BZ}(V)^{R\la}$ and \cite[Lem.~4.4.18]{heyer20246functor} implies that $V$ is compact. Lemma~\ref{lem:compact-induction-Smith}(3) then implies that $\mathrm{PD}_f(V)$ is nuclear. By \cite[Lem.~4.4.4]{heyer20246functor}, $V\simeq \mathrm{PD}_f(\mathrm{PD}_f(V))$ and $\mathrm{PD}_f(V)$ is $f$-prim. Hence $V$ is nuclear.

	(2) Suppose that $V$ admits a compact presentation $\{V_{r_k}^{\bullet}\}$. Then $V$ is compact by Lemma~\ref{lem:compact-induction-Smith}. Consider the trace-class functor $(-)^{\trace}$ in \cite[Def.~A.13]{bosco2023rationalpadicHodge}. The quasi-isomorphism $\gamma_k$ factors through a morphism $V_{r_k}^{\bullet}\rightarrow V_{r_{k+1}}^{\bullet,\trace}\rightarrow V_{r_{k+1}}^{\bullet}$. Hence $V$ is a retract of $V^{\trace}$, so $V$ is nuclear and $V\simeq V^{\trace}$. We conclude that $V$ is $f$-prim by (1).
\end{proof}
\begin{corollary}
	In the situation of Theorem~\ref{theoremprimobject}, if $G$ is compact, then $V\in\Rep^{\la}_{E_{\square}}(G)$ is $f$-prim if and only if its underlying solid $E$-vector space is perfect. 
\end{corollary}
\begin{proof}
	By \cite[Prop.~2.37]{mikami2026finiteness}, $V$ is compact if and only if the underlying $E$-space is compact. The statement then follows from that $V\in D(E_{\square})$ is nuclear and compact if and only if $V$ is dualizable if and only if $V$ is perfect \cite[Prop.~5.37, Cor.~5.51.1]{andreychev2021pseudocoherent}. See also \cite[Cor.~6.2.52]{camargo2026locallyanalyticnotes}.
\end{proof}
\begin{remark}
	By Lemma~\ref{lem:Banachification-analytic}, if $V$ admits a compact presentation, then $V\simeq V_{r_k}^{\bullet,\trace}$ is a bounded complex of finite sums of representations of the form $\cInd_K^GW$ where $W$ is Banach. However, we do not know whether $\bbD_{\BZ}(V)$ is of compact presentation, even though the natural map $V\to\bbD_{\BZ}(\bbD_{\BZ}(V))$ is an equivalence by the theorem above. 
\end{remark}
\subsection{Finiteness of $\Ext$ groups}
We prove the main finiteness result for locally analytic representations of compact presentation (Theorem~\ref{thm:finite-presentation-perfect}).
\begin{lemma}\label{lemmahomsmithadmissible}
Let $K$ be an open compact subgroup of $G$ and let $f:W_{r_0}\to W_{r_1}$ be a $K$-equivariant morphism of locally analytic representations whose underlying solid $E$-vector spaces are Smith. We assume that $f$ is trace-class.  If $\Pi\in \Rep_{E_{\square}}^{\la}(G)^{\heartsuit}$ is an admissible locally analytic representation of $K$, then the induced morphism
\[
 R\underline{\Hom}_{\calD^{\la}(K,E)}(W_{r_1},\Pi)
 \longrightarrow
 R\underline{\Hom}_{\calD^{\la}(K,E)}(W_{r_0},\Pi).
\]
in $D(E_{\square})$ is represented by a morphism of bounded complexes of
Banach spaces which is degreewise compact.
\end{lemma}

\begin{proof}
	We follow the notation in the proof of Lemma~\ref{lem:Banachification-analytic}. Let 
	\[\Pi^{R(h^+)-\mathrm{an}}=R\underline{\Hom}_{\calD^{\la}(K,E)}(\calD^{(h^+)}(K,E),\Pi)\] be the derived $h^+$-analytic vectors of $\Pi$ (\cite[Def. 4.29, Thm. 4.36]{rodrigues2022solid}). We claim that $\Pi^{R(h^+)-\mathrm{an}}\simeq \Pi^{(h^+)-\mathrm{an}}$ (cf. \cite[Prop.~3.27]{agrawal2024resolutions}). Since $\Pi$ is admissible, $\Pi=R\uHom_E(\Pi^{\vee},E)^{R\la}$ by \cite[Lem. 4.2.6]{JC23}, and as in the proof of \textit{loc. cit.} (and also in \cite[Prop. 4.48]{rodrigues2022solid}), 
	\begin{align*}
		\Pi^{R(h^+)-\mathrm{an}}
		&\simeq R\underline{\Hom}_{\calD^{\la}(K,E)}(\calD^{(h^+)}(K,E),R\underline{\Hom}_{E}(\Pi^{\vee},E))\\
		&\simeq R\underline{\Hom}_{E}(\calD^{(h^+)}(K,E)\otimes_{\calD^{\la}(K,E)}^{L\square}\Pi^{\vee},E)\\
		&\simeq \underline{\Hom}_{E}(\calD^{(h^+)}(K,E)\otimes_{\calD^{\la}(K,E)}^{\square}\Pi^{\vee},E)\\
		&\simeq \Pi^{(h^+)-\mathrm{an}}
	\end{align*}
	where for the third equality, we used \cite[Cor. 4.2.4]{JC23} and that $\calD^{(h^+)}(K,E)\otimes_{\calD^{\la}(K,E)}^{\square}\Pi^{\vee}$ is isomorphic to a finite complex of LS spaces of compact type (see discussions below).
	
	Hence
	\[R\underline{\Hom}_{\calD^{\la}(K,E)}(W_{r_k},\Pi)=R\underline{\Hom}_{\calD^{\la}(K,E)}(W_{r_k},\Pi^{R(h^+)-\mathrm{an}})=R\underline{\Hom}_{\calD^{\la}(K,E)}(W_{r_k},\Pi^{(h^+)-\mathrm{an}}).\] 
	Choose $h'>h$, then $\uHom_E(\calD_{(h')}(K,E)\otimes_{\calD(K,E)}^{\square}\Pi^{\vee},E)$ is $h^+$-analytic. Since $\Pi^{\vee}$ is coadmissible over $\calD^{\la}(K,E)$, $\calD_{(h')}(K,E)\otimes_{\calD(K,E)}^{\square}\Pi^{\vee}$ is finitely generated over $\calD_{(h')}(K,E)$. By \cite[Thm.~8.9]{ST03}, the distribution algebra $\calD_{(h')}(K,E)(*)$ is Auslander-regular with finite global projective dimension. Thus we may choose a resolution of $\calD_{(h')}(K,E)\otimes_{\calD(K,E)}^{\square}\Pi^{\vee}$ by finite projective $\calD_{(h')}(K,E)$-modules $M_{(h'),d}\rightarrow \cdots \rightarrow M_{(h'),1}\rightarrow M_{(h'),0}$. Tensoring with the LS algebras $\calD^{(h^+)}(K,E)$ \cite[Thm.~4.45]{rodrigues2022solid}, we obtain a resolution of $\Pi^{(h^+)-\mathrm{an},\vee}$ by $M^{(h^+)}_{\bullet}=\calD^{(h^+)}(K,E)\otimes_{\calD_{(h')}(K,E)} M_{(h'),\bullet}$. Taking the $E$-dual which is exact for these LS spaces \cite[Thm. 3.40]{rodrigues2022solid}, we see $\Pi^{(h^+)-\mathrm{an}}$ is isomorphic to a complex $M^{(h^+),\vee}_{\bullet}$ where each $M^{(h^+),\vee}_{i}$ is a direct summand of 
	\[\calD^{(h^+)}(K,E)^{\vee,\oplus n_i}\simeq \calC^{(h^+)-\mathrm{an}}(K,E)^{\oplus n_i}.\]
	
	Let $h$ be large enough such that $W_{r_k}, k=0,1$ are $(h^+)$-analytic. As in the proof of Proposition~\ref{BZdual}, $R\underline{\Hom}_{\calD^{\la}(K,E)}(W_{r_k}, M_{i}^{(h^+),\vee})$ is a direct summand of 
	\[R\underline{\Hom}_{\calD^{\la}(K,E)}(W_{r_k}, \calC^{(h^+)-\mathrm{an}}(K,E)^{\oplus n_i})\simeq W_{r_k}^{\vee, \oplus n_i}.\] 
	Hence 
	\[R\underline{\Hom}_{\calD^{\la}(K,E)}(W_{r_k}, M_{i}^{(h^+),\vee})=\underline{\Hom}_{\calD^{\la}(K,E)}(W_{r_k}, M_{i}^{(h^+),\vee})\] 
	and the induced map $\underline{\Hom}_{\calD^{\la}(K,E)}(W_{r_1},M_{i}^{(h^+),\vee})\rightarrow \underline{\Hom}_{\calD^{\la}(K,E)}(W_{r_0},M_{i}^{(h^+),\vee})$ is a compact morphism between Banach spaces (as a direct summand of the compact morphism $W_{r_1}^{\vee,\oplus n_i}\rightarrow W_{r_0}^{\vee,\oplus n_i}$). The result follows.
\end{proof}
The following lemma is a derived variant of Proposition~\ref{propfinite}. 
\begin{lemma}\label{lem:strict-filtered-compact}
Let $C_{k,i},k,i=0,1,\cdots,N$ be objects in $D(E_{\square})$ with morphisms $g_{k,i}:C_{k,i}\rightarrow C_{k,i+1}$ and $\gamma_{k,i}:C_{k,i}\rightarrow C_{k+1,i}$. Assume that $C_{k,0}=0$ and write $C_k:=C_{k,N}$ for all $k$. Let $\mathrm{gr}_{k,i}:=\Cone(g_{k,i})\in D(E_{\square})$ for $i\leq N-1$. We assume that for each $k,i$, there is a homotopy $g_{k+1,i}\circ\gamma_{k,i}
\simeq \gamma_{k,i+1}\circ g_{k,i}$ such that the induced map $\gamma_{k,i}:\mathrm{gr}_{k,i}\rightarrow \mathrm{gr}_{k+1,i}$ in $D(E_{\square})$ can be represented by a degreewise compact morphism between bounded complexes of Banach $E$-spaces.

If the composition $\gamma=\gamma_{N-1,N}\circ\cdots\circ\gamma_{1,N}\circ \gamma_{0,N}: C_{0}\rightarrow C_{N}$ is a quasi-isomorphism, then $C_0$ is a perfect object in $D(E_{\square})$.
\end{lemma}
\begin{proof}
An object in $D(E_{\square})$ is perfect if and only if $C$ is dualizable if and only if $C$ is nuclear and compact by 
\cite[Prop.~9.3]{scholze2022complex} and
\cite[Lem.~5.46, Cor.~5.51.1]{andreychev2021pseudocoherent}. Since Banach spaces are nuclear \cite[Prop.~3.29]{rodrigues2022solid} and by our assumption $\mathrm{gr}_{k,i}$ is isomorphic to a bounded complex of Banach spaces, it is nuclear. It follows from \cite[Lem.~5.40]{andreychev2021pseudocoherent} and an induction on $i$ using the triangles $C_{k,i}\rightarrow C_{k,i+1}\rightarrow \mathrm{gr}_{k,i}\rightarrow$ that $C_{k,i}$ are nuclear for all $k,i$. We remain to show that $C_0$ is compact.

For the rest of the proof, we write $D(E_{\square})$ for its homotopy category which is a triangulated category. Write $D(E_{\square})^{\omega}$ for the full subcategory of compact objects in $D(E_{\square})$ consisting of objects that are isomorphic to a retract of a bounded complex of Smith spaces \cite[Lem.~5.36]{andreychev2021pseudocoherent}. Then $D(E_{\square})^{\omega}$ is a triangulated subcategory of $D(E_{\square})$ \cite[Lem.~5.40]{andreychev2021pseudocoherent} and is thick, thus we can consider the triangulated quotient category $Q(E_{\square}):=D(E_{\square})/D(E_{\square})^{\omega}$ \cite[\href{https://stacks.math.columbia.edu/tag/05RA}{Tag 05RA}]{stacks-project}. Let 
\[P:D(E_{\square})\rightarrow Q(E_{\square})\] 
be the quotient functor so that $P(C)\simeq 0$ if $C\in D(E_{\square})^{\omega}$ \cite[\href{https://stacks.math.columbia.edu/tag/05RK}{Tag 05RK}]{stacks-project}. 

If $f:G_1\rightarrow G_2$ is a degreewise compact morphism between bounded complexes of Banach spaces, then $f^{\vee}: G_2^{\vee}\rightarrow G_1^{\vee}$ is a degreewise trace-class morphism between bounded complexes of Smith spaces (applying the exact dual in \cite[Thm.~3.40(1)]{rodrigues2022solid}) which factors through $G_1^{\vee,\mathrm{tr}}$ \cite[Lem.~A.14]{bosco2023rationalpadicHodge}. Here $G_1^{\vee,\mathrm{tr}}$ is obtained by applying the trace-class functor $(-)^{\trace}$ \cite[Def.~A.13]{bosco2023rationalpadicHodge} to $G_1^{\vee}$ and is a complex of Banach spaces, see the proof of \cite[Lem.~A.58]{bosco2021drinfeldproetale}. Hence the map $f$ factors through $G_1^{\vee,\mathrm{tr},\vee}\in D(E_{\square})^{\omega}$ a bounded complex of Smith spaces. Therefore, $P(f):P(G_1)\rightarrow P(G_2)$ is a zero morphism.

Thus, by our assumption, the maps $P(\gamma_{k,i}):P(\mathrm{gr}_{k,i})\rightarrow P(\mathrm{gr}_{k+1,i})$ are zero for all $k,i$. Since $P$ is exact, the triangle $P(C_{k+1,i})\rightarrow P(C_{k+1,i+1})\rightarrow P(\mathrm{gr}_{k+1,i})\rightarrow $ is exact. Hence the map $P(\gamma_{k,i+1}):P(C_{k,i+1})\rightarrow P(C_{k+1,i+1})$ factors through a map $P(C_{k,i+1})\rightarrow P(C_{k+1,i})$ (see the proof of \cite[\href{https://stacks.math.columbia.edu/tag/0149}{Tag 0149}]{stacks-project}). Consequently, the composition $P(\gamma)=P(\gamma_{N-1,N})\circ\cdots\circ P(\gamma_{1,N})\circ P(\gamma_{0,N}):P(C_0)\rightarrow P(C_N)$ factors through a map $P(C_0)=P(C_{0,N})\rightarrow P(C_{1,N-1})\rightarrow \cdots \rightarrow P(C_{N-1,1})\rightarrow P(C_{N,0})=0$. In particular, $P(\gamma)$ is a zero morphism. Since $\gamma$ is an isomorphism by our assumption, so is $P(\gamma)$. This forces $P(C_0)=0$. Hence $C_0\in D(E_{\square})^{\omega}$ by \cite[\href{https://stacks.math.columbia.edu/tag/05RK}{Tag 05RK}]{stacks-project}, i.e., $C_0$ is compact. The result follows.
\end{proof}
\begin{theorem}\label{thm:finite-presentation-perfect}
	Let $V$ be a locally analytic representation of $G$ of compact presentation in the sense of Definition \ref{deffinitepresentation},
	and let $\Pi\in \Rep^{\la}(G)^{\heartsuit}$ be an admissible locally analytic representation of $G$. Then $R\underline{\Hom}_{\calD^{\la}(G,E)}(V,\Pi)$ is a perfect object in $D(E_\square)$. In particular,
	\[\dim_E\Ext^i_{\calD^{\la}(G,E)}(V,\Pi)<\infty\] 
	for every $i\geq0$.
\end{theorem}
\begin{proof}
Take a compact presentation 
$V_{r_k}^\bullet,k\in\bbN$ of $V$ as in Definition~\ref{deffinitepresentation} and write $\sigma_{\leq i}V_{r_k}^\bullet$ for the stupid truncation of $V_{r_k}^{\bullet}$ so that there are exact sequences $0\rightarrow V_{r_k}^i[-i]\rightarrow \sigma_{\leq i}V_{r_k}^{\bullet}\rightarrow \sigma_{\leq i-1}V_{r_k}^{\bullet}\rightarrow 0$ of complexes \cite[\href{https://stacks.math.columbia.edu/tag/0118}{Tag 0118}]{stacks-project}. Set 
\begin{equation}
	\begin{aligned}
		C_{k,i}&:=R\underline{\Hom}_{\calD^{\la}(G,E)}(\sigma_{\leq i}V_{r_k}^\bullet,\Pi),\\
		\mathrm{gr}_{k,i}&:=R\underline{\Hom}_{\calD^{\la}(G,E)}(V_{r_k}^{i+1}[-i-1],\Pi)
	\end{aligned}
\end{equation}
so that $\mathrm{gr}_{k,i}\simeq \Cone(C_{k,i}\rightarrow C_{k,i+1})$. Denote again $\gamma_k: C_{k+1,i}\rightarrow C_{k,i}$ for the natural map induced by the truncation of the map $\gamma_k: V_{r_k}\rightarrow V_{r_{k+1}}$ and the functor $R\underline{\Hom}_{\calD^{\la}(G,E)}(-,\Pi)$. Then $\gamma_k$ induces a map $\mathrm{gr}_{k+1,i}\rightarrow \mathrm{gr}_{k,i}$. 

By the assumptions in Definition~\ref{deffinitepresentation}, $V_{k}^i\simeq \oplus_{j=1}^{n_i}\cInd_{K_{i,j}}^{G}W_{r_k,i,j}$ for some Smith space representations of open compact subgroups $K_{i,j}$ of $G$. By Frobenius reciprocity, we have 
\[R\underline{\Hom}_{\calD^{\la}(G,E)}(V_{k}^i,\Pi)\simeq \oplus_{j=1}^{n_i}R\underline{\Hom}_{\calD^{\la}(K_{i,j},E)}(W_{r_k,i,j},\Pi).\]
The map $\gamma_k^i:V_{k}^i\rightarrow V_{k+1}^i$ is induced by maps $\gamma_{k,i,j}:W_{r_k,i,j}\rightarrow W_{r_{k+1},i,j}$. Thus the map $\mathrm{gr}_{k+1,i}\rightarrow \mathrm{gr}_{k,i}$ is the direct sum of maps 
\[R\underline{\Hom}_{\calD^{\la}(K_{i,j},E)}(W_{r_{k+1},i,j},\Pi)\xrightarrow{R\underline{\Hom}_{\calD^{\la}(K_{i,j},E)}(\gamma_{k,i,j},\Pi)} R\underline{\Hom}_{\calD^{\la}(K_{i,j},E)}(W_{r_{k},i,j},\Pi).\] 
By Lemma \ref{lemmahomsmithadmissible}, the above maps can be represented by degreewise compact morphisms between bounded complexes of Banach spaces. So the same holds for $\mathrm{gr}_{k+1,i}\rightarrow \mathrm{gr}_{k,i}$. Since there exists only finitely many $i$ such that $V_{r_k}^i\neq 0$, $C_{k,i}\simeq R\underline{\Hom}_{\calD^{\la}(G,E)}(V,\Pi)$ and $\gamma_{k}: C_{k+1,i}\rightarrow C_{k,i}$ are quasi-isomorphisms for $i$ large enough and for all $k$, and $C_{k,i}=0$ for $i$ small enough. After picking finitely many indices $k, i$ and renumbering them suitably, the assumptions in Lemma \ref{lem:strict-filtered-compact} are satisfied for the collection $\{C_{k,i}\}_{k,i}$. The result follows.
\end{proof}
\begin{remark}
	We explain the homological algebra issue underlying the proof of Theorem~\ref{thm:finite-presentation-perfect} and Lemma~\ref{lem:strict-filtered-compact} related to Definition~\ref{deffinitepresentation}. Let $V\in \Rep^{\la}_{E_\square}(G)$ with a compact presentation $V_{r_k}^{\bullet}$. We take as example the case where $V_{r_k}^{i}=0$ if $i\neq 0,-1$. Applying $R\underline{\Hom}_{\calD^{\la}(G,E)}(-,\Pi)$ to $V_{r_k}^{\bullet}$ produces a commutative diagram in $D(E_{\square})$: 
\begin{equation}\label{equationhomologicalcommutative}
	\begin{tikzcd}
		0\ar[r] & \mathrm{gr}_{1,0} \ar[r,"f_{1,0}"] \ar[d, "\gamma_{0,0}"'] & \mathrm{gr}_{1,-1} \ar[r] \ar[d, "\gamma_{0,-1}"'] & 0  \\
		0 \ar[r] & \mathrm{gr}_{0,0}\ar[r,"f_{0,0}"] & \mathrm{gr}_{0,-1} \ar[r] & 0
		\end{tikzcd}
\end{equation}
where $\mathrm{gr}_{k,i}=R\underline{\Hom}_{\calD^{\la}(G,E)}(V_{r_k}^i,\Pi)$. Moreover, $R\underline{\Hom}_{\calD^{\la}(G,E)}(V,\Pi)\simeq \Cone(f_{i,0})[-1]$. By Lemma~\ref{lemmahomsmithadmissible} the map $\gamma_{0,i}:\gr_{1,i}\rightarrow \gr_{0,i}$ can be represented by a degreewise compact morphism between bounded complexes of Banach spaces. If the whole diagram \eqref{equationhomologicalcommutative} could be represented by a commutative diagram of bounded chain complexes of Banach spaces with vertical maps degreewise compact, then the quasi-isomorphism $\gamma_0: \Cone(f_{1,0})\rightarrow \Cone(f_{0,0})$ would satisfy the assumption in Proposition~\ref{propfinite} and we may conclude that $R\underline{\Hom}_{\calD^{\la}(G,E)}(V,\Pi)$ is perfect. 

However, a priori, we may only assume that \eqref{equationhomologicalcommutative} is represented by a diagram of bounded complexes of Banach spaces which is commutative up to homotopy. Suppose $\gr_{k,i}\simeq \gr_{k,i}^{\bullet}$ are bounded complexes of Banach spaces and the arrows in \eqref{equationhomologicalcommutative} are chain maps with vertical maps degreewise compact. Then there exists a homotopy $h:\gr_{1,0}^{\bullet}\rightarrow \gr_{0,-1}^{\bullet-1}$ such that $\gamma_{0,-1}\circ f_{1,0}-f_{0,0}\circ \gamma_{0,0}=d\circ h+h\circ d$ so that the quasi-isomorphism $\Cone(f_{1,0})\xrightarrow{\sim} \Cone(f_{0,0})$ is represented by a chain map $\gamma_0^{\bullet}$ with $\gamma_0^{i}=\begin{pmatrix}
\gamma_{0,-1}^{i}&h^{i+1}\\&\gamma_{0,0}^{i+1}
\end{pmatrix}:\gr_{1,-1}^{i}\oplus \gr_{1,0}^{i+1}\rightarrow \gr_{0,-1}^{i}\oplus \gr_{0,0}^{i+1}$ in degree $i$ \cite[\href{https://stacks.math.columbia.edu/tag/014F}{Tag 014F}]{stacks-project}. Since $h^i$ may not be compact, we cannot apply Proposition~\ref{propfinite} directly for the chain map $\gamma_0^{\bullet}$. The homotopical issue can be resolved using the third resolution $V_{r_2}^{\bullet}$. Then one obtains similarly $f_{2,0}:\gr_{2,0}\rightarrow \gr_{2,-1}$, a homotopy $h':\gr_{2,0}^{\bullet}\rightarrow \gr_{1,-1}^{\bullet-1}$ and a chain map $\gamma_1^{\bullet}:\Cone(f_{2,0})\rightarrow \Cone(f_{1,0})$ using $h'$. We observe that the composition $\gamma_0^{\bullet}\circ\gamma_1^{\bullet}:\Cone(f_{2,0})\rightarrow \Cone(f_{0,0})$ satisfies the assumption of Proposition~\ref{propfinite} since the composition of a compact morphism with other morphisms of Banach spaces is again a compact morphism, hence $R\underline{\Hom}_{\calD^{\la}(G,E)}(V,\Pi)$ is perfect.
\end{remark}
\begin{remark}\label{remarktworesolution}
	It should be true that the existence of two resolutions $V_{r_0}^{\bullet}$ and $V_{r_1}^{\bullet}$ (instead of infinitely many) of $V$ in Definition~\ref{deffinitepresentation} is sufficient to deduce the same conclusion of Theorem~\ref{thm:finite-presentation-perfect}. In fact, with a more careful argument, \eqref{equationhomologicalcommutative} can be represented by a strictly commutative diagram of bounded complexes of Banach spaces with vertical maps degreewise compact in the specific situation. Moreover, $V$ is already prim with the two resolutions by the proof of Theorem~\ref{theoremprimobject}. In all examples we know, representations of $G$ with two resolutions $V_{r_k}^{\bullet},k=0,1$ in Definition~\ref{deffinitepresentation} are of compact presentation. 
\end{remark}
\section{Dualities of supercuspidal representations of
\texorpdfstring{$\GL_2(\bbQ_p)$}{GL2(Qp)}}\label{sec:duality-supercuspidal}
In this section, we first recall some locally analytic representations of $\GL_2(\mathbb{Q}_p)$ which admits geometric realization using the Drinfeld tower. Then we compute the Schneider--Teitelbaum dual and the Bernstein--Zelevinsky dual of these representations. Finally, we give some applications by combining these two duality results: we discuss the resulting Serre-type duality and the trace pairing on the related extension groups.

Let $E/\bbQ_p$ be a sufficiently large finite extension and write $G=\PGL_2(\bbQ_p)$. Let $\bbQ_p^{\mathrm{ur}}$ be the maximal unramified extension of $\bbQ_p$ and let $\breve{\mathbb Q}_p$ be its completion. We normalize the local reciprocity map so that $p\in\bbQ_p^{\times}$ is sent to a geometric Frobenius element. Our convention is that the Hodge--Tate weight of the $p$-adic cyclotomic character $\epsilon:\Gal(\overline{\mathbb{Q}}_p/\bbQ_p)\to\bbQ_p^\times$ is $1$. If no further comment is made, locally analytic representations are solid. If $V$ is a solid $E$-space, we write $V^{\vee}:=\uHom_E(V,E)$ for its (underived) $E$-linear dual and $V(*)$ for its underlying compactly generated topological space.

\subsection{The supercuspidal locally analytic representations}\label{subsec:supercuspidal-representations}
We recall the definition of some locally analytic representations which can be constructed using the de Rham complex of the Drinfeld tower in \cite{dospinescuLeBras2017revetements}.  

Let $D_p$ be the unique 
quaternion division algebra over $\bbQ_p$, and let $D_p^\times$ be its group of units. Let $\rho:\Gal(\overline{\mathbb{Q}}_p/\bbQ_p)\to \GL_2(E)$ be a $2$-dimensional
de Rham representation, with Hodge--Tate weights $0,1$ and $\det\rho\cong\epsilon$. We let $M$ be the Fontaine's rank two $(\varphi,N,\Gal(\overline{\mathbb{Q}}_p/\bbQ_p))$-module $D_{\rm pst}(\rho\otimes \epsilon^{-1})$ over $E\otimes_{\bbQ_p}\bbQ_p^{\mathrm{ur}}$ attached to $\rho$, and write $\WD(\rho)=\WD(M)$ for the Weil--Deligne representation associated to $M$ (see \cite[\S~4]{breuilSchneider2007}). We assume that $\WD(\rho)$ is absolutely irreducible. 

Let $\Pi(\rho)$ be the admissible unitary Banach representation of $\GL_2(\bbQ_p)$ attached to $\rho$ under the $p$-adic local Langlands correspondence by Colmez \cite{colmez2010representations}, and write $\pi(\rho):=\Pi(\rho)^{\la}$ for its subspace of locally analytic vectors; see also \cite[\S0.2]{colmezDospinescuNiziol2020}. Our normalization follows \cite{colmez2010representations,colmezDospinescuNiziol2020}: the central character of $\Pi(\rho)$ and $\pi(\rho)$ is $\det\rho\otimes\epsilon^{-1}$ which is trivial. Hence the action of $\GL_2(\bbQ_p)$ on $\pi(\rho)$ factors through $G$. 

The results of Dospinescu--Le Bras in \cite[Thm.~1.2 \& Thm.~1.4]{dospinescuLeBras2017revetements} give a geometric realization of
$\pi(\rho)$ using the de Rham complex of the Drinfeld tower.  Let $\pi_{\sm}$ be the smooth supercuspidal representation of $\GL_2(\bbQ_p)$ attached to $\WD(\rho)$ by the local Langlands correspondence, and let $\tau_{\sm}$ be its Jacquet--Langlands transfer to $D_p^\times$. By our assumption, the central characters of $\pi_{\sm}$ and $\tau_{\sm}$ are trivial.

Consider the Drinfeld tower $\{\Sigma_{m,E}\}_{m\in\bbN}$ in \eqref{equationDrinfeldtowerrational} of dimension $1$ for $\bbQ_p$, see also \cite[\S~3.1]{dospinescuLeBras2017revetements}. Define
\begin{align*}
 \calO[M]&:=\uHom_{E[D_p^\times]}(\tau_{\sm},
 \varinjlim_m H^0(\Sigma_{m,E},\calO_{\Sigma_{m,E}})),\\
 \Omega^1[M]&:=\uHom_{E[D_p^\times]}(\tau_{\sm},
 \varinjlim_m H^0(\Sigma_{m,E},\Omega^1_{\Sigma_{m,E}}))
\end{align*}
which are Fr\'echet $\calD^{\la}(G,E)$-modules. 
\begin{definition}\label{definitionpippic}
	We define solid $\calD^{\la}(G,E)$-modules attached to $M$:
	\begin{equation}
		\pi_c:=\calO[M]^{\vee},\quad\widetilde{\pi}:=\Omega^1[M]^{\vee},
	\end{equation}
	which are admissible locally analytic representations by \cite[Thm.~1.2, Rem.~1.3(b)]{dospinescuLeBras2017revetements}.
\end{definition}
Put $M_{\dR}:= (\overline{\mathbb{Q}}_p\otimes_{\bbQ_p^{\mathrm{ur}}}M)^{\Gal(\overline{\mathbb{Q}}_p/\bbQ_p)}\simeq D_{\rm dR}(\rho\otimes \epsilon^{-1})$. This is a two-dimensional $E$-vector space. The de Rham differential induces a strict exact sequence
\begin{equation}\label{eq:universal-sequencedual}
	0\rightarrow\calO[M]\xrightarrow{d}\Omega^1[M]
 \rightarrow M_{\dR}\otimes_E\pi_{\sm}^{\vee}\rightarrow0.
\end{equation}
See \cite[Thm.~1.4]{dospinescuLeBras2017revetements} and
\cite[\S~2.1.5]{colmezDospinescuNiziol2020}. Taking the exact dual between Fr\'echet spaces and LS spaces (\cite[Thm.~3.40]{rodrigues2022solid}), one has a short exact sequence:
\begin{equation}\label{eq:universal-sequence} 0\longrightarrow \pi_{\sm}^{\oplus 2}
 \longrightarrow\widetilde\pi\longrightarrow\pi_c\longrightarrow0
\end{equation}
of admissible locally analytic representations of $G$. The representation $\pi(\rho)$ is related to the exact sequence \eqref{eq:universal-sequence} as follows. 	We fix an isomorphism $D_{\rm dR}(\epsilon^{-1})\simeq E$ so that we can identify $M_{\rm dR}\simeq D_{\rm dR}(\rho\otimes_E\epsilon^{-1})\simeq D_{\rm dR}(\rho)\otimes_ED_{\rm dR}(\epsilon^{-1})\simeq D_{\rm dR}(\rho)$.
\begin{theorem}[{\cite[Thm.~1.2 \& Thm.~1.4]{dospinescuLeBras2017revetements}}]\label{theoremdospinescuLeBras}
	Let $\calL=\Fil^{0}D_{\rm dR}(\rho)\subset D_{\dR}(\rho)$ be the $1$-dimensional Hodge filtration, viewed as an $E$-line of $M_{\rm dR}$ under the isomorphism $M_{\rm dR}\simeq D_{\rm dR}(\rho)$. Let $ M_{\rm dR}^{\vee}\rightarrow \calL^{\vee}$ be the quotient map. 	Then the following diagram is a pushout:
	\begin{equation}\label{equationpushoutdospinescuLeBras}
		\begin{tikzcd}0 \arrow[r] & M_{\rm dR}^{\vee}\otimes_E\pi_{\sm} \arrow[r] \arrow[d, two heads] & \widetilde{\pi} \arrow[r] \arrow[d, two heads] & \pi_c \arrow[r] \arrow[d, equal] & 0 \\
			0 \arrow[r] & \calL^{\vee}\otimes_E\pi_{\sm} \arrow[r] & \pi(\rho) \arrow[r] & \pi_c \arrow[r] & 0.
		\end{tikzcd}
	\end{equation}
\end{theorem}
\begin{proof}
	For the $E$-line $\calL\subset M_{\rm dR}$, let $W_{M,\calL}\subset\Omega^1[M]$ be the inverse image of $\calL\otimes_E\pi_{\sm}^{\vee}\subset M_{\dR}\otimes_E\pi_{\sm}^{\vee}$. By \cite[Prop.~2.9]{colmezDospinescuNiziol2020} (essentially \cite[Thm.~1.4]{dospinescuLeBras2017revetements}), there is a $\GL_2(\bbQ_p)$-equivariant isomorphism $(W_{M,\calL})^{\vee}\cong\pi(\rho)$. The result follows by taking $E$-duals.
\end{proof}
Hence $\pi(\rho)$ is an extension of $\pi_c$ by $\pi_{\rm sm}$. For later use, we also need to consider other such extensions, classified in \cite{ding2022Ext1}. 
\begin{definition}\label{definitionpiL}
	For each $E$-line $\calL\subset M_{\rm dR}$, we denote by $\pi_{\calL}$ be the pushout of the map $M_{\rm dR}^{\vee}\otimes_E\pi_{\sm}\rightarrow \widetilde{\pi}$ along the map $M_{\rm dR}^{\vee}\otimes_E\pi_{\sm}\rightarrow \calL^{\vee}\otimes_E\pi_{\sm}$ as in \eqref{equationpushoutdospinescuLeBras}. 
\end{definition}
In particular, Theorem~\ref{theoremdospinescuLeBras} gives $\pi(\rho)\cong \pi_{\mathcal{L}}$ with $\mathcal{L}=\Fil^{0} D_{\dR}(\rho)$.	
\begin{theorem}[{\cite[Thm.~2.5]{ding2022Ext1}}]\label{theoremDingext1}
	There is a natural isomorphism 
	\[M_{\dR}\cong \Ext^1_{\calD^{\la}(G,E)}(\pi_c,\pi_{\sm})\] 
	such that an $E$-line $\calL\subset M_{\dR}$ is sent to the class of the extension $\pi_{\calL}$. Moreover, $\pi_{\calL}\cong\pi_{\calL'}$ if and only if $\calL=\calL'$.
\end{theorem}
\begin{proof}
	The result in \emph{loc.\ cit.} computed the extension group in the category of admissible $\GL_2(\bbQ_p)$-representations. Since there is no $\GL_2(\bbQ_p)$-equivariant map from $\pi_c$ to $\pi_{\sm}$ by, for example, \cite[Prop.~2.5.2]{su2025locally}, any such extension must have trivial central character. Hence the extension group coincides with the group of $G$-extensions, which also coincides with the extension group in the category of locally analytic representations of $G$ over $E$ by Proposition~\ref{prop:classical-solid-Hom-Ext} below.  For the last assertion, see for example \cite[Prop.~2.6.3]{su2025locally}. See also \cite[Thm.~2.5.7 \& Thm.~3.2.2 \& Rem.~3.2.4]{su2025locally} for more discussions.
\end{proof}
\begin{remark}\label{rem:general-integral-weights}
The representations above have analogues for arbitrary regular integral
Hodge--Tate weights, using the corresponding equivariant vector bundles
and dual BGG complexes (\cite[Rem.~1.3(a) and Rem.~11.12]{dospinescuLeBras2017revetements},
\cite[\S2.2]{su2025locally}), or using Colmez's construction (\cite{colmez2018poids}, \cite[Thm.~1.1]{ding2023change}). Similar discussions hold for these representations and the compatibility of the duality
functors in this section with translation functors and fixing central characters is discussed in Appendix~\ref{sec:translation-duality}, in particular
Proposition~\ref{prop:fixed-central-translation}.
\end{remark}

For later use, we also record here a compatibility result between the classical and solid extension groups of locally analytic representations. In particular, when discussing $\Hom$ or $\Ext^1$ spaces, we can freely use results in classical settings and will omit further discussion of this point in the sequel.
\begin{proposition}\label{prop:classical-solid-Hom-Ext}
Let $V,W$ be admissible locally analytic representations. Taking the underlying topological spaces induces natural isomorphisms of $E$-vector spaces:
\[
\begin{aligned}
	\Hom_{\calD^{\la}(G,E)}(V,W)&\xrightarrow{\sim} \Hom_{G}(V(*),W(*))\\
	\Ext^1_{\calD^{\la}(G,E)}(V,W)&\xrightarrow{\sim} \Ext^1_{G}(V(*),W(*))
\end{aligned}
\]
where the groups on the right are computed in the abelian category of
classical admissible locally analytic representations.
\end{proposition}

\begin{proof}
The map on $\Hom$ is an isomorphism by \cite[Prop.~3.7]{rodrigues2022solid} since the underlying topological spaces of  admissible locally analytic representations are of compact type. Hence the abelian category of classical admissible locally analytic representations is embedded fully faithfully in $\Rep^{\la}_{E_\square}(G)^{\heartsuit}$ and we write the category as $\Rep^{\la}_{E_\square}(G)^{\heartsuit,\mathrm{adm}} $.

The $\Ext^1$-groups in the abelian categories can be computed by Yoneda extensions \cite[\href{https://stacks.math.columbia.edu/tag/06XU}{Tag 06XU}]{stacks-project}. By the fully faithful embedding, it suffices to show if there is an extension 
\[0\rightarrow W\rightarrow T\rightarrow V\rightarrow 0\]
in $\Rep^{\la}_{E_\square}(G)^{\heartsuit} $ such that if $W,V\in \Rep^{\la}_{E_\square}(G)^{\heartsuit,\mathrm{adm}} $, then $T\in \Rep^{\la}_{E_\square}(G)^{\heartsuit,\mathrm{adm}}$. Admissibility of $G$ depends only on the restriction of $T$ to some compact open subgroup $H\subset G$. Since $V$ is of compact type, one has $R\uHom_{E}(V,E)=V^{\vee}$ \cite[Lem.~3.27]{rodrigues2022solid}. Hence the dual sequence 
\[0\rightarrow V^{\vee}\rightarrow T^{\vee}\rightarrow W^{\vee}\rightarrow 0\]
is exact. If we can show that $T^{\vee}$ arises from a classical coadmissible $\calD^{\la}(H,E)(*)$-module, then taking the exact dual in \cite[Thm.~3.40]{rodrigues2022solid}, we see $T^{\vee\vee}$ is an extension of $V$ by $W$ and we can conclude that $T\xrightarrow{\sim}T^{\vee\vee}$ by the five lemma. Then $T$ is admissible.

To show $T^{\vee}$ is coadmissible, write $\calD^{\la}(H,E)=\varprojlim_n \calD_{n}(H,E)$ as an inverse limit of Noetherian Banach $E$-algebras along (classically) flat and injective maps as in \cite{ST03}. Finitely presented $\calD_n(H,E)(*)$ are Banach modules with the canonical topology. It's then easy to see that the classical and solid finitely presented $\calD_n(H,E)$-modules are equivalent and finitely presented modules in $\Mod_{\square}(\calD_n(H,E))$ are stable under extensions. Moreover, by the flatness results in \cite[Rem.~3.2]{ST03} and arguments in \cite[Prop.~5.33]{bode2026sixoperations} using Bar complexes, one can show that $\calD_n(H,E)\otimes_{\calD^{\la}(H,E)}^{L\square}W^{\vee}\simeq \calD_n(H,E)\otimes_{\calD^{\la}(H,E)}^{\square}W^{\vee}$ is finitely presented over $\calD_n(H,E)$. We can then argue in the same way as in \cite[Lem.~2.1.1]{BreuilExt1} to conclude that $T^{\vee}\simeq \varprojlim_n \calD_n(H,E)\otimes_{\calD^{\la}(H,E)}^{\square}T^{\vee}$ and each $\calD_n(H,E)\otimes_{\calD^{\la}(H,E)}^{\square}T^{\vee}$ is finitely presented. The result follows.
\end{proof}

\subsection{Schneider--Teitelbaum duality}\label{subsectionSTdual}
In this subsection, we compute the Schneider--Teitelbaum dual of $\pi(\rho)$ and representations $\pi_c,\widetilde{\pi}$ in Definition~\ref{definitionpippic} (Proposition~\ref{proposition:STdualpic}). The results below are known to experts but we can not find in the literature. And our methods seem to be new.

For $V\in\Rep^{\la}_{E_\square}(G)$, the Schneider--Teitelbaum dual of $V$ is defined as 
\[\bbD_{\ST}(V):=R\uHom(V,E)^{R\la}\] 
where $(-)^{R\la}$ denotes the derived locally analytic vectors with respect to the action of $G$ \cite[Def.~4.40]{rodrigues2022solid}; see \cite[Prop.~4.2.7]{JC23} for the compact case. We collect some basic properties of this functor in Appendix~\ref{sec:STdual}, especially its relation with Schneider--Teitelbaum's original definition \cite{ST05}: For a coadmissible $\calD^{\la}(G,E)$-module $M$, set
\[\bbD(M):=R\uHom_{\calD^{\la}(G,E)}(M,\overline{\calD}^{\la}(G,E)[\dim G]),\]
see Definition~\ref{defSTdual} for the notation and note that the character $\mathfrak d_G$ is trivial for out $G$ by \cite[Rem.~3.4]{ST05}. Proposition~\ref{prop:ST-linear-dual} gives an isomorphism 
\[ \bbD(V^{\vee})\simeq R\uHom_E(\bbD_{\ST}(V),E).\]
\begin{remark}\label{remarSTdualindependent}
	By the definition, the representation $\bbD_{\rm ST}(V)$ is independent of whether we consider $V$ as a $\PGL_2(\bbQ_p)$-representation or a $\GL_2(\bbQ_p)$-representation.
\end{remark}
We first determine the dual of some topologically irreducible representations. Let $\Delta$ be a $(\varphi,\Gamma)$-module of rank two over the Robba ring $\calR_E$ for $\bbQ_p$ with the coefficients in $E$. We assume that $\Delta$ is indecomposable and let $\pi(\Delta)$ be the associated admissible locally analytic representation of $\GL_2(\bbQ_p)$ in \cite[Thm.~0.1]{colmez2016representations}. Let $\Delta^{\vee}$ be the dual $(\varphi,\Gamma)$-module and write $\check\Delta:=\Delta^{\vee}(\epsilon)=\Delta^{\vee}\otimes_{\calR_E}D_{\rm rig}(\epsilon)$. 
\begin{proposition}\label{propositionSTselfdual}
	Suppose that $\pi(\Delta)(*)$ (resp. $\pi(\check\Delta)(*)$) is topologically irreducible and $\pi(\Delta)^{\vee}(*)$ (resp. $\pi(\check\Delta)^{\vee}(*)$) has (Gelfand--Kirillov) dimension $1$ (\cite[\S~8]{ST03}), then 
	\[\bbD_{\ST}(\pi(\Delta))\simeq \pi(\Delta^{\vee}(\epsilon))[-1]\]
	in $\Rep^{\la}_{E_\square}(\GL_2(\bbQ_p))$.
\end{proposition}
\begin{proof}
	Let $H\subset \GL_2(\bbQ_p)$ be an open compact subgroup. The dimension $1$ assumption on $\pi(\Delta)^{\vee}(*)$ implis that 
	\begin{align*}
		\bbD(\pi(\Delta)^{\vee}(*))&:=R\Hom_{\calD^{\la}(\GL_2(\bbQ_p),E)(*)}( \pi(\Delta)^{\vee}(*), \overline{\calD}^{\la}(\GL_2(\bbQ_p),E)(*))[4]\\
		&\simeq R\Hom_{\calD^{\la}(H,E)(*)}( \pi(\Delta)^{\vee}(*), \overline{\calD}^{\la}(H,E)(*))[4]	
	\end{align*}
	concentrates in degree $-1$ and $0$. Since $\pi(\Delta)(*)$ has no dimension $0$ quotient, it is pure in the sense of \cite[Prop.~8.7]{ST03}. By the argument of \cite[Lem.~6.2]{huWang2024quaternion} we find $\bbD(\pi(\Delta)^{\vee}(*))$ concentrates in degree $-1$. By \cite[Lem.~8.4]{ST03}, \cite[Rem.~4.2.8]{JC23}, $\bbD(\pi(\Delta)^{\vee}(*))$ or $\bbD(\pi(\Delta)^{\vee})$ are both inverse limits of extension groups over Banach algebras. One can compare solid and classical $\Ext$-groups between finitely presented modules over Noetherian Banach algebras as in \cite[\S~5.4]{bode2026sixoperations}. Consequently, we can deduce that $R\uHom_E(\bbD_{\ST}(\pi(\Delta)),E)\simeq \bbD(\pi(\Delta)^{\vee})\simeq \underline{\bbD(\pi(\Delta)^{\vee})(*)}$ concentrates in a single degree. Since $\bbD_{\ST}(\pi(\Delta))$ is admissible in the sense of \cite[Def.~4.2.1]{JC23}, $H^i(\bbD_{\ST}(\pi(\Delta)))$ are admissible, thus are acyclic for $R\uHom_{E}(-,E) $\cite[Lem.~3.27]{rodrigues2022solid}. Then the spectral sequence for hypercohomology gives $H^i(R\uHom_E(\bbD_{\ST}(\pi(\Delta)),E))\simeq H^{-i}(\bbD_{\ST}(\pi(\Delta)))^{\vee}$. We deduce that $\bbD_{\ST}(\pi(\Delta))$ concentrates in degree $1$. Similarly for $\bbD_{\ST}(\pi(\check\Delta))$. 

	We need to construct a non-zero morphism $\pi(\Delta)\rightarrow \bbD_{\ST}(\pi(\check\Delta))[1]$. Colmez constructed an extension 
	\begin{align}\label{equationDboxtimesP1}
		0\rightarrow \pi(\check\Delta)^{\vee}(*)\rightarrow \Delta\boxtimes_{\omega} \bbP^1\rightarrow \pi(\Delta)(*)\rightarrow  0
	\end{align}
	of topological $\calD^{\la}(\GL_2(\bbQ_p),E)(*)$-modules in \cite[Thm.~0.1]{colmez2016representations}. We claim that the extension \eqref{equationDboxtimesP1} is non-split. Otherwise, since $\pi(\Delta)(*)$ is locally analytic, there exist non-zero locally analytic vectors in $\Delta\boxtimes_{\omega} \bbP^1$. But this is not possible, since one can check that the $\calD^{\la}(\begin{pmatrix}
		1 & \bbZ_p \\ 0 & 1 \end{pmatrix},E)$-module $\Delta\boxtimes \bbZ_p=\Delta\simeq \calR_E^{\oplus 2}$ has no non-zero locally analytic vectors for the action of $\begin{pmatrix}
		1 & \bbZ_p \\ 0 & 1 \end{pmatrix}$ (using for example \cite[Thm.~4.1.5]{JC23} and that $\calR_E$ is the ring of functions on the boundary of the open unit disk).

	By Colmez's construction $\Delta\boxtimes_{\omega} \bbP^1\simeq \Delta\boxtimes \bbZ_p\oplus \Delta\boxtimes p\bbZ_p$ is an LF space since $\calR_E$ is. Using \cite[Lem.~A.33]{bosco2021drinfeldproetale} and \cite[Lem.~3.2]{steingart2024higher}, it is not hard to see that the condensification of \eqref{equationDboxtimesP1} remains exact. The extension is also non-split by the same argument above, thus gives a non-zero element in 
	\begin{align*}
		\Ext^1_{\calD^{\la}(\GL_2(\bbQ_p),E)}(\pi(\Delta),\pi(\check\Delta)^{\vee})&=\Hom_{\calD^{\la}(\GL_2(\bbQ_p),E)}(\pi(\Delta),\pi(\check\Delta)^{\vee}[1])\\
		&=\Hom_{\calD^{\la}(\GL_2(\bbQ_p),E)}(\pi(\Delta),\bbD_{\ST}(\pi(\check\Delta))[1])
	\end{align*}
	where the last identity is the adjunction in \cite[Cor.~3.3.7]{JC23}. Let $f:\pi(\Delta)\rightarrow \bbD_{\ST}(\pi(\check\Delta))[1]$ be the non-zero morphism given by the extension. Then $f$ is an injection since both the source and target are admissible (hence $f$ is classical) and $\pi(\Delta)(*)$ is topologically irreducible. Applying $\bbD_{\ST}(-)[1]$, the non-zero map $\bbD_{\ST}(f)[1]:\pi(\check\Delta)\rightarrow \bbD_{\ST}(\pi(\Delta))[1]$ is also injective by the same reason, which implies that $f$ is a surjection. The result follows. 
\end{proof}
\begin{remark}
	Proposition~\ref{propositionSTselfdual} applies for all indecomposable $\Delta$ unless $\Delta$ is trianguline or $\Delta$ is de Rham with regular Hodge--Tate weights up to twist \cite[Thm.0.3]{colmez2018poids} (the dimension $1$ assumption is always satisfied by \cite[Thm.~1.5]{DPS2024GKJL}). If $\Delta$ is trianguline, $\bbD_{\ST}(\pi(\Delta))$ can be computed using \cite[\S~6]{ST05}. The result for the proposition also holds for regular de Rham $\Delta$ by Proposition~\ref{proposition:STdualpic} below.
\end{remark}
\begin{remark}
	Proposition~\ref{propositionSTselfdual} says that $\pi(\Delta)$ is self-dual up to a twist under $\bbD_{\ST}(-)$. When the representations appear in the completed cohomology or its patched versions, such self-dual result can often be obtained using the Poincar\'e duality of the completed cohomology or the big patching module, see for example \cite[Thm.~8.2]{huWang2022modp} and \cite[Appendix D]{Ding2024wall}. 
\end{remark}
\begin{lemma}\label{lemma:ST-smooth-dual}
	We have $\bbD_{\ST}(\pi_{\sm})\simeq\pi_{\sm}^{\vee,\sm}[0]$ where $ \pi_{\sm}^{\vee,\sm}\simeq \pi_{\sm}$ is the smooth contragredient of $\pi_{\sm}$.
\end{lemma}
\begin{proof}
	By \cite[Prop.~5.3.2]{JC23}, $\pi_{\sm}^{\vee,R\sm}\simeq \pi_{\sm}^{\vee,R\la}$ where $(-)^{R\sm}$ is defined for $\calD^{\sm}(G,E)$-modules in \textit{loc.cit.}. Since $(-)^{\sm}$ is exact \cite[Lem.~5.1.9]{JC23}, $\pi_{\sm}^{\vee,R\la}=\pi_{\sm}^{\vee,\sm}$. Since $\pi_{\sm}$ has trivial central character, its contragredient is $\pi_{\sm}$ itself, by
	\cite[\S~33.4, Prop. (1) \& Exercise]{bushnellHenniart2006}.
\end{proof}
Next, we determine $\bbD_{\ST}(\pi_c),\bbD_{\ST}(\widetilde{\pi})$ using the results in \cite{ding2023change}. Let $D=D_{\rm rig}(\rho)$ be the $(\varphi,\Gamma)$-module over $\calR_E$ associated to the Galois representation $\rho$ in the last subsection. Let $\Delta$ be the $p$-adic differential equation associated to $D_{\rm pst}(\rho)$. In other words, $\Delta$ is the unique $(\varphi,\Gamma)$-submodule of Hodge--Tate weights $0$ of $D_{\rm rig}(\rho)[\frac{1}{t}]$ given by \cite[Prop.~2.11]{ding2023change}, the change of weights of $D_{\rm rig}(\rho)$. Let $\lambda=(0,0)$ and $\mu=(-1,0)$ considered as the weights of the diagonal Cartan of $\mathfrak{gl}_2$. By \cite{dospinescu2012actions}, the infinitesimal characters of the Harish-Chandra center of $U(\mathfrak{gl}_2)$ of $\pi(\rho),\pi(\Delta)$ are given by that of the Verma modules (with respect to the Borel subalgebra of the upper-triangular matrices) of the highest weights $\lambda$ and $\mu$ respectively. Let $T_{\mu}^{\lambda},T_{\lambda}^{\mu}$ be the translation functor recalled in \S~\ref{sec:translation-duality}. By \cite[Thm.~1.1(2)]{ding2023change} (see also \cite[Prop.~2.4.2]{su2025locally}), there exists a short exact sequence
\begin{equation}\label{equation:translation-sequence}
	0\rightarrow \pi_c\rightarrow T_{\mu}^{\lambda}\pi(\Delta)\rightarrow \widetilde{\pi}\rightarrow 0
\end{equation}
of locally analytic representations. Moreover, one has $T_{\lambda}^{\mu}\widetilde{\pi}\simeq T_{\lambda}^{\mu}\pi_c\simeq \pi(\Delta)$.
\begin{proposition}\label{proposition:STdualpic}
	We have
	\[
 	\bbD_{\ST}(\pi_c)\simeq\widetilde\pi[-1],
 	\qquad\bbD_{\ST}(\widetilde\pi)\simeq\pi_c[-1],
\]
and 
\[\bbD_{\ST}(\pi(\rho))\simeq \pi(\rho)[-1].\]
\end{proposition}
\begin{proof}
	A similar proof of \cite[Thm.~3.11]{li2026locallyanalyticvectorscompleted}, along with \eqref{eq:universal-sequence}, shows that $\pi_c,\pi(\Delta)$ as well as $\widetilde{\pi}$ have Gelfand--Kirillov dimension $1$. Since $\pi_c(*),\pi(\Delta)(*)$ are topologically irreducible by \cite[Cor.~0.4]{colmez2018poids} and the sequence \eqref{eq:universal-sequence} is the universal extension of $\pi_c$ by $\pi_{\sm}$ in Theorem~\ref{theoremDingext1}, all these representations are pure. Hence by the arguments in the proof of Proposition~\ref{propositionSTselfdual}, $\bbD_{\ST}(\pi_c),\bbD_{\ST}(\widetilde{\pi}),\bbD_{\ST}(\pi(\Delta))$ concentrates in degree $1$. Moreover, applying $\bbD_{\ST}$ to \eqref{eq:universal-sequence} and using Lemma~\ref{lemma:ST-smooth-dual} give a short exact sequence
	\begin{equation}\label{equationSTdualpic1}
		0\rightarrow\pi_{\rm sm}^{\oplus 2}\rightarrow\bbD_{\ST}(\pi_c)[1]\rightarrow \bbD_{\ST}(\widetilde{\pi})[1]\rightarrow 0.
	\end{equation}
	On the other hand, by Proposition~\ref{propositionSTselfdual}, we have $\bbD_{\ST}(\pi(\Delta)\simeq\pi(\Delta^{\vee}(\epsilon))[-1]$. Since $\det \Delta= t^{-1}\calR(\epsilon)$, $\Delta^{\vee}(\epsilon)\simeq t\Delta$ so that $\pi(\Delta^{\vee}(\epsilon))\simeq \pi(\Delta)\otimes \det$. By \cite[Cor.~3.2.14(3)]{JC23} or using Proposition~\ref{prop:translation-duality}, $\bbD_{\ST}$ is compatible with the translation functors. Hence applying $\bbD_{\ST}$ to \eqref{equation:translation-sequence} gives another short exact sequence
	\[0\rightarrow \bbD_{\ST}(\widetilde{\pi})[1]\rightarrow T_{\mu}^{\lambda}\pi(\Delta)\rightarrow \bbD_{\ST}(\pi_c)[1]\rightarrow 0.\]
	Combined with \eqref{equationSTdualpic1} considering the irreducible constituents of $T_{\mu}^{\lambda}\pi(\Delta)(*)$, one has $\bbD_{\ST}(\widetilde{\pi})[1]\simeq \pi_c$. Hence $\bbD_{\rm ST}(\pi_c)\simeq \bbD_{\ST}(\bbD_{\ST}(\widetilde{\pi})[1])\simeq \widetilde{\pi}[-1]$.

	Note that since $\det \rho\simeq \epsilon$, we have $\rho^{\vee}\otimes\epsilon\simeq \rho$. With the determination of $\bbD_{\ST}(\pi_c)$ and the methods in Proposition~\ref{propositionSTselfdual}, the last isomorphism for $\pi(\rho)$, which is a non-split extension of $\pi_c$ by $\pi_{\sm}$, is an easy exercise. We leave it to the reader.
\end{proof}
\subsection{Bernstein--Zelevinsky duality}\label{subsectionBZdual}
In this subsection, we compute the Bernstein--Zelevinsky dual \eqref{definitionBZdual} of $\pi(\rho)$ and some relevant locally analytic representations of $G=\PGL_2(\bbQ_p)$, using Theorem~\ref{maintheorem}.
\begin{remark}
	Unlike the Schneider--Teitelbaum duality, $\bbD_{\BZ}(\pi(\rho))$ depends on whether we consider $\pi(\rho)$ as a $G$-representation or a $\GL_2(\bbQ_p)$-representation (differed by a degree shift), see Remark~\ref{remarSTdualindependent} and Remark~\ref{remarkothercentralcharacter}. 
\end{remark}
\begin{lemma}\label{lem:BZ-smooth}
There exists an isomorphism $\bbD_{\BZ}(\pi_{\sm})\simeq\pi_{\sm}[0]$.
\end{lemma}
\begin{proof}
By \cite[Rem.~5.5]{strauch2025bernsteinzelevinskydualitylocallyanalytic}, $\bbD_{\BZ}(\pi_{\sm})$ coincides with the classical Bernstein--Zelevinsky dual of $\pi_{\sm}$. The result now follows from \cite[Rem.~III.3.4]{schneiderStuhler1997representation} and the self-duality of $\pi_{\sm}$ for the smooth contragredient in Lemma~\ref{lemma:ST-smooth-dual}.
\end{proof}
We keep the notation in \S~\ref{subsec:supercuspidal-representations}. for the de Rham Galois representation $\rho$ and the Weil--Deligne representation $M$. We associated a smooth representation $\tau_{\sm}$ of $D_p^{\times}$. As in Lemma~\ref{lemma:ST-smooth-dual}, since $\rho\simeq \rho^{\vee}\otimes \epsilon$, one has $\tau_{\sm}^{\vee}\simeq \tau_{\sm}$ (cf. \cite[\S~56.1, (56.1.4)]{bushnellHenniart2006}). Let \[\calF_M:=\calF_{\tau_{\sm}}\simeq \calF_{\tau_{\sm}^{\vee}}\simeq \calF_{M}^{\vee}\] 
be the $G$-equivariant vector bundle on $X:=\Omega_{\mathrm{Dr},E}=(\mathbb P^{1,\mathrm{an}}_{\bbQ_p}\setminus\mathbb P^1(\bbQ_p))_E$ attached to $\tau_{\rm sm}$ by Definition~\ref{definitionFM} and Lemma~\ref{lemmaFMvectorbundle}. Let $\Omega^1_{X}=\calO_X(-2):=\calO_{\bbP_E^{1,\an}}(-2)|_X$ be the dualizing line bundle on $X$. By Serre duality, \'etale descent, projection formula, and Definition~\ref{definitionpippic}, we have
\begin{equation}\label{equationsheafpic}
	\begin{aligned}
		\widetilde\pi&\simeq\Gamma(X,\calF_M^\vee\otimes_{\calO_X}\Omega_X^1)^{\vee}\simeq  H_c^1(X,\calF_M)\\
		\pi_c&\simeq \Gamma(X,\calF_M^\vee)^{\vee}\simeq H_c^1(X,\calF_M\otimes_{\calO_X}\Omega_X^1).
	\end{aligned}	
\end{equation}
Moreover, by the discussion around \eqref{equation:translation-sequence}, together with the results on the translation of equivariant sheaves in \cite[\S~2]{su2025locally}, we have
\begin{align}\label{equationsheafpidelta}
	\pi(\Delta)\simeq \Gamma(X,\calF_M^\vee(-1))^{\vee}\simeq H_c^1(X,\calF_M(-1)),
\end{align}
where $(-1)$ denotes the twist by the tautological line bundle on $\bbP^{1,\an}_E$.
\begin{proposition}\label{prop:BZ-pic-pitilde}
We have isomorphisms
\[
 \theta_{\widetilde{\pi}}:\bbD_{\BZ}(\widetilde\pi)\simeq\pi_c[-1],
 \qquad \theta_{\pi_c}:\bbD_{\BZ}(\pi_c)\simeq\widetilde\pi[-1],
\]
and 
\[\theta_{T_{\lambda}^{\mu}\pi(\Delta)}:\bbD_{\BZ}(T_{\lambda}^{\mu}\pi(\Delta))\simeq T_{\lambda}^{\mu}\pi(\Delta)[-1]\]
in $\Rep^{\la}_{E_\square}(G)$.
\end{proposition}
\begin{proof}
The first two isomorphisms follow from Corollary~\ref{corBZdualDrinfeld}, \eqref{equationsheafpic}, \eqref{equationsheafpidelta} and the self-duality of $\calF_M$. The last is by the compatibility between the translation functor and $\pi(\Delta)=T_{\lambda}^{\mu}\pi_c$. 
\end{proof}
Let $d:\calO_{\Sigma_{m,E}}\rightarrow \Omega_{\Sigma_{m,E}}^1$ be the differential map (see \eqref{equationDrinfeldtowerrational} for $\Sigma_{m,E}$) and also write $d:\calF_M\rightarrow \calF_M\otimes_{\calO_X}\Omega_X^1$ for the induced map by descent.
\begin{lemma}\label{lemmaDBZd}
	The following diagram commute:
	\begin{equation}\label{equationdiagramDBZd}
		\begin{tikzcd}
			\widetilde{\pi}\arrow["\sim",d,"{\theta_{\widetilde{\pi}}}"']\arrow[r,"-d"]&\pi_c\arrow["\sim",d,"{\theta_{\pi_c}}"']\\
			\bbD_{\BZ}(\pi_c)[1]\arrow[r,"{\bbD_{\BZ}(d)[1]}"]&\bbD_{\BZ}(\widetilde{\pi})[1].
		\end{tikzcd}
	\end{equation}
\end{lemma}
\begin{proof}
	Applying the Serre duality, a calculation using the residue pairing in \cite[Lem.~2.16]{rodrigues2022solid} shows that the map $\bbD_{\GS}(d)[-1]:\calF_M\rightarrow \calF_M\otimes_{\calO_X}\Omega_X^1$ coincides with $-d$. Let $U\subset X$ be an affinoid, the following diagram commutes:
	\[\begin{tikzcd}H^0(U,\calF_M) \arrow[r,"d"] \arrow[d,"\sim"] & H^0(U,\calF_M\otimes_{\calO_X}\Omega_X^1) \arrow[d,"\sim"']\\ 
		H_c^1(U,\calF_M\otimes_{\calO_X}\Omega_X^1)^{\vee} \arrow[r,"-d^{\vee}"] & H_c^1(U,\calF_M)^{\vee}.
	\end{tikzcd}
	\]
	Now the commutativity of \eqref{equationdiagramDBZd} follows from the proof of Theorem~\ref{maintheorem} using resolutions by compactly induced representations in Proposition~\ref{propositionresolutions}. 
\end{proof}
We define
\begin{equation}\label{equationdefinitiondeRhamcohomology}
	R\Gamma_{\dR,c}(X,\calF_M):=R\Gamma_c(X, [\calF_M\xrightarrow{d}\calF_M\otimes_{\calO_X}\Omega_X^1])
\end{equation}
where the compactly supported cohomology is taken in the classical sense (cf. \S~\ref{subseccompare}). Then the dual of \eqref{eq:universal-sequencedual} induces an isomorphism
\[H^1_{\dR,c}(X,\calF_M)=\ker(\widetilde{\pi}\xrightarrow{d}\pi_c)\simeq M_{\dR}^{\vee}\otimes_{E}\pi_{\sm}.\]
Applying $\bbD_{\BZ}$, we obtain an isomorphism 
\begin{equation}\label{equationmapDBZHdR}
	\begin{aligned}
		H^1_{\dR,c}(X,\calF_M)\simeq \fib(\widetilde{\pi}\rightarrow \pi_c)\xrightarrow{\text{\eqref{equationdiagramDBZd}}} \fib(\bbD_{\BZ}(\pi_c)\rightarrow \bbD_{\BZ}(\widetilde{\pi}))[1]&\simeq \bbD_{\BZ}(\fib(\widetilde{\pi}\rightarrow \pi_c))\\
		&\simeq \bbD_{\BZ}(H^1_{\dR,c}(X,\calF_M))
	\end{aligned}
\end{equation}
which is uniquely defined up to scalars. By the definition of $\bbD_{\GS}$ and tensor-Hom adjunction (and the $E_{\square}$-flatness of the LS spaces), the above isomorphism induces a pairing 
\begin{equation}\label{equationpairingdeRham}
	H^1_{\dR,c}(X,\calF_M)\otimes_{E_{\square}}H^1_{\dR,c}(X,\calF_M)\to \calC_c^{\la}(G,E)
\end{equation}
equivariant for the actions of $G\times G$ on the two sides. Since 
\[\uHom_{\calD^{\la}(G\times G,E)}(\pi_{\sm}\otimes_{E_{\square}}\pi_{\sm},\calC_c^{\la}(G,E))\simeq \uHom_{\calD^{\la}(G,E)}(\pi_{\sm},\bbD_{\BZ}(\pi_{\sm}))\simeq E,\] 
the map \eqref{equationpairingdeRham} induces a pairing
\begin{equation}\label{equationpairingMdR}
	(-,-):M_{\dR}^{\vee}\otimes_{E_{\square}}M_{\dR}^{\vee}\rightarrow \Hom_{\calD^{\la}(G\times G,E)}(\pi_{\sm}\otimes_{E_{\square}}\pi_{\sm},\calC_c^{\la}(G,E))\simeq E
\end{equation}
which is perfect.
\begin{proposition}\label{propositionalternative}
	The pairing \eqref{equationpairingMdR} is alternating: $(x,y)=-(y,x)$ for all $ x,y\in M_{\dR}^{\vee}$.
\end{proposition}
\begin{proof}
	Suppose that there is a map $\alpha:\pi_1\rightarrow \bbD_{\BZ}(\pi_2)$, which induces by adjunction a map 
	\[\pi_1\otimes_{E_{\square}}^L\pi_2\rightarrow \calC_c^{\la}(G,E).\]
	Since $\otimes_{E_{\square}}$ is symmetric, there is a natural isomorphism $\mathrm{sw}:\pi_2\otimes_{E_{\square}}^L\pi_1\simeq \pi_1\otimes_{E_{\square}}^L\pi_2$. There is also an involution $\iota:\calC_c^{\la}(G,E)\to \calC_c^{\la}(G,E)$ induced $g\mapsto g^{-1}$. By the tensor Hom adjunction we obtain a map $\alpha^t:\pi_2\rightarrow\bbD_{\BZ}(\pi_1)$ which is the transpose of $\alpha$ with respect $\bbD_{\BZ}$. Take $\pi_1=\bbD_{\BZ}(\pi_2)$ this gives the canonical map $\eta_{\pi_2}:\pi_2\to \bbD_{\BZ}(\bbD_{\BZ}(\pi_2))$. In general $\alpha^{t}=\bbD_{\BZ}(\alpha)\circ \eta_{\pi_2}$. It's not hard to show that \eqref{equationpairingMdR} is alternating if and only if $\alpha^{t}=-\alpha$ where $\alpha$ denotes the map in \eqref{equationmapDBZHdR}. 

	We need to be extremely careful about the sign. Write $\epsilon_{\pi,s}:\bbD_{\BZ}(\pi)[-s]\xrightarrow{\sim}\bbD_{\BZ}(\pi[s])$ for the canonical isomorphism. Under the convention of \cite[\href{https://stacks.math.columbia.edu/tag/0FNG}{Tag~0FNG}, (9)]{stacks-project}, $\epsilon_{\pi,s}$ sends an element of degree $m$ in the Hom complex by $(-1)^{ms}$. Hence we have
	\begin{equation}\label{equationdualshiftcomposition}
	    \epsilon_{\pi[b],a}\circ(\epsilon_{\pi,b}[-a])=(-1)^{ab}\epsilon_{\pi,a+b}.
	\end{equation}
	For $\alpha:\pi_1\to\bbD_{\BZ}(\pi_2)$ and $s\in\bbZ$, we have
	\begin{equation}\label{equationtransposeshiftexplicit}
	    (\epsilon_{\pi_2,-s}\circ(\alpha[s]))^t=(-1)^s\epsilon_{\pi_1,s}\circ(\alpha^t[-s]).
	\end{equation}
	This is because the following diagram
	\[
	\begin{tikzcd}[column sep=large]
	    \pi_2[-s]\otimes_{E_\square}^L\pi_1[s]
	    \arrow[r,"\mathrm{sw}"]\arrow[d,"\sim"']
	    & \pi_1[s]\otimes_{E_\square}^L\pi_2[-s]
	    \arrow[d,"\sim"]\\
	    \pi_2\otimes_{E_\square}^L\pi_1
	    \arrow[r,"\mathrm{sw}"]
	    & \pi_1\otimes_{E_\square}^L\pi_2.
	\end{tikzcd}
	\]
	commutes up to $(-1)^{s^2}=(-1)^s$ under the convention of \cite[\href{https://stacks.math.columbia.edu/tag/0FNG}{Tag~0FNG}]{stacks-project}.

	We now argue as in \cite[Thm.~1.6]{balmer2000triangularWittI}. Consider the following commutative diagram
	\begin{equation}\label{equationBalmerGamma}
		\begin{tikzcd}[column sep=4.2em]
			\widetilde{\pi}\arrow[r,"d"]\arrow[d,"\theta_{\widetilde{\pi}}"',"\sim"]
			& \pi_c\arrow[r,"u_1"]\arrow[d, "\theta_{\pi_c}"',"\sim"]
			& H^1_{\dR,c}(X,\calF_M)[1]\arrow[r,"u_2"]\arrow[d,"\sim","\psi"']
			& \widetilde{\pi}[1]\arrow[d,"{\theta_{\widetilde{\pi}}[1]}"',"\sim"]
			\\
			\bbD_{\BZ}(\pi_c)[1]\arrow[r,"{-\bbD_{\BZ}(d)[1]}"]
			& \bbD_{\BZ}(\widetilde{\pi})[1]\arrow[r,"{\bbD_{\BZ}(u_2)[2]}"]
			& \bbD_{\BZ}(H^1_{\dR,c}(X,\calF_M))[1]\arrow[r,"{\bbD_{\BZ}(u_1)[2]}"]
			& \bbD_{\BZ}(\pi_c)[2].
		\end{tikzcd}
	\end{equation}
	The bottom is obtained from the top triangle by applying $\bbD_{\BZ}(-)$ and rotating. The left square commutes by \eqref{equationdiagramDBZd}. The existence of $\psi$ follows from the axiom TR 3 of the triangulated category. In our case, the map $\psi$ is uniquely defined by the above diagram using that the map $\bbD_{\BZ}(u_1)[1]$ is injective. By our definition, $\psi[-1]$ coincides with \eqref{equationmapDBZHdR}. Thus we only need to show that $\psi[-1]^{t}=-\psi[-1]$.

	Applying $\bbD_{\BZ}(-)$ to \eqref{equationBalmerGamma} and using $\eta_{\pi_c[-2]}$, $\eta_{H^1_{\dR,c}(X,\calF_M)[-1]}$, $\eta_{\widetilde{\pi}[-1]}$ and $\eta_{\pi_c[-1]}$, we obtain the following transpose of \eqref{equationBalmerGamma}:
	\begin{equation}\label{equationBalmerGammatranspose}
		\begin{tikzcd}[column sep=4.2em]
			\pi_c[-2]\arrow[r,"{u_1[-2]}"]\arrow[d, "{(\theta_{\widetilde{\pi}}[1])^t}"',"\sim"]
			& H^1_{\dR,c}(X,\calF_M)[-1]\arrow[r,"{u_2[-2]}"]\arrow[d,"\sim","{\psi^t}"']
			& \widetilde{\pi}[-1]\arrow[r,"{-d[-1]}"]\arrow[d,"{\theta_{\pi_c}^t}"',"\sim"]
			& \pi_c[-1]\arrow[d,"{\theta_{\widetilde{\pi}}^t}"',"\sim"]
			\\
			\bbD_{\BZ}(\widetilde{\pi})[-1]\arrow[r,"{\bbD_{\BZ}(u_2)}"]
			& \bbD_{\BZ}(H^1_{\dR,c}(X,\calF_M))[-1]\arrow[r,"{\bbD_{\BZ}(u_1)}"]
			& \bbD_{\BZ}(\pi_c)\arrow[r,"{\bbD_{\BZ}(d)}"]
			& \bbD_{\BZ}(\widetilde{\pi}).
		\end{tikzcd}
	\end{equation}
	Here we omitted the isomorphisms $\epsilon_{\pi,s}$. We have an identity $(\theta_{\pi_c}[-1])^t=\epsilon_{\pi_c,-1}\circ\theta_{\widetilde{\pi}}:\widetilde{\pi}\to\bbD_{\BZ}(\pi_c[-1])\simeq \bbD_{\BZ}(\pi_c)[1]$. To see this, we need know the isomorphism $\theta_{\pi_c}[-1]:\pi_c[-1]\simeq \bbD_{\BZ}(\widetilde{\pi})$ and the pairing $\widetilde{\pi}\otimes_{E_{\square}}^L\pi_c[-1]\to \calC_{c}^{\la}(G,E)$ from Proposition~\ref{prop:BZ-pic-pitilde} give the natural pairing $\widetilde{\pi}\otimes_{E_{\square}}^L\bbD_{\BZ}(\widetilde{\pi})\to \calC_{c}^{\la}(G,E)$. Taking a resolution of $\widetilde{\pi}$ in (1) of Proposition~\ref{propositionresolutions}, this follows from the proof of Theorem~\ref{maintheorem}.

	Applying \eqref{equationtransposeshiftexplicit} and using \eqref{equationdualshiftcomposition}, we get
	\[
	\begin{aligned}
	\theta_{\pi_c}^t&=-\epsilon_{\pi_c[-1],1}\circ((\theta_{\pi_c}[-1])^t[-1])\\
	&=-\epsilon_{\pi_c[-1],1}\circ (\epsilon_{\pi_c,-1}[-1])\circ(\theta_{\widetilde{\pi}}[-1])\\
	&=\theta_{\widetilde{\pi}}[-1].
	\end{aligned}
	\]
	Similarly, we have $\theta_{\widetilde{\pi}}^t=\theta_{\pi_c}[-1]$. By the commutativity of the second square of \eqref{equationBalmerGammatranspose} and the right square of \eqref{equationBalmerGamma}, we have
	\begin{equation}\label{equationBalmerGammacomparison}
		\begin{aligned}
			\bbD_{\BZ}(u_1)[1]\circ\psi^t[1]&=\theta_{\widetilde{\pi}}\circ u_2[-1],\\
			\bbD_{\BZ}(u_1)[1]\circ\psi[-1]&=\theta_{\widetilde{\pi}}\circ u_2[-1].
		\end{aligned}
	\end{equation}
	Since $\bbD_{\BZ}(u_1)[1]$ is injective, we conclude that $\psi^t[1]=\psi[-1]$. Applying \eqref{equationtransposeshiftexplicit}, we get $(\psi[-1])^t=-\psi^t[1]=-\psi[-1]$, as desired.
\end{proof}
\begin{theorem}\label{theoremDBZpirho}
	We have 
	\[\bbD_{\BZ}(\pi(\rho))\simeq \pi(\rho)[-1].\]
\end{theorem}
\begin{proof}
	Let $\calL\subset M_{\rm dR}$ and let $\pi_{\calL}$ be the $G$-representation constructed in Definition~\ref{definitionpiL}. By Lemma~\ref{lem:BZ-smooth}, Proposition~\ref{prop:BZ-pic-pitilde}, applying $\bbD_{\BZ}$ to \eqref{equationpushoutdospinescuLeBras}, we get a diagram:
\begin{equation}\label{equationdiagramtheoremDBZpirho}
\begin{tikzcd}0 \arrow[r] & \calL\otimes_E\pi_{\sm} \arrow[r,"i"] \arrow[d, hook] & \widetilde{\pi} \arrow[r] \arrow[d, equal] & \bbD_{\BZ}(\pi_{\calL})[1] \arrow[r] \arrow[d, two heads] & 0 \\
	0 \arrow[r] & M_{\rm dR}\otimes_E\pi_{\sm}\simeq\bbD_{\BZ}(M_{\rm dR}^{\vee}\otimes_E\pi_{\sm}) \arrow[r] & \widetilde{\pi}=\bbD_{\BZ}(\pi_c)[1] \arrow[r] & \pi_c\simeq \bbD_{\BZ}(\widetilde{\pi})[1] \arrow[r] & 0.
\end{tikzcd}
\end{equation}
It follows that $\bbD_{\BZ}(\pi_{\calL})[1]\simeq \pi_{\calL'}$ for some $E$-line $\calL'\subset M_{\rm dR}$ by Theorem~\ref{theoremDingext1}. The remaining task is to show $\calL'=\calL$. Notice that the second row of the above diagram is obtained by applying $\bbD_{\BZ}$ to the universal extension and corresponds to the map $M_{\dR}\simeq \Ext^1_{\calD^{\la}(G,E)}(\pi_c,\pi_{\sm})\to \Hom_{\calD^{\la}(G,E)}(\pi_{\sm},\widetilde{\pi})\simeq M_{\rm dR}^{\vee}$, or a pairing 
\begin{equation}\label{equationpairing1}
	(-,-):M_{\dR}^{\vee}\times M_{\dR}^{\vee}\to E. 
\end{equation}
The desired identity $\calL'=\calL$ is equivalent to that the composition $\calL\hookrightarrow M_{\dR}\rightarrow M_{\dR}^{\vee}\rightarrow \calL^{\vee}$ is zero. This is equivalent to that $(x,x)=0$ for all $x\in \calL$. Let $\calL$ vary, we need prove that the pairing \eqref{equationpairing1} is alternating. Since \eqref{equationmapDBZHdR} induced by $\bbD_{\BZ}$ also sends $H_{\dR,c}^1(X,\calF_M)=M_{\dR}^{\vee}\otimes_E\pi_{\sm}\twoheadrightarrow \calL^{\vee}\otimes_{E}\pi_{\sm}$ to the map $\calL\otimes_{E}\pi_{\sm}\hookrightarrow \bbD_{\BZ}(H_{\dR,c}^1(X,\calF_M))[1]=\fib(\bbD_{\BZ}(\pi_c)\rightarrow \bbD_{\BZ}(\widetilde{\pi}))[1]$ in \eqref{equationdiagramtheoremDBZpirho}, the pairing \eqref{equationpairing1} coincides with \eqref{equationpairingMdR} and we conclude by Proposition~\ref{propositionalternative}.
\end{proof}
\begin{remark}
	There is another proof for Theorem~\ref{theoremDBZpirho} using results in \S~\ref{subsectionSTdual}. For $V\in \Rep^{\la}_{E_\square}(G)$, the inclusion $\calC_c^{\la}(G,E)\hookrightarrow \calC^{\la}(G,E)=(\Ind_{\{1\}}^GE)^{\la}$ induces a map 
	\[\eta_{V}:\bbD_{\BZ}(V)^{R\la}\rightarrow \bbD_{\ST}(V)\]
	via 
	\[ R\uHom_{\calD^{\la}(G,E)}(V,\calC^{\la}(G,E))\simeq R\uHom_{E}(V,E).\]
	If $V$ is smooth and supercuspidal, one can check that $\eta_V$ is an isomorphism. Using Proposition~\ref{proposition:STdualpic}, we get a diagram
	\[\begin{tikzcd}{\pi_{\sm}} \arrow[r]\arrow[d,"{\eta_{\pi_{\sm}}}"'] & \widetilde\pi \arrow[r]\arrow[d,"{\eta_{\pi_c}[1]}"] & {\bbD_{\BZ}(\pi(\rho))[1]} \arrow[r]\arrow[d,"{\eta_{\pi(\rho)}[1]}"]
		& {\pi_{\sm}[1]}\arrow[d,"{\eta_{\pi_{\sm}}[1]}"]\\ 
		{\pi_{\sm}} \arrow[r]& \widetilde\pi
		\arrow[r] & \pi(\rho) \arrow[r]& {\pi_{\sm}[1].}
	\end{tikzcd}
	\]
	Then the first vertical arrow is an isomorphism. This forces the second vertical arrow non-zero. But $\End_{\calD^{\la}(G,E)}(\widetilde{\pi})=E$ by \cite[Prop.~2.5.2 \& Prop.~2.6.1 \& Rem.~2.5.9]{su2025locally}, we see $\eta_{\pi_c}$ is an isomorphism. We conclude that $\eta_{\pi(\rho)}$ is also an isomorphism.
\end{remark}

\subsection{Serre functor and trace pairings}
In this subsection, we combine the Schneider--Teitelbaum and Bernstein--Zelevinsky duality to establish some Serre-type duality on extension groups of locally analytic representations in \S~\ref{subsec:supercuspidal-representations}. 

Let $\langle \pi_{\sm},\pi_c\rangle\subset \Rep^{\la}_{E_\square}(G)$ be the full triangulated subcategory generated by $\pi_{\sm},\widetilde{\pi},\pi_{c}$ in the sense of \cite[\href{https://stacks.math.columbia.edu/tag/0ATG}{Tag 0ATG}]{stacks-project}. Objects in $\langle \pi_{\sm},\pi_c\rangle$ can be obtained from $\pi_{\sm},\pi_c$ via finitely many operations consisting of cones/extensions, shifts and retracts. In particular, $\pi(\rho)\in \langle \pi_{\sm},\pi_c\rangle$.

We first collect the following finiteness result.
\begin{proposition}\label{prop:four-perfectness}
For any $\pi_1,\pi_2\in \langle \pi_{\sm},\pi_c\rangle$, the object $R\uHom_{\calD^{\la}(G,E)}(\pi_1,\pi_2)$ is perfect in $D(E_\square)$.
\end{proposition}

\begin{proof}
By \cite[\S15.5]{bushnellHenniart2006}, $\pi_{\sm}\cong\cInd_J^{G}\sigma$ for a compact open subgroup $J\subset G$ and a finite-dimensional representation $\sigma$ of $J$. This implies that $\pi_{\sm}$ is of compact presentation; see also Example~\ref{rem:smooth-finite-presentation}. By Proposition~\ref{propositionDrinfeldcompactpresentation} and \eqref{equationsheafpic}, $\pi_c$ is of compact presentation. Moreover, $\pi_{\sm},\pi_c$ are admissible. Hence $R\uHom_{\calD^{\la}(G,E)}(\pi_1,\pi_2)$ is perfect if $\pi_1,\pi_2\in\{\pi_c,\pi_{\sm}\}$ by Theorem~\ref{thm:finite-presentation-perfect}. The result now follows by d\'evissage.
\end{proof}

\begin{proposition}\label{prop:BZ-tensor-Hom}
Let $\pi\in\langle \pi_{\sm},\pi_c\rangle$ or $\pi$ is of compact presentation, and let $V\in\Rep^{\la}_{E_\square}(G)$ such that the underlying solid $E$-vector space of $V$ is nuclear \cite[Def.~A.10]{bosco2023rationalpadicHodge}. There is a natural equivalence
\[
 \bbD_{\BZ}(\pi)\otimes^{L\square}_{\calD^{\la}(G,E)}V[-3]\simeq R\uHom_{\calD^{\la}(G,E)}(\pi,V).
\]
\end{proposition}
\begin{proof}
We first construct a natural map between the two sides of the desired isomorphism. Let $K\subset G$ be an open compact subgroup. By \cite[Prop.~3.2.12(2)\& Prop.~A.2.1(4)]{JC23}, one has (see Lemma~\ref{lem:Rlaastensorkernel} and Proposition~\ref{prop:DlaGRla} for the $G$-equivariance):
\begin{align*}
	\calC_c^{\la}(G,E)\otimes_{\calD^{\la}(G,E)}^{L\square}V[-3]&\simeq \calC_c^{\la}(K,E)\otimes_{\calD^{\la}(K,E)}^{L\square}V[-3]\\
	&\simeq R\uHom_{\calD^{\la}(K,E)}(\calD^{\la}(K,E),V)\\
	&\simeq V.
\end{align*} 
For any $\pi\in \Rep^{\la}_{E_\square}(G)$, the natural map
\[ \pi\otimes_E^{L\square} \bbD_{\BZ}(\pi)\otimes^{L\square}_{\calD^{\la}(G,E)}V[-3]\to \calC_c^{\la}(G,E)\otimes^{L\square}_{\calD^{\la}(G,E)}V[-3]\simeq V\]
is $\calD^{\la}(G,E)$-equivariant, by tensor-Hom adjunction, inducing a natural map
\[f_\pi:\bbD_{\BZ}(\pi)\otimes^{L\square}_{\calD^{\la}(G,E)}V[-3]\rightarrow R\uHom_{\calD^{\la}(G,E)}(\pi,V).\] 

Since $\pi_c,\pi_{\sm}$ are of compact presentation, by d\'evissage, it suffices to prove that $f_\pi$ is an equivalence when $\pi=\cInd_K^{G}W$ where $K\subset G$ is compact open and $W$ is a locally analytic $K$-representation whose underlying solid $E$-space is Smith. By \cite[Prop.~A.2.1(4)]{JC23},
\begin{align*}
 R\uHom_{\calD^{\la}(G,E)}(\pi,V) \simeq R\uHom_{\calD^{\la}(K,E)}(W,V) \simeq R\uHom_{\calD^{\la}(K,E)}
 (E,R\uHom_{E}(W,V)).
\end{align*}
Since $V$ is nuclear and $W$ is Smith, the evaluation map $W^\vee\otimes^L_{E_\square}V\to R\uHom_{E}(W,V)$ is an equivalence by \cite[Prop.~5.35]{andreychev2021pseudocoherent} and \cite[Lem.~3.8(2) \& Lem.~3.10]{rodrigues2022solid}. On the other hand, Proposition \ref{BZdual} gives
\[
 \bbD_{\BZ}(\pi)\simeq
 W^\vee\otimes^{L\square}_{\calD^{\la}(K,E)}\calD^{\la}(G,E).
\]
By \cite[Prop.~3.2.12(2)]{JC23} again, we get
\[
 R\uHom_{\calD^{\la}(K,E)}
 (E,W^\vee\otimes^{L}_{E_\square}V)
 \simeq(W^\vee\otimes^{L\square}_{\calD^{\la}(K,E)}V)[-3].
\]
This gives the desired equivalence for $f_{\pi}$. 
\end{proof}
\begin{remark}\label{BZ-tensor-Hom}
	The same proof as above shows that for a general $p$-adic Lie group $G$, $ \pi,V\in \mathrm{Rep}_{E_{\square}}^{\la}(G)$ with $\pi \in \mathrm{Rep}_{E_{\square}}^{\la}(G)^{\omega}$ and with the underlying $E$-space of $V$ nuclear, we have a natural isomorphism
	\begin{align*}
		(\bbD_{\BZ}(\pi)\otimes \frd_G^*)\otimes^{L\square}_{\calD^{\la}(G,E)}V[-\dim G]\simeq R\uHom_{\calD^{\la}(G,E)}(\pi,V).
	\end{align*}
	where $\frd_G^*$ is the right character of $G$ in \S~\ref{sec:STdual}. Moreover, by \cite[Lem.~4.4.6]{heyer20246functor}, the above isomorphism holds for any prim $\pi$ and arbitrary $V$, cf. Theorem~\ref{theoremprimobject}.
\end{remark}
\begin{definition}\label{def:serre-functor}
	We define the ``Serre functor" $\calS:\Rep^{\la}_{E_\square}(G)\to \Mod_{\square}(\calD^{\la}(G,E))$ by
	\[
	V\mapsto \calS(V):=\bbD_{\ST}(\bbD_{\BZ}(V))[\dim G].
	\]
\end{definition}
\begin{corollary}\label{cor:ST-BZ-comparison}
For $\pi_1,\pi_2\in\langle \pi_{\sm},\pi_c\rangle$, there is a natural equivalence of perfect complexes in $D(E_\square)$:
\[
 R\uHom_{E}(
 R\uHom_{\calD^{\la}(G,E)}(\pi_1,\pi_2),E)
 \simeq R\uHom_{\calD^{\la}(G,E)}
 (\pi_2,\calS(\pi_1)).
\]
\end{corollary}
\begin{proof}
	The representations $\pi_{\sm}$ and $\pi_c$ are admissible, in particular they are nuclear by \cite[Prop.~A.11]{bosco2023rationalpadicHodge} being LS of compact type. Hence any object of $\langle \pi_{\sm},\pi_c\rangle$ is nuclear by \cite[Lem.~5.40]{andreychev2021pseudocoherent} and is admissible in the sense of \cite[Def.~4.2.1]{JC23}. Then one has
\[
\begin{aligned}
 R\uHom_{E}(
 R\uHom_{\calD^{\la}(G,E)}(\pi_1,\pi_2),E)&
 \simeq R\uHom_{E}(\bbD_{\BZ}(\pi_1)\otimes_{\calD^{\la}(G,E)}^{L\square}\pi_2,E)[3]\\
 &\simeq R\uHom_{\calD^{\la}(G,E)}(\bbD_{\BZ}(\pi_1),R\uHom_E(\pi_2,E))[3]\\
 &\simeq R\uHom_{\calD^{\la}(G,E)}
 (\bbD_{\BZ}(\pi_1),\bbD_{\ST}(\pi_2))[3]\\
 &\simeq R\uHom_{\calD^{\la}(G,E)}
 (\pi_2,\bbD_{\ST}(\bbD_{\BZ}(\pi_1))[3])\\
&=R\uHom_{\calD^{\la}(G,E)}(\pi_2,\calS(\pi_1)).
\end{aligned}
\]
where the first isomorphism is Proposition~\ref{prop:BZ-tensor-Hom}, the third one is the adjunction in \cite[Cor.~3.3.7]{JC23}, and the fourth is by Lemma~\ref{lem:ST-transpose}.
\end{proof}

We now define trace pairings on extension groups using Corollary~\ref{cor:ST-BZ-comparison}. For each $V\in\{\pi_{\sm},\widetilde\pi,\pi_c,\pi(\rho)\}$, results in \S~\ref{subsectionSTdual} and \S~\ref{subsectionBZdual} give an isomorphism $\bbD_{\ST}(\bbD_{\BZ}(V))\simeq V$ and we fix one for each $V$. Let
\[
 \Phi_V:R\Hom_{\calD^{\la}(G,E)}(V,V[3])\xrightarrow{\sim}R\Hom_E(R\uHom_{\calD^{\la}(G,E)}(V,V),E)
\]
be the inverse of the equivalence in Corollary~\ref{cor:ST-BZ-comparison}, after taking the underlying complex of abstract $E$-vector spaces. Define the trace map associated to $V$ by
\[
 \trace_V:\Ext^3_{\calD^{\la}(G,E)}(V,V)\longrightarrow E,\qquad x\longmapsto H^0(\Phi_V)(x)(\mathrm{id}_V).
\]
For $V,W\in\{\pi_{\sm},\widetilde\pi,\pi_c,\pi(\rho)\}$ and $i\in\mathbb Z$, the \emph{trace pairing} is
\[
 \Ext^i_{\calD^{\la}(G,E)}(V,W)\times\Ext^{3-i}_{\calD^{\la}(G,E)}(W,V)\longrightarrow E,\qquad(f,g)\longmapsto\trace_V(g\circ f),
\]
where $g\circ f$ denotes the Yoneda product. 
\begin{proposition}\label{prop:trace-pairing-perfect}
The trace pairing above is perfect for every $V,W\in\{\pi_{\sm},\widetilde\pi,\pi_c,\pi(\rho)\}$ and $i\in\mathbb Z$.
\end{proposition}
\begin{proof}
Fix $V,W\in\{\pi_{\sm},\widetilde\pi,\pi_c,\pi(\rho)\}$. By Corollary~\ref{cor:ST-BZ-comparison}, we obtain a natural equivalence
\[
 \Phi_{V,W}:R\uHom_{\calD^{\la}(G,E)}(W,V[3])\xrightarrow{\sim}R\uHom_E(R\uHom_{\calD^{\la}(G,E)}(V,W),E).
\]
Taking cohomology in degree $-i$ and the underlying abstract $E$-spaces (using Proposition~\ref{prop:four-perfectness}) give an isomorphism
\[
 \Ext^{3-i}_{\calD^{\la}(G,E)}(W,V)\xrightarrow{\sim}\Hom_E(\Ext^i_{\calD^{\la}(G,E)}(V,W),E).
\]

We check that this isomorphism is induced by the trace pairing. For $f:V\to W[i]$ and $g:W\to V[3-i]$, naturality of $\Phi_{V,W}$ gives a commutative diagram:
\[
\begin{tikzcd}[column sep=5.2em,row sep=large]
	\Hom(W[i],V[3]) \arrow[r,"{H^{0}(\Phi_{V,W[i]})}"] \arrow[d,"{-\circ f}"'] & \Hom_{E}(\Hom(V,W[i]),E) \arrow[d,"{(f\circ-)^{\vee}}"] \\
	\Hom(V,V[3]) \arrow[r,"{H^{0}(\Phi_{V})}"] & \Hom_{E}(\Hom(V,V),E).
\end{tikzcd}
\]
Hence
\[
 H^0(\Phi_{V,W[i]})(g[i])(f)=H^0(\Phi_V)(g\circ f)(\mathrm{id}_V)=\trace_V(g\circ f).
\]
Thus the isomorphism $H^0(\Phi_{V,W[i]})$ sends $g[i]$ to $f\mapsto \trace_V(g\circ f)$. Since $H^0(\Phi_{V,W[i]})$ is an isomorphism between finite-dimensional vector spaces, the trace pairing is perfect.
\end{proof}

\clearpage \appendix
\section{Schneider--Teitelbaum duality for locally analytic representations}\label{sec:STdual}
In this appendix we collect some results on Schneider--Teitelbaum duality \cite{ST05} for solid locally analytic representations, using the machinery developed in \cite{rodrigues2022solid,JC23}.

Let $E$ be a finite extension of $\bbQ_p$.  Let $G$ be a $p$-adic Lie group of dimension $d$.  Similarly as $\calC^{\la}(G,E),\calC_c^{\la}(G,E),\calD^{\la}(G,E)$, we write $\calC^{\sm}(G,E)$, $\calC_c^{\sm}(G,E)$, and $\calD^{\sm}(G,E)$ for the space of smooth functions, the space of compactly supported smooth functions, and the smooth distribution algebra in \cite[\S~5]{JC23}, respectively. Let
\[
 \overline{\calD}^{\la}(G,E)
 :=\underline{\Hom}_{E_\square}(\calC_c^{\la}(G,E),E)
\]
be the space of locally analytic distributions on $G$, with its natural $\calD^{\la}(G,E)$-bimodule structure. We work in $D_{\square}(\calD^{\la}(G,E))$, the stable $\infty$-category of solid $\calD^{\la}(G,E)$-modules. 

Let $\mathfrak g=E\otimes_{\bbQ_p}\Lie G$ and $\Delta_G=\bigwedge_E^d\mathfrak g$, with the adjoint $G$-action.  Let $\delta_G$ be the smooth modulus character, normalized by $\mu_G(gXg^{-1})=\delta_G(g)\mu_G(X)$ for compact open subsets $X\subset G$, where $\mu_G$ is a fixed $E$-valued left Haar measure on $G$, and put $\mathfrak d_G=\Delta_G\otimes_E\delta_G$. For $\chi\in\{\Delta_G,\delta_G,\mathfrak d_G\}$, we view $\chi$ as a left representation and write $\chi^*$ for its $E$-linear dual, regarded as a right representation by $(\lambda\cdot g)(v)=\lambda(gv)$. Through the involution $g\mapsto g^{-1}$, $\chi^*$ becomes the usual left representation $\chi^{-1}$, which is the inverse of $\chi$. 
\begin{remark}
	The right character $\mathfrak d_G^*$ coincides with the right character with the same notation in \cite[Prop.~3.5]{ST05}. Note that in \cite[p.~301]{ST05}, $\delta_G$ denotes the same modulus character but $G$ acts non-trivially on the right.
\end{remark}
In the next statements, $(-)^{R\la}$ on $\calD^{\la}(G,E)$ is taken for the left regular action.
\begin{lemma}\label{lem:Rlaastensorkernel}
For $M\in D_{\square}(\calD^{\la}(G,E))$, the natural morphism
\[
 \calD^{\la}(G,E)^{R\la}
 \otimes^L_{\calD^{\la}(G,E)}M
 \longrightarrow M^{R\la}
\]
is an equivalence.
\end{lemma}

\begin{proof}
This is essentially a reformulation of \cite[Cor.~3.2.14(1)]{JC23}. The natural map $\calD^{\la}(G,E)^{R\la}\allowbreak\to\calD^{\la}(G,E)$ induces a morphism $\calD^{\la}(G,E)^{R\la}\otimes^L_{\calD^{\la}(G,E)}M\to \calD^{\la}(G,E)\otimes^L_{\calD^{\la}(G,E)}M\simeq M$. Its source is derived locally analytic, so by \cite[Cor.~3.3.7]{JC23} this morphism factors through $M^{R\la}$. By \cite[Prop.~3.2.6(3)]{JC23}, both sides of the morphism preserve colimits. It is therefore enough to take $M=\calD^{\la}(G,E)\otimes_{E_\square}E_\square[S], S\in\mathrm{ExtDisc}$ to check the equivalence, which follows from \cite[Cor.~3.2.15]{JC23}.
\end{proof}
\begin{lemma}\label{lemmadualcharacter}
	There exists a natural isomorphism 
	\[R\uHom_{\calD^{\la}(G,E)}(E,\overline{\calD}^{\la}(G,E))\simeq \mathfrak{d}_G^*[-d]\]
	of right $\calD^{\la}(G,E)$-modules.
\end{lemma}
\begin{proof}
	If $G$ is compact, this is \cite[Lem.~3.2.10]{JC23}. For non-compact $G$, this is a consequence of \cite[Cor.~3.7]{ST05} whose proof can be adapted to the solid setting. 
\end{proof}
We prove the following proposition, which is \cite[Cor.~3.2.14(1)]{JC23} when $G$ is compact, in order to identify the Schneider--Teitelbaum duality for locally analytic representations with the duality for solid $\calD^{\la}(G,E)$-modules. See Proposition~\ref{prop:ST-linear-dual} and \cite[Rem.~1.2.3(1)]{JC23}.
\begin{proposition}\label{prop:DlaGRla}
There is a natural equivalence of $\calD^{\la}(G,E)$-bimodules
\[
 \calD^{\la}(G,E)^{R\la}
 \simeq
 \calC_c^{\la}(G,E)\otimes_{E_\square}\mathfrak d_G^*[-d].
\]
Put $l_{g_1}f(g)=f(g_1^{-1}g)$ and $r_{g_2}f(g)=f(gg_2^{-1})$. Then the regular actions of $g$ on $\calD^{\la}(G,E)$ correspond to the actions
\[
 g\cdot(f\otimes\lambda)=l_gf\otimes\lambda,
 \qquad
 (f\otimes\lambda)\cdot g=r_gf\otimes(\lambda\cdot g),
\]
on the right, where $g\in G$, $f\in \calC_c^{\la}(G,E)$, and $\lambda\in\mathfrak d_G^*$.
\end{proposition}

\begin{proof}
Let $K\subset G$ be an open compact subgroup. By \cite[Cor.~3.2.15]{JC23}, there is an isomorphism of $\calD^{\la}(K,E)$-bimodules
	\[
	 \calD^{\la}(K,E)^{R\la} \simeq \calC^{\la}(K,E)\otimes_{E_\square}\Delta_G^*[-d].
	\]
	Under the above isomorphism, the actions on $\calC^{\la}(K,E)$ are given by $l_g$ and $r_g$, by its proof using \cite[Cor.~3.2.14 \& Prop.~A.2.1(3)]{JC23}. Applying $\cInd_K^G$ we obtain a left $\calD^{\la}(G,E)$-linear equivalence
	\begin{equation}\label{equationDlaGRla1}
		\calD^{\la}(G,E)^{R\la}
	 \simeq \calC_c^{\la}(G,E)\otimes_{E_\square}\Delta_G^*[-d],	
	\end{equation}
	which is also right $\calD^{\la}(K,E)$-linear. Then it suffices to identify the right $\calD^{\la}(G,E)$-action on the right-hand side. 

	Consider the space $\calC^{\la}(G,\calD^{\la}(G,E))$ with $G$-actions. If $f: G\to \calD^{\la}(G,E)$ and $g_i\in G_i=G, i=1,2,3,4$, then $(g_1,g_2,g_3,g_4).f(-)=g_3f(g_1^{-1}-g_2)g_4$.  By \cite[Def.~3.2.3(3)]{JC23}, $\calD^{\la}(G,E)^{R\la}\simeq R\uHom_{\calD^{\la}(G,E)}(E, \calC^{\la}(G,\calD^{\la}(G,E))_{*_{1,3}})$ where $*_{1,3}$ denotes the diagonal action of $G\subset G_1\times G_3$. 

	If $K\subset G$ is compact open, one has $\calC^{\la}(K,\calD^{\la}(G,E))\simeq \calC^{\la}(K,E)\otimes_{E_{\square}}\calD^{\la}(G,E)$. For $g\in G$, let $\delta_g$ be the Dirac measure. There is a twisting isomorphism 
	\[\calC^{\la}(K,E)_0\otimes_{E_{\square}}\delta_g\calD^{\la}(K,E)\simeq \calC^{\la}(gK,E)\otimes_{E_{\square}}\delta_g\calD^{\la}(K,E), f\otimes \delta_g\mu\mapsto \sum_i \delta_g\mu_{1,i}.f\otimes \delta_g\mu_{2,i}\] 
	if $\nabla(\mu)=\sum_i\mu_{1,i}\otimes\mu_{2,i}$ with the inverse $f\otimes \delta_g\mu\mapsto\sum_i S(\mu_{1,i}) \delta_g^{-1}f\otimes \delta_g\mu_{2,i}$ ($\nabla$ is the comultiplication, and $S$ is the antipode of $\calD^{\la}(K,E)$, cf. \cite[Prop.~A.2.1(3)]{JC23}). Here $\calC^{\la}(K,E)_0$ denotes the same space of $\calC^{\la}(K,E)$ with the trivial left actions of $G$. We obtain a map 
\[\calC^{\la}(K,E)_0\otimes_{E_{\square}}\delta_g\calD^{\la}(K,E)\rightarrow \calC^{\la}(gK,\calD^{\la}(G,E))\] 
for each $g\in G/K$. Since $\calC^{\la}(G,\calD^{\la}(G,E))=\prod_{g\in G/K}\calC^{\la}(gK,\calD^{\la}(G,E))$ and $\overline{\calD}^{\la}(G,E)=\prod_{g\in G/K}\delta_g\calD^{\la}(K,E)$, taking direct product, we obtain a morphism 
\begin{align*}
	\calC^{\la}(K,E)_0\otimes_{E_{\square}}\overline{\calD}^{\la}(G,E)\rightarrow \calC^{\la}(G,\calD^{\la}(G,E)).
\end{align*} 
The map is $K$-equivariant for the right translations on $\calC^{\la}(K,E)_0$ and the action of $K\subset G_2$ on the target. By the Frobenius reciprocity, we obtain the following map 
\begin{align}\label{equationDlaGRla2}
	\Phi:\calC_c^{\la}(G,E)_0\otimes_{E_{\square}}\overline{\calD}^{\la}(G,E)\rightarrow \calC^{\la}(G,\calD^{\la}(G,E)),
\end{align}
which extends the more transparent map where $\overline{\calD}^{\la}(G,E)$ is replaced by $\calD^{\la}(G,E)$. Informally, for $h\in\calC^{\la}(G,E)$ and $x\in G$, one has $\Phi(f\otimes \mu)(x)(h)=\mu(y\mapsto f(y^{-1}x)h(y))$. The map $\Phi$ is equivariant for the following actions: $g_2,g_4\in G_2\times G_4$ acts on the left hand side by $f\otimes \mu\mapsto f(g_4-g_2)\otimes \mu\delta_{g_4}$, and the $*_{1,3}$-action on the right corresponds to the left action on the $\overline{\calD}^{\la}(G,E)$-factor. Taking derived $G$-invariants of \eqref{equationDlaGRla2} for the last actions and using Lemma~\ref{lemmadualcharacter}, we get a $G_2\times G_4$-equivariant map
\begin{align}\label{equationDlaGRla3}
	\calC_c^{\la}(G,E)\otimes_{E_{\square}}\mathfrak{d}_G^*[-d]\to \calD^{\la}(G,E)^{R\la}.
\end{align} 
Using the isomorphism $f\to (g\mapsto f(g^{-1}))$ of $\calC^{\la}_c(G,E)$ we recover the desired actions of $G$ on the left hand side. Comparing with \eqref{equationDlaGRla1} and restricting \eqref{equationDlaGRla3} to $\calC_c^{\la}(K,E)\otimes_{E_{\square}}\mathfrak{d}_K^*[-d]\to \calD^{\la}(K,E)^{R\la}$ we see \eqref{equationDlaGRla3} is an isomorphism.
\end{proof}

For a compact open subgroup $K\subset G$, let $\calD^h(K,E)=\calD^{(h^+)}(K,E)$ be the $h$-analytic distribution algebras and $\mathrm{Mod}^{\mathrm{qc}}_{E_\square}(\calD^{\la}(K,E))$ the category of quasi-coherent modules over $\calD^{\la}(K,E)$ of \cite[Def.~4.1.1]{JC23}.  The objects of $\mathrm{Mod}^{\mathrm{qc}}_{E_\square}(\calD^{\la}(K,E))$ are compatible systems $(M_h)_h$ of solid $\calD^h(K,E)$-modules in $D_{\square}(\calD^h(K,E))$.  There is a localization functor:
\[
	D_{\square}(\calD^{\la}(K,E))\to \mathrm{Mod}^{\mathrm{qc}}_{E_\square}(\calD^{\la}(K,E)),\quad M\mapsto j^*M=(\calD^h(K,E)\otimes^{L\square}_{\calD^{\la}(K,E)}M)_h,
\]
which admits  a fully faithful right adjoint $(M_h)_h\mapsto j_*(M_h)_h=R\varprojlim_hM_h$ and a fully faithful left adjoint $(M_h)_h\mapsto j_!(M_h)_h=(R\varprojlim_hM_h)^{R\la}$ \cite[Lem.~4.1.4 \& Thm.~4.1.5]{JC23}.  We write $\Rep^{\mathrm{ad}}_{E_\square}(G)\subset \Rep^{\la}_{E_\square}(G)$ for the full subcategory consisting of admissible representations, whose restriction to $K$ is admissible in the sense of \cite[Def.~4.2.1]{JC23}.  Let $\bbD_{\mathrm{qc}}$ denote the duality on $\mathrm{Mod}^{\mathrm{qc}}_{E_\square}(\calD^{\la}(K,E))$ of \cite[Def.~4.1.11]{JC23}.

\begin{definition}\label{defSTdual}
For $V\in D_{\square}(\calD^{\la}(G,E))$, let $V^\vee:=R\underline{\mathrm{Hom}}_{E_\square}(V,E)$ with the usual dual action. For $V\in\Rep^{\la}_{E_\square}(G)$, put
\[
 \bbD_{\ST}(V):=(V^\vee)^{R\la}.
\]
For $M\in D_{\square}(\calD^{\la}(G,E))$, put
\[
 \bbD(M):= R\underline{\mathrm{Hom}}_{\calD^{\la}(G,E)}
 (M,\overline{\calD}^{\la}(G,E) \otimes_{E_{\square}}\mathfrak d_G[d]).
\]
Here $\overline{\calD}^{\la}(G,E) \otimes_{E_{\square}}\mathfrak d_G[d]\simeq (\calC_c^{\la}(G,E)\otimes_{E_\square}\mathfrak d_G^*[-d])^{\vee}$ is equipped with the (left and right) dual actions of $G$ in Proposition~\ref{prop:DlaGRla} and $\bbD(M)$ is a left $\calD^{\la}(G,E)$-module induced by the right action on $\overline{\calD}^{\la}(G,E) \otimes_{E_{\square}}\mathfrak d_G[d]$ and the involution  $g\mapsto g^{-1}$ of $G$.
\end{definition}

\begin{proposition}\label{prop:ST-linear-dual}
For every solid $\calD^{\la}(G,E)$-module $M$, there is a natural equivalence of solid $\calD^{\la}(G,E)$-modules
\[
 (M^{R\la})^\vee\simeq\bbD(M).
\]
In particular, $\bbD_{\ST}(V)^\vee\simeq\bbD(V^\vee)$ for $V\in\Rep^{\la}_{E_\square}(G)$.
\end{proposition}

\begin{proof}
Lemma~\ref{lem:Rlaastensorkernel}, tensor-Hom adjunction, and Proposition~\ref{prop:DlaGRla} give
\begin{align*}
 (M^{R\la})^\vee
 &\simeq R\underline{\mathrm{Hom}}_{\calD^{\la}(G,E)}
 (M,(\calD^{\la}(G,E)^{R\la})^\vee)\\
 &\simeq R\underline{\mathrm{Hom}}_{\calD^{\la}(G,E)}
 (M,\overline{\calD}^{\la}(G,E)\otimes_{E_{\square}}\mathfrak d_G[d]),
\end{align*}
which is $G$-equivariant.
\end{proof}

We now collect some adjunction and biduality properties of the Schneider--Teitelbaum duality functor $\bbD_{\ST}$.
\begin{lemma}\label{lem:ST-transpose}
For $V,W\in\Rep^{\la}_{E_\square}(G)$, there is a natural equivalence in $D(E_\square)$:
\[
 R\underline{\mathrm{Hom}}_{\calD^{\la}(G,E)}
 (V,\bbD_{\ST}(W))
 \simeq
 R\underline{\mathrm{Hom}}_{\calD^{\la}(G,E)}
 (W,\bbD_{\ST}(V)).
\]
\end{lemma}

\begin{proof}
The adjunction of $(-)^{R\la}$ in \cite[Cor.~3.3.7]{JC23} and tensor-Hom adjunction give
\begin{align*}
 R\underline{\mathrm{Hom}}_{\calD^{\la}(G,E)}
 (V,\bbD_{\ST}(W))
 &\simeq
 R\underline{\mathrm{Hom}}_{\calD^{\la}(G,E)}(V,W^\vee)\\
 &\simeq
 R\underline{\mathrm{Hom}}_{\calD^{\la}(G,E)}(W,V^\vee)\\
 &\simeq
 R\underline{\mathrm{Hom}}_{\calD^{\la}(G,E)}
 (W,\bbD_{\ST}(V)).
\end{align*}
Here, the second equivalence follows because \cite[Prop.~A.2.1(2),(4)]{JC23} identifies both sides with  $R\underline{\mathrm{Hom}}_{\calD^{\la}(G,E)}(V\otimes^L_{E_\square}W,E)\simeq R\underline{\mathrm{Hom}}_{\calD^{\la}(G,E)}(E,R\uHom_E(V\otimes^L_{E_\square}W,E))$.
\end{proof}
For $V\in \Rep_{E_{\square}}^{\la}(G)$, the canonical map $V\to V^{\vee\vee}$ together with $\bbD_{\ST}(V)\to V^{\vee}$ induces a map $V\longrightarrow V^{\vee\vee}\longrightarrow \bbD_{\ST}(V)^\vee$
which factors through a canonical map 
\[V\longrightarrow \bbD_{\ST}(\bbD_{\ST}(V))=:\bbD_{\ST}^2(V)\]
by the adjunction of $(-)^{R\la}$.
\begin{proposition}\label{prop:ST-biduality-admissible}
For $V\in\Rep^{\mathrm{ad}}_{E_\square}(G)$, the representation $\bbD_{\ST}(V)$ is admissible and the natural morphism $V\to\bbD_{\ST}^2(V)$ is an equivalence.
\end{proposition}

\begin{proof}
It suffices to prove the result after restricting to an open compact subgroup $K\subset G$. Let $j_!,j_*,j^*$ and $\bbD_{\mathrm{qc}}$ be the functors in \cite[\S~4]{JC23} recalled above. Since $V$ is admissible, $(M_h)_h:=j^*V$ is coherent in the sense of \cite[Def.~4.2.1]{JC23} and $V=j_!j^*V$. Each $M_h$ is perfect, so $\bbD_{\mathrm{qc}}((M_h)_h)=(R\uHom_{\calD^h(K,E)}(M_h,\calD^h(K,E)\otimes_{E_{\square}}\mathfrak{d}_G[d]))_h$ is coherent and the biduality for perfect objects gives $\bbD_{\mathrm{qc}}^2(j^*V)\simeq j^*V$. By \cite[Prop.~4.2.7]{JC23}, we have $\bbD_{\ST}(V)\simeq j_!\bbD_{\mathrm{qc}}j^*V$, which is admissible. Using the same proposition and the full faithfulness of $j_!$, we get
\begin{align*}
 \bbD_{\ST}^2(V)
 &\simeq j_!\bbD_{\mathrm{qc}}j^*j_!\bbD_{\mathrm{qc}}j^*V\\
 &\simeq j_!\bbD_{\mathrm{qc}}^2j^*V\\
 &\simeq j_!j^*V=V.
\end{align*}
The result follows.
\end{proof}

\begin{corollary}
The functor $\bbD_{\ST}$ restricts to an anti-equivalence
\[
 \bbD_{\ST}:(\Rep^{\mathrm{ad}}_{E_\square}(G))^{\mathrm{op}}
 \xrightarrow{\sim}
 \Rep^{\mathrm{ad}}_{E_\square}(G),
\]
whose quasi-inverse is $\bbD_{\ST}$. Let $V\in\Rep^{\la}_{E_\square}(G)$ and $W\in\Rep^{\mathrm{ad}}_{E_\square}(G)$.  Applying $\bbD_{\ST}$ induces a natural equivalence in $D(E_\square)$
\[
 R\underline{\mathrm{Hom}}_{\calD^{\la}(G,E)}(V,W)
 \simeq
 R\underline{\mathrm{Hom}}_{\calD^{\la}(G,E)}
 (\bbD_{\ST}(W),\bbD_{\ST}(V)).
\]
\end{corollary}

\begin{proof}
Proposition~\ref{prop:ST-biduality-admissible} shows that $\bbD_{\ST}$ preserves admissibility and that the natural morphism $\mathrm{id}\to\bbD_{\ST}^2$ is an equivalence. The isomorphism for $R\uHom_{\calD^{\la}(G,E)}(-,-)$ follows from Lemma~\ref{lem:ST-transpose} and $W\simeq \bbD_{\ST}^2(W)$.
\end{proof}

\section{Translation functors and duality functors}\label{sec:translation-duality}
In this appendix, we give a solid treatment of translation functors for locally analytic representations, extending the usual formalism in \cite{jena2021translation}. We discuss their compatibility with Schneider--Teitelbaum and Bernstein--Zelevinsky duality. We will also give similar statements for representations with fixed central character.

We start with the general framework for translation functors. Let $L$ be a finite extension of $\bbQ_p$. Let $\mathbb G$ be a connected split reductive group over $L$, put $G=\mathbb G(L)$ and see it as a locally $\bbQ_p$-analytic group, and fix a split maximal torus $T=\bbT(L)$ contained in a fixed Borel subgroup $B$ of $\mathbb G$. Take  $E$ large enough such that $\Hom_{\bbQ_p}(L,E)=[L:\bbQ_p]$ and write $\mathfrak g=E\otimes_{\bbQ_p}\Lie G$ and $\frz=Z(U(\mathfrak g))$, the Harish--Chandra center of the universal enveloping algebra. The image of $\frz$ is central in $\calD^{\la}(G,E)$ \cite[Prop.~3.7]{ST02}.  For $V\in D_{\square}(\calD^{\la}(G,E))$ and $\mathfrak{m}\in\mathrm{MaxSpec}(\frz)$ (the spectrum of maximal ideals of $\frz$), put
\[
 V_{\mathfrak{m}}:=\colim_{n\geq1}
 R\uHom_{\frz}(\frz/\mathfrak{m}^n,V).
\]
Then we get a natural map $V_{\mathfrak{m}}\to V$. Write $\pr_{\mathfrak{m}}(V):=V_{\mathfrak{m}}$.

\begin{definition}\label{def:infinitesimal-character-finiteness}
We say that $V$ is $\frz$-finite if there is a finite subset $\Sigma\subset\mathrm{MaxSpec}(\frz)$ for which the natural morphism
\[
 \bigoplus_{\mathfrak{m}\in\Sigma} V_{\mathfrak{m}}\longrightarrow V
\]
is an equivalence and, for every $\mathfrak m\in\Sigma$, there exists a positive integer $N$ such that
$\mathfrak m^{N}$ acts by zero on $V_{\mathfrak m}$, i.e. $V_{\frm}$ is in the image $D_{\square}(\calD^{\la}(G,E)/\frm^{N})\to D_{\square}(\calD^{\la}(G,E))$ via the forgetful functor.
\end{definition}
\begin{remark}
	One can show that $\calD^{\la}(G,E)/\frm^N\simeq \calD^{\la}(G,E)\otimes^{L\square}_{\frz}\frz/\frm^N$ using Kostant's theorem that $U(\frg)$ is free over $\frz$ and the arguments for \cite[Prop.~4.8]{SchmidtStrauch2016dimensions}.
\end{remark}
Let $\rho$ be the half-sum of the positive roots. For a weight $\lambda\in(E\otimes_{\bbQ_p}\Lie T)^*$, write $|\lambda|$ for the orbit of $\lambda$ under the dot action $w\mathbin{\cdot}\lambda=w(\lambda+\rho)-\rho$. Let $\chi_\lambda:\frz\to E$ be the infinitesimal character of the Verma module of highest weight $\lambda$.  Put $\mathfrak{m}_\lambda=\ker(\chi_\lambda)$ which depends only on $|\lambda|$. Write $\pr_{|\lambda|}=\pr_{\mathfrak{m}_\lambda}$.

\begin{definition}\label{def:infinitesimal-character-block}
We say that an object $V\in D_{\square}(\calD^{\la}(G,E))$ lies in the $|\lambda|$-block if the natural map $\pr_{|\lambda|}(V)\to V$ is an equivalence.
\end{definition}

\begin{lemma}\label{lem:infinitesimal-character-block}
Suppose that $V\in D_{\square}(\calD^{\la}(G,E)/\frm_{\lambda}^{N})$ for some $\lambda$ and $N\geq 1$. Then $V$ lies in the $|\lambda|$-block.
\end{lemma}

\begin{proof}
By the Harish--Chandra isomorphism, $\frz$ is Noetherian. Since $\mathfrak{m}_\lambda^N$ acts by zero on $V$, all its cohomology modules are $\mathfrak{m}_\lambda$-torsion. Thus the natural map $\pr_{|\lambda|}(V)\to V$ is an equivalence by the arguments for \cite[\href{https://stacks.math.columbia.edu/tag/0954}{Tag 0954} \& \href{https://stacks.math.columbia.edu/tag/0955}{Tag 0955}]{stacks-project}.
\end{proof}

Let $\lambda$ and $\mu$ be compatible, meaning that $\nu:=\mu-\lambda$ lies in the image of $X^*(\mathrm{Res}_{L/\bbQ_p}(\bbT)_E)\to(E\otimes_{\bbQ_p}\Lie T)^*, \eta\mapsto d\eta$; see \cite[Defs.~2.4.2 \& 2.4.3]{jena2021translation}.  If $\bar\nu\in W.\nu$ is the dominant weight in the Weyl orbit of $\nu$, let $L(\bar\nu)$ be the corresponding algebraic representation of $\mathrm{Res}_{L/\bbQ_p}(\bbG)_E$ and of $G$.  
\begin{lemma}
	Let $V\in D_{\square}(\calD^{\la}(G,E))$ be $\frz$-finite. Then $L(\bar\nu)\otimes_{E_{\square}}^{L}V$ is also $\frz$-finite.
\end{lemma}
\begin{proof}
	We can assume $V=V_{\frm}$ for some $\frm$. Then $V\in D_{\square}(\calD^{\la}(G,E)/\frm^N)$ for some $N\geq 1$, and $L(\bar\nu)\otimes_{E_{\square}}^{L}V=(L(\bar\nu)\otimes_{E_{\square}}^{L}U(\frg)/\frm^N)\otimes^{L\square}_{U(\frg)/\frm^N}V$. Here we view $L(\bar\nu)\otimes_{E_{\square}}^{L}U(\frg)/\frm^N$ as a $U(\frg)$-bimodule as in \cite[\S~2.2]{bernstein1980tensor}. The desired statement now follows from the proof of \cite[Cor.~2.6]{bernstein1980tensor}.
\end{proof}
As in \cite[Def.~2.4.5]{jena2021translation}, for every $\frz$-finite solid $\calD^{\la}(G,E)$-module $V$ put
\[
 \mathrm T_\lambda^\mu(V):=
 \pr_{|\mu|}(
 L(\bar\nu)\otimes^L_{E_\square}\pr_{|\lambda|}(V)).
\]

Let $w_0$ be the longest Weyl element and put $\lambda^\vee=-w_0\lambda$.  If $\lambda$ is dominant algebraic, then $\lambda^\vee$ is the highest weight of the dual representation $L(\lambda)^\vee$.  Let $\iota$ be the antipode of $U(\mathfrak g)$, namely the anti-automorphism determined by $\iota(X)=-X$ for $X\in\mathfrak g$.  Then $\chi_{\lambda^\vee}=\chi_\lambda\circ\iota$ (see for example \cite[Prop.~7.4.9]{Dixmier1977enveloping}).

\begin{lemma}\label{lem:central-actions-analytic-kernel}
For the left and right $\calD^{\la}(G,E)$-actions $(g_1\cdot f\cdot g_2)(g)=f(g_1^{-1}gg_2^{-1})$ on $\calC_c^{\la}(G,E)$, the left and right actions of $\frz$ agree.
\end{lemma}

\begin{proof}
First, we have \cite[II.4.20, Remark]{Jantzen2003representations}
\[
 \calO(\mathrm{Res}_{L/\bbQ_p}(\bbG)_E)
\simeq\bigoplus_\eta L(\eta)^\vee\otimes_E L(\eta),
\]
where $\eta$ ranges over the dominant weights in $\mathrm{Res}_{L/\bbQ_p}(\bbT)_E$. The left action is on $L(\eta)^\vee$, while the right action on $L(\eta)$ is $v\cdot g_2=g_2^{-1}v$. Thus both actions of $z\in\frz$ on this summand are given by multiplying $\chi_{\eta^\vee}(z)=\chi_\eta(\iota(z))$.

The subspace of locally algebraic functions of $\calC_c^{\la}(G,E)$, spanned by functions of the form $\varphi f$, with $\varphi\in \calC_c^{\sm}(G,E)$ and $f\in \calO(\mathrm{Res}_{L/\bbQ_p}(\bbG)_E)$, is dense in $\calC_c^{\la}(G,E)$ for its locally convex topology. By continuity we get the claim.
\end{proof}

\begin{proposition}\label{prop:translation-duality}
For a $\frz$-finite $V\in\Rep^{\la}_{E_\square}(G)$ in the $|\lambda|$-block, the objects $V^\vee=R\uHom_E(V,E)$, $\bbD_{\ST}(V)$, and $\bbD_{\BZ}(V)$ lie in the $|\lambda^\vee|$-block, and there are natural equivalences in $D_{\square}(\calD^{\la}(G,E))$
\[
\begin{aligned}
 (\mathrm T_\lambda^\mu(V))^\vee
 &\simeq
 \mathrm T_{\lambda^\vee}^{\mu^\vee}(V^\vee),\\
 \bbD_{\ST}(\mathrm T_\lambda^\mu(V))
 &\simeq
 \mathrm T_{\lambda^\vee}^{\mu^\vee}(\bbD_{\ST}(V)),\\
 \bbD_{\BZ}(\mathrm T_\lambda^\mu(V))
 &\simeq
 \mathrm T_{\lambda^\vee}^{\mu^\vee}(\bbD_{\BZ}(V)).
\end{aligned}
\]
\end{proposition}

\begin{proof}
We have $(L(\bar\nu)\otimes_E V)^\vee \simeq L(\bar\nu)^\vee\otimes_E V^\vee$. Also for any weight $\lambda$, we have $(\pr_{|\lambda|}V)^\vee\simeq \pr_{|\lambda^\vee|}(V^\vee)$ under the $\frz$-finite assumption. Together with $L(\bar\nu)^\vee\simeq L(-w_0\bar\nu)$, this proves the first equivalence in the statement.

By definition, $\bbD_{\ST}(W)=(W^\vee)^{R\la}$. Applying $(-)^{R\la}$ to the first equivalence and using \cite[Cor.~3.2.14(3)]{JC23} gives the second equivalence.

For Bernstein--Zelevinsky duality, let $W\in\Rep^{\la}_{E_\square}(G)$ be $\frz$-finite. We may assume $W\in D_{\square}(\calD^{\la}(G,E)/\frm_{\lambda}^{N}),N\geq 1$. Then we have 
\begin{align*}
	R\uHom_{\calD^{\la}(G,E)}(W,\calC_c^{\la}(G,E))&\simeq R\uHom_{\calD^{\la}(G,E)/\frm_{\lambda}^N}(W,R\uHom_{\calD^{\la}(G,E)}(\calD^{\la}(G,E)/\frm_{\lambda}^N,\calC_c^{\la}(G,E)))\\
	&\simeq R\uHom_{\calD^{\la}(G,E)/\frm_{\lambda}^N}(W,R\uHom_{\frz}(\frz/\frm_{\lambda}^N,\calC_c^{\la}(G,E))).
\end{align*}
Lemma~\ref{lem:central-actions-analytic-kernel} shows that $R\uHom_{\frz}(\frz/\frm_{\lambda}^N,\calC_c^{\la}(G,E))$ is a $\frz/\frm_{\lambda}^N$-module for the right action of $G$. Converting this right action to a left action via the involution $\iota$, using Lemma~\ref{lem:infinitesimal-character-block}, we see $\bbD_{\BZ}(W)$ is $\frz$-finite and 
\begin{equation}\label{eq:BZ-primary-projection}
 \bbD_{\BZ}(\pr_{|\lambda|}(W))
 \simeq\pr_{|\lambda^\vee|}(\bbD_{\BZ}(W)).
\end{equation}
Finally, we have
\begin{equation}\label{eq:BZ-algebraic-tensor}
\begin{aligned}
 \bbD_{\BZ}(L(\bar\nu)\otimes^L_{E_\square}W)&\simeq
 R\uHom_{\calD^{\la}(G,E)}
 (W,L(\bar\nu)^\vee\otimes_{E_\square}\calC_c^{\la}(G,E))\\
 &\simeq L(\bar\nu)^\vee\otimes^L_{E_\square}\bbD_{\BZ}(W).
\end{aligned}
\end{equation}
The second equivalence is induced by sending $f\in L(\bar\nu)^\vee\otimes_{E_\square}\calC_c^{\la}(G,E)\simeq \calC_c^{\la}(G,L(\bar\nu)^\vee)$ to $g\mapsto g^{-1}f(g)$. We get the desired equivalence by \eqref{eq:BZ-primary-projection}.
\end{proof}
Let $Z_G$ be the center of $G$.  If $\zeta_{\bar\nu}:Z_G\to E^\times$ is the central character of $L(\bar\nu)$, then $\mathrm T_\lambda^\mu$ sends representations with central character $\omega$ to representations with central character $\omega\zeta_{\bar\nu}$. Let $\omega:Z_G\to E^\times$ be a locally analytic character.  Put
\[
 \calD^{\la}_\omega(G,E):=
 \calD^{\la}(G,E)
 \otimes^L_{\calD^{\la}(Z_G,E),\omega}E.
\]
Since the quotient $G\to G/Z_G$ admits locally analytic sections, one directly checks that the object $\calD^{\la}_\omega(G,E)$ is concentrated in degree $0$. Write $\Rep^{\la}_{E_\square,\omega}(G)$ for the full subcategory of $D_{\square}(\calD^{\la}_\omega(G,E))$ consisting of objects which are derived locally analytic. Let $\calC_{c,\omega}^{\la}(G,E)$ be the space of locally analytic functions on $G$ with compact support modulo $Z_G$ satisfying $f(zg)=\omega(z)^{-1}f(g)$ for $z\in Z_G$ and $g\in G$. The action $(g_1\cdot f\cdot g_2)(g)=f(g_1^{-1}gg_2^{-1})$ makes it a $\calD^{\la}_\omega(G,E)$-bimodule.

\begin{definition}
The Schneider--Teitelbaum duality with central character $\omega$ is
\[
\begin{aligned}
 \bbD_{\ST,\omega}: (\Rep^{\la}_{E_\square,\omega}(G))^{\mathrm{op}} &\longrightarrow \Rep^{\la}_{E_\square,\omega^{-1}}(G),\\
 V&\longmapsto(V^\vee)^{R\la}.
\end{aligned}
\]
The Bernstein--Zelevinsky duality with central character $\omega$ is
\[
\begin{aligned}
 \bbD_{\BZ,\omega}:
 (\Rep^{\la}_{E_\square,\omega}(G))^{\mathrm{op}}
 &\longrightarrow
 D_{\square}(\calD^{\la}_{\omega^{-1}}(G,E)),\\
 V&\longmapsto
 R\uHom_{\calD^{\la}_\omega(G,E)}
 (V,\calC_{c,\omega}^{\la}(G,E)).
\end{aligned}
\]
\end{definition}

\begin{proposition}\label{prop:fixed-central-translation}
Let $\lambda$ and $\mu$ be compatible weights, $\nu=\mu-\lambda$, and let $V\in\Rep^{\la}_{E_\square,\omega}(G)$ be $\frz$-finite and in the $|\lambda|$-block. Then, in $D_{\square}(\calD^{\la}_{\omega^{-1}\zeta_{\bar\nu}^{-1}}(G,E))$, there are natural equivalences, for $?\in\{\ST,\BZ\}$,
\[
 \bbD_{?,\omega\zeta_{\bar\nu}}
 (\mathrm T_\lambda^\mu(V))\simeq
 \mathrm T_{\lambda^\vee}^{\mu^\vee}
 (\bbD_{?,\omega}(V)).
\]
\end{proposition}

\begin{proof}
The equivalence for the Schneider--Teitelbaum duality follows from the proof of Proposition~\ref{prop:translation-duality}. For Bernstein--Zelevinsky duality, similarly as above, we have
\begin{equation}\label{eq:fixed-BZ-algebraic-tensor}
\begin{aligned}
 \bbD_{\BZ,\omega\zeta_{\bar\nu}}
 (L(\bar\nu)\otimes^L_{E_\square}V)&\simeq
 R\uHom_{\calD^{\la}_\omega(G,E)}
 (V,L(\bar\nu)^\vee\otimes_{E_\square}
 \calC_{c,\omega\zeta_{\bar\nu}}^{\la}(G,E))\\
 &\simeq L(\bar\nu)^\vee\otimes^L_{E_\square}\bbD_{\BZ,\omega}(V).
\end{aligned}
\end{equation} 
Then we can conclude as in the proof of Proposition~\ref{prop:translation-duality}.
\end{proof}

\setcounter{tocdepth}{2}
\bibliographystyle{alpha}
\bibliography{BZdualDrinfeld}
\end{document}